\documentclass[10pt]{article}

\newtheorem{theorem}{Theorem}[section]
\newtheorem{lemma}[theorem]{Lemma}
\newtheorem{proposition}[theorem]{Proposition}
\newtheorem{corollary}[theorem]{Corollary}

\newtheorem{definition}[theorem]{Definition}

\newtheorem{remark}[theorem]{Remark}

\newenvironment{poof}{\textit{Proof:  }}{
~\hfill\rule{2mm}{3mm}\vspace{.2in}}
 
\newenvironment{poofof}[1]{\textit{Proof of #1:  }}{
~\hfill\rule{2mm}{3mm}\vspace{.2in}}

\newcommand{\pa}[1]{\partial}

\usepackage{times}
\usepackage{amsfonts}
\usepackage{amssymb,amsmath,tabularx,graphicx,float,enumitem}
\usepackage{graphicx}
\usepackage{bm}
\usepackage{marginnote}
\usepackage{color}
\usepackage{stmaryrd}

\def\ds{\displaystyle}

\newcommand{\abs}[1]{\left\lvert{#1}\right\rvert}
\newcommand{\norm}[1]{\left\lVert{#1}\right\rVert}
\newcommand{\snorm}[1]{\left[{#1}\right]}
\newcommand{\hnorm}[3]{\norm{#1}_{C^{#2, #3}}}
\newcommand{\chnorm}[2]{\norm{#1}_{C^{#2}}}
\newcommand{\hdotnorm}[3]{\norm{#1}_{C^{#2, #3}}}
\newcommand{\hflnorm}[2]{\norm{#1}_{C^{\floor{#2}, #2 - \floor{#2}}}}

\newcommand{\hsnorm}[2]{\snorm{#1}_{C^{0, #2}}}

\newcommand{\starnorm}[1]{\abs{#1}_{*}}

\newcommand{\Hil}[1]{\mathcal{H}({#1})}
\newcommand{\mc}[1]{\mathcal{#1}}
\newcommand{\p}{\partial}
\newcommand{\D}[2]{\frac{d{#1}}{d{#2}}}
\newcommand{\PD}[2]{\frac{\p{#1}}{\p{#2}}}

\newcommand{\paren}[1]{\left({#1}\right)}
\newcommand\floor[1]{\lfloor#1\rfloor}
\newcommand{\diff}{\triangle}
\newcommand{\trl}{\mc{T}}
\newcommand{\wt}[1]{\widetilde{#1}}
\newcommand{\wh}[1]{\widehat{#1}}
\newcommand{\mbs}{\mathbb{S}^1}
\newcommand{\jump}[1]{\llbracket{#1}\rrbracket}
\newcommand{\dual}[2]{\langle{#1},{#2}\rangle}

\begin{document}

\title{Well-posedness and global behavior of the Peskin problem \\ of an immersed elastic filament in Stokes flow}
\author{Yoichiro Mori, Analise Rodenberg, Daniel Spirn
\footnote{The authors thank Peter Pol\'a\v{c}ik whose lectures on evolution equations, attended by the first two authors, led to the use of semigroup theory in H\"older spaces for this paper. The authors thank Charlie Peskin for allowing us to name this problem after him, and the first author thanks him for introducing the immersed boundary method to the first author as a graduate student; the discussions then led to this project. The first author was supported in part by NSF grant DMS-1516978 and DMS-1620316, and the third author by  NSF grant DMS-1516565. The authors also thank the hospitality of the IMA where most of this work was performed.}\\ \textit{\small School of Mathematics, University of Minnesota, Minneapolis, MN 55455}}
\date{\today}

\maketitle
\begin{abstract}
We consider the problem of a one dimensional elastic filament immersed in a two dimensional steady Stokes fluid.  Immersed boundary problems in which a thin elastic structure interacts with a surrounding fluid are prevalent in science and engineering, a class of problems for which Peskin has made pioneering contributions. Using boundary integrals, we first reduce the fluid equations to an evolution equation solely for the immersed filament configuration. We then establish local well-posedness for this equation with initial data in low-regularity H\"older spaces.  This is accomplished by first extracting the principal linear evolution by a small scale decomposition and then establishing precise smoothing estimates on the nonlinear remainder.  Higher regularity of these solutions is established via commutator estimates with error terms generated by an explicit class of integral kernels.  Furthermore, we show that the set of equilibria consists of uniformly parametrized circles and prove nonlinear stability of these equilibria with explicit exponential decay estimates, the optimality of which we verify numerically.  Finally, we identify a quantity which respects the symmetries of the problem and controls global-in-time behavior of the system.
\end{abstract}

\section{Introduction}\label{sect:intro}
Fluid structure interaction (FSI) problems, in which an elastic structure interacts with a surrounding fluid, abound in science and engineering. 
In this paper, we consider the problem of an elastic filament immersed in a two dimensional Stokes fluid. 
This problem, inspired by Peskin's immersed boundary (IB) method \cite{mittal2005immersed,1peskinIBM,2peskinIBM}, 
is arguably one of the simplest of FSI problems. 
We shall call our problem the {\em Peskin problem} in honor of his seminal contributions. 
The aim of our paper is to study the well-posedness of the Peskin problem 
and to address issues of regularity, stability of equilibria and questions on global behavior of solutions. 
An analytical study of  this problem is important not only because of its simplicity among FSI problems,
but also because such a study may form the basis for the numerical analysis of various FSI algorithms, 
including the IB method.

\begin{figure}
\begin{center}
\includegraphics[width=2.5in]{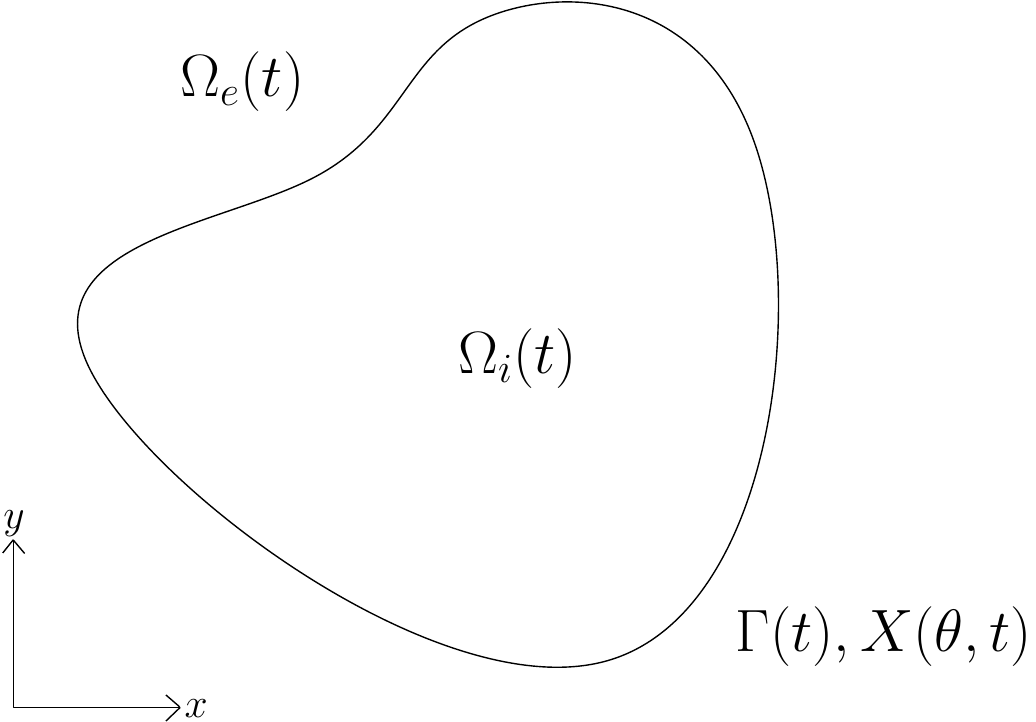}
\end{center}
\caption{\label{fig:config}Configuration of the problem.}
\end{figure}

\begin{figure}
\begin{center}
\includegraphics[width=\textwidth]{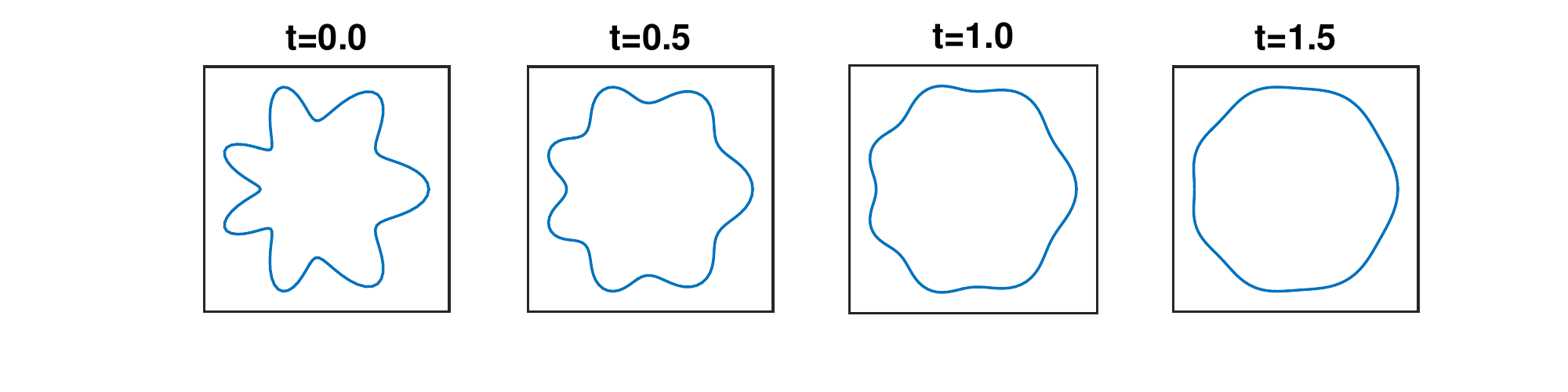}
\end{center}
\caption{\label{fig:demo}A sample simulation of the Peskin problem. Snapshots at $t=0, 0.5, 1$ and $1.5$ are shown.
The initial data is $\bm{X}_0(\theta)=((1+\cos(7\theta)/4)\cos(\theta)+\cos(2\theta)/8,(1+\cos(7\theta)/4)\sin(\theta)+\sin(2\theta)/8)^{\rm T}$.
Simulation was performed using the numerical algorithm discussed in Section \ref{sect:numerics}.}
\end{figure}

Let $\Gamma$ be a time-dependent simple closed curve in $\mathbb{R}^2$
representing the moving immersed elastic filament. The curve $\Gamma$ separates $\mathbb{R}^2$ into the interior domain $\Omega_{\rm i}$ and the exterior 
unbounded domain $\Omega_{\rm e}=\mathbb{R}^2\backslash(\Omega_{\rm i} \cup \Gamma)$. 
A schematic diagram is given in Figure \ref{fig:config}.
The curve $\Gamma$ is parameterized by $\theta\in \mbs$ (here and elsewhere we write $\mbs$ to mean $\mathbb{R}/2\pi\mathbb{Z} $) so that $\bm{X}(\theta,t)$ traces the curve $\Gamma$ (in the counter-clockwise direction) 
for each fixed $t$. The parametrization $\theta$ is taken to be the material or Lagrangian coordinate, so that $\bm{X}(\theta,t)$ for fixed $\theta$ moves with the 
local fluid velocity. The equations satisfied by $\bm{X}(\theta,t)$, the fluid velocity $\bm{u}(\bm{x},t), \bm{x}\in \mathbb{R}^2\backslash \Gamma$ 
and the pressure $p(\bm{x},t)$ are:
\begin{align}
\Delta \bm{u}-\nabla p&=0 \text{ in } \mathbb{R}^2\backslash \Gamma \label{stokes_eqn},\\
\nabla \cdot \bm{u}&=0 \text{ in } \mathbb{R}^2\backslash \Gamma \label{incomp},\\
\jump{\bm{u}}&=0 \text{ on } \Gamma \label{ucont},\\
\jump{(\nabla \bm{u}+(\nabla \bm{u})^{\rm T}-pI)\bm{n}}&=\p^2_\theta \bm{X} \abs{\p_\theta \bm{X}}^{-1} \text{ on } \Gamma,
\label{stress_jump}\\
\p_t \bm{X}&=\bm{u}(\bm{X},t). \label{Xt=u}
\end{align}
Equation \eqref{stokes_eqn} and \eqref{incomp} states that the fluids in the interior and exterior domains satisfy 
the incompressible Stokes equation. 
The viscosities of the fluid in the two domains are equal, and we have normalized this value to be $1$.
The interfacial conditions for $\bm{u}$ and $p$ at the curve $\Gamma$ are given in \eqref{ucont} and \eqref{stress_jump}.
For any quantity of interest $w$ defined on $\Omega\backslash \Gamma$, $\jump{w}$ denotes the jump in its value 
across $\Gamma$:
\begin{equation}
\jump{w}=\left.w\right|_{\Gamma_{\rm i}}-\left.w\right|_{\Gamma_{\rm e}}
\end{equation}
where $\left.w\right|_{\Gamma_{\rm i}}$ denotes the limiting value of $w$ evaluated on $\Gamma$ from the $\Omega_{\rm i}$
side, and likewise for $\left.w\right|_{\Gamma_{\rm e}}$. Condition \eqref{ucont} states that 
the fluid velocity is continuous at the interface $\Gamma$. 
This continuity implies that we may in fact extend $\bm{u}(\bm{x},t)$ up to the interface $\Gamma$.
Equation \eqref{stress_jump} is the stress jump condition; 
$\nabla \bm{u}$ is the rate of deformation matrix, $(\nabla \bm{u})^{\rm T}$ its transpose, 
$I$ is the identity matrix, and $\bm{n}$ is the outward unit normal on $\Gamma$
(pointing from $\Omega_{\rm i}$ to $\Omega_{\rm e}$). The notation $\p_\theta \bm{X}$ denotes the derivative (and $\p_\theta^k \bm{X}$
its $k$-th derivative) of $\bm{X}$ with respect to $\theta$, and $\abs{\cdot}$ is the Euclidean length. 
We have taken the force density (per unit material coordinate) of the immersed elastic filament to be proportional to  
$\p_\theta^2 \bm{X}$ (and normalized this coefficient to be $1$). Although this is the elastic force density we consider in this paper, 
we shall later comment on the other common choices for the elastic force. The Jacobian factor $\abs{\p_\theta \bm{X}}$ is present for
conversion from Lagrangian to arclength coordinates.
Equation \eqref{Xt=u} states that the elastic structure moves with the local fluid velocity. Note that evaluation of $\bm{u}$ on $\Gamma$ ($\bm{x}=\bm{X}(\theta,t)$)
is made possible by the continuity of $\bm{u}$ across $\Gamma$ as expressed in \eqref{ucont}.
Finally, we must impose conditions on the behavior of $\bm{u}$ and $p$ at infinity:
\begin{equation}\label{upinfty}
\bm{u}\to 0 \text{ as } \abs{\bm{x}}\to \infty, \; p\text{ is  bounded}.
\end{equation}
This concludes the statement of the Peskin problem. 
A sample simulation is given in Figure \ref{fig:demo}.

The {\em immersed boundary (IB)} reformulation of the above is to write the fluid equations together with the interfacial 
conditions in distributional form. 
We replace \eqref{stokes_eqn}-\eqref{stress_jump} with:
\begin{equation}\label{IB_stokes}
-\Delta \bm{u}+\nabla p=\int_{\mbs} \p_{\theta}^2\bm{X}\delta(\bm{x}-\bm{X}(\theta,t))d\theta,\quad \nabla \cdot \bm{u}=0,
\end{equation}
where $\delta(\bm{x}-\bm{X})$ is the Dirac delta function located at the point $\bm{X}$. The above equations
are to be satisfied in a distributional sense in $\mathbb{R}^2$.
Note that the interfacial conditions are now expressed in the form of a distributional body force on 
the right hand side of the Stokes equation. All other equations remain the same.
We shall refer to the use of \eqref{stokes_eqn}-\eqref{stress_jump} as the {\em jump} formulation of the problem.
If the functions $\bm{u},p$ and $\bm{X}$ are sufficiently smooth, one can show that the two 
formulations are equivalent \cite{lai2001remark}.  The IB formulation serves as the basis of the {\em immersed boundary (IB)} method.
The fluid domain and the immersed elastic structure are discretized independently of each other, 
and communication between the two takes place through \eqref{Xt=u} and the body force term in \eqref{IB_stokes}.
Ease of implementation and robustness of the algorithm have enabled the simulation of challenging FSI problems 
and have made the IB method among the most 
popular numerical methods for FSI problems. We refer the reader to the review articles \cite{mittal2005immersed,2peskinIBM} for details.

Finally, we consider the {\em boundary integral (BI)} formulation of the Peskin problem. 
Equation \eqref{IB_stokes} (or equivalently \eqref{stokes_eqn}-\eqref{stress_jump}), together with 
condition \eqref{upinfty} can be used to solve for $\bm{u}$ and $p$ to yield:
\begin{align}\label{bndry_integral_u}
\bm{u}(\bm{x},t)&=\int_{\mbs} G(\bm{x}-\bm{X}(\theta',t))\p_{\theta'}^2\bm{X}(\theta',t)d\theta',\\
\nonumber G(\bm{x})&=\frac{1}{4\pi}\paren{-\log \abs{\bm{x}} I+\frac{\bm{x}\otimes \bm{x}}{\abs{\bm{x}}^2}}\\
&=\frac{1}{4\pi}\paren{-\log\abs{\bm{x}}I+\frac{1}{\abs{\bm{x}}^2}\begin{pmatrix}x^2 & xy\\ xy & y^2\end{pmatrix}}, \; 
\bm{x}=(x,y)^{\rm T},\\
\label{bndry_integral_p}
p(\bm{x},t)&=\int_{\mbs} \bm{P}_{\rm st}(\bm{x}-\bm{X}(\theta',t))\cdot \p_{\theta'}^2\bm{X}(\theta',t)d\theta',\quad 
\bm{P}_{\rm st}(\bm{x})=\frac{\bm{x}}{2\pi \abs{\bm{x}}^2}.
\end{align}
Function $G$ is the Stokeslet, the fundamental solution of the Stokes equation in $\mathbb{R}^2$.
We note here that we do not suffer from the Stokes paradox of logarithmic growth of the velocity field $\bm{u}$ at infinity; 
this is thanks to the fact that the integral of $\bm{F}=\p_\theta^2 \bm{X}$ over $\theta$ is equal to $0$.
Substituting the above into \eqref{Xt=u}, we obtain the following closed equation for the evolution of $\bm{X}$.
\begin{equation}\label{BIform}
\p_t \bm{X}(\theta,t)=\int_{\mbs} G(\bm{X}-\bm{X}')\p_{\theta'}^2\bm{X}'d\theta'
=\int_{\mbs} G(\bm{X}(\theta,t)-\bm{X}(\theta',t))\p_{\theta'}^2\bm{X}(\theta',t)d\theta'.
\end{equation}
In the above and henceforth we write $\bm{X}'=\bm{X}(\theta',t)$, and we use similar notation for other primed quantities.
The BI formulation makes clear that the only initial condition that needs to be supplied to this problem is the 
initial configuration $\bm{X}(\theta,0)=\bm{X}_0(\theta)$.

The three formulations of the Peskin problem, the jump, IB and BI formulations, are equivalent assuming sufficient 
smoothness of the solutions $\bm{u}, p$ and $\bm{X}$. 
Certain analytical information is easier or more difficult to obtain depending on the specific formulation of the problem.
Furthermore, all of the three formulations are the basis of computational 
methods for this problem. The jump formulation is used in the immersed interface method \cite{li2006immersed, li2001immersed} and moving mesh 
methods such as the Arbitrary Lagrangian Eulerian (ALE) method \cite{donea1982arbitrary}, the IB formulation in the immersed boundary and related methods
\cite{hou2012numerical,mittal2005immersed,2peskinIBM,sethian2003level,tryggvason2001front} 
and the BI formulation can be used as a starting point for a boundary element/collocation method \cite{hou2012numerical,Poz}.
Establishing sufficient smoothness of the solution, therefore, is important from both analytic and numerical points of view.
From a numerical standpoint, unless some smoothness is established, 
it may not be clear whether the various methods are approximating the same solution. 
The wealth of numerical methods that can be used to tackle this problem 
has the potential to make the Peskin problem a standard testbed for the numerical analysis of FSI problems \cite{beale2015uniform,MoriIBMConv}.

We now state some important properties of  solutions to the Peskin problem. First of all, 
we have area conservation of the region $\Omega_{\rm i}$ which follows from the incompressibility condition \eqref{incomp}
and condition \eqref{Xt=u}. More concretely, 
\begin{equation}\label{area_conservation}
\D{}{t}\abs{\Omega_{\rm i}}=0, \; \abs{\Omega_{\rm i}}=\frac{1}{2}\int_{\mbs} \paren{X\p_\theta Y-Y\p_\theta X}d\theta, \; \bm{X}=(X,Y)^{\rm T}.
\end{equation}

 Solutions to the Peskin problem also satisfy the following energy identity, 
which states that the elastic energy $\mc{E}$ of the filament is lost through viscous dissipation $\mc{D}$ in the fluid:
\begin{equation}\label{energy_identity}
\D{\mc{E}}{t}=-\mc{D}, \; \mc{E}=\frac{1}{2}\int_{\mbs}\abs{\p_\theta \bm{X}}^2 d\theta, \; \mc{D}=\int_{\mathbb{R}^2} \abs{\nabla \bm{u}}^2 d\bm{x},
\end{equation}
which may be most easily seen from \eqref{IB_stokes}.  One can multiply the Stokes equation by $\bm{u}$,  integrate by parts, 
and use \eqref{Xt=u}. 
It is easy to justify both \eqref{area_conservation} and \eqref{energy_identity} given sufficient smoothness 
of the solution.

A third property of  solutions  is dilation invariance.
\begin{equation}\label{dilation_invariance}
\text{If } \bm{X}(\theta,t) \text{ is a solution, so is } a\bm{X}(\theta,t) \text{ for any } a>0;
\end{equation}
this can be seen most easily from \eqref{BIform}. 
We are not aware of any previous work that mentions this property.
Together with translation and rotation invariance ($\bm{X}$ rotated and translated is also a solution), 
we have a four dimensional group of symmetries
acting on the space of solutions.

To construct our solution theory, we shall use the BI formulation of the problem given in \eqref{BIform}.
First, note that we may formally integrate by parts in \eqref{BIform} to obtain:
\begin{equation}\label{pvBI}
\p_t \bm{X}(\theta,t)=-{\rm p.v.}\int_{\mbs} \p_{\theta'}G(\bm{X}-\bm{X}')\p_{\theta'}\bm{X}'d\theta',
\end{equation}
where the above integral is to be understood in the principal value sense. We may compute the kernel above as follows,
\begin{equation}\label{pthetaG}
-\p_{\theta'}G(\bm{X}-\bm{X}')=-\frac{1}{4\pi}\paren{\frac{\Delta \bm{X}\cdot \p_{\theta'} \bm{X}'}{\abs{\Delta \bm{X}}^2} I+\p_{\theta'}\paren{\frac{\Delta \bm{X}\otimes \Delta \bm{X}}{\abs{\Delta \bm{X}}^2}}}, \; \Delta \bm{X}=\bm{X}-\bm{X}'.
\end{equation}
The necessity to interpret \eqref{pvBI} in the principal value sense comes from the fact that:
\begin{equation*}
\frac{\Delta \bm{X}\cdot \p_{\theta'} \bm{X}'}{\abs{\Delta \bm{X}}^2}\approx\frac{1}{\theta-\theta'} \text{ when } \abs{\theta-\theta'}\ll 1.
\end{equation*}
This in turn suggests that we may extract the Hilbert transform from the kernel. 
Recall that, for a function $w$ defined on $\mbs$, its Hilbert transform $\mc{H}w$ is given by:
\begin{equation*}
(\mc{H}w)(\theta)=\frac{1}{2\pi}{\rm p.v.}\int_{\mbs} \cot\paren{\frac{\theta-\theta'}{2}}w(\theta')d\theta'.
\end{equation*}
We may thus rewrite \eqref{BIform} as follows:
\begin{equation}\label{SSD}
\begin{split}
\p_t \bm{X}&=\Lambda \bm{X}+\mc{R}(\bm{X}), \; \Lambda \bm{X}=-\frac{1}{4}\mc{H}(\p_\theta \bm{X}),\\
\mc{R}(\bm{X})&=-\frac{1}{4\pi}\int_{\mbs}\paren{
\paren{\frac{\Delta \bm{X}\cdot \p_{\theta'} \bm{X}'}{\abs{\Delta \bm{X}}^2}-\frac{1}{2}\cot\paren{\frac{\theta-\theta'}{2}}} I
+\p_{\theta'}\paren{\frac{\Delta \bm{X}\otimes \Delta \bm{X}}{\abs{\Delta \bm{X}}^2}}}\p_{\theta'}\bm{X}'d\theta'.
\end{split}
\end{equation}
The hope then is that our problem can be understood as a perturbation of the linear evolution driven by $\Lambda$. 
This approach of extracting the principal linear part in interfacial fluid problems is known as the {\em small scale decomposition}
($\Lambda$ should control behavior at higher spatial wave number, or behavior at small spatial scales, and hence its name) 
and was introduced in \cite{beale1993growth,hou1994removing} for the study of Hele-Shaw and water wave problems. 
In the context of numerical computation, the small scale decomposition allows for the removal 
of numerical stiffness; the stiff principal linear part is treated with an implicit numerical scheme whereas 
the remainder term is treated explicitly. Application of the small scale decomposition to IB problems can be 
found in \cite{hou2008efficient,hou2008removing}, although the small scale decomposition found in these papers 
seems to be slightly different from the one used in 
this manuscript, even taking into account the fact that they deal with the dynamic Stokes/Navier Stokes system.
In Section \ref{sect:numerics}, we shall use the small scale decomposition to develop a numerical 
scheme to computationally verify some of our theoretical results. The sample simulation in Figure \ref{fig:demo}
was generated using this algorithm.
 
Our approach to proving well-posedness is to turn \eqref{SSD} into an integral equation, a standard technique 
used in the study of semilinear parabolic equations \cite{henry2006geometric,lunardi,sell2013dynamics}:
\begin{equation}\label{mild_soln}
\bm{X}(t)=e^{t\Lambda} \bm{X}_0+\int_0^t e^{(t-s)\Lambda}\mc{R}(\bm{X}(s))ds.
\end{equation}
In the above, $e^{t\Lambda}$ is the semigroup generated by $\Lambda$, and $\bm{X}_0$ is the initial value.
We have suppressed the $\theta$ dependence of $\bm{X}$ and $\bm{X}_0$. 
The success of this approach hinges upon whether the extraction of $\Lambda$ is indeed a small 
scale decomposition in the sense that 
$\mc{R}$ is lower order (in some precise sense) compared with the principal linear part $\Lambda$.

We mention prior numerical and semi-analytical evidence suggesting that this may indeed be the case.
The operator $\Lambda$ can be viewed as the square-root of the Laplacian, or the Dirichlet-to-Neumann map, as we shall 
see in Section \ref{SGestSec}. As such, the effect of $\Lambda$ is to take one derivative in $\theta$. 
If $\Lambda$ is indeed the principal linear part, a numerical method based on explicit time-stepping should 
face a CFL type time-step restriction; to avoid numerical instability, 
$\Delta t$, the time step, should be refined proportionally to $\Delta \theta$,
the spatial discretization of the material coordinate. Studies on the stability of IB methods \cite{stockie1999analysis}, as well as our own numerical 
experiments,  strongly indicate the presence of such a CFL condition. 

From an analytical point of view, whether $\mc{R}$ is of lower order 
depends on the choice of function space. In this paper, we work in H\"older spaces $C^{k,\gamma}(\mbs)$, 
where $k$ is a non-negative integer and $0<\gamma<1$. We shall use the same notation $C^{k,\gamma}(\mbs)$
for functions with values in $\mathbb{R}$ as well as $\mathbb{R}^2$.
The standard definitions of H\"older norms are recalled in Section \ref{fnspace_defs}. 

Before we can state the definition of a solution to the Peskin problem,
we must introduce the following quantity defined for a function $\bm{X}(\theta)\in C^1(\mbs)$ which was introduced in the thesis work of \cite{bertozzi1991existence} and saw continued usage in \cite{majdabertozzi}:
\begin{equation}  \label{e:rRT}
\starnorm{\bm{X}}=\inf_{\theta,\theta' \in \mbs, \theta\neq \theta'} \frac{\abs{\bm{X}(\theta) - \bm{X}(\theta')}}{\abs{\theta-\theta'}}.
\end{equation}
Note that $\starnorm{\bm{X}}=0$ if and only if $\abs{\p_\theta \bm{X}}=0$ at some point or if the curve self-intersects, i.e. $\bm{X}(\theta)=\bm{X}(\theta')$ for some $\theta\neq \theta'$. 
Thus, $\bm{X}$ defines a non-degenerate simple closed curve if and only if $\starnorm{\bm{X}}>0$.
Let $C^n([0,T];C^{k,\gamma}(\mbs))$ be the space of $C^n$ functions of $t$, $0\leq t\leq T$ with values in $C^{k,\gamma}(\mbs)$. 
We define two notions of solutions to the Peskin problem.

\begin{definition}[Mild Solution]
Let $\bm{X}(t)\in C([0,T]; C^{1,\gamma}(\mbs)), 0<\gamma<1$ and $\starnorm{\bm{X}(t)}>0$ for $0\leq t\leq T$.
Then, $\bm{X}$ is a {\em mild solution} to the Peskin problem with initial value $\bm{X}(0)=\bm{X}_0$ 
if it satisfies equation \eqref{mild_soln} for $0<t\leq T$ and $\bm{X}(t) \to \bm{X}_0$ in $C^{1,\gamma}(\mbs)$ as $t\to 0$. 
\end{definition}
\begin{definition}[Strong Solution] Let $\bm{X}(t)\in C([0,T]; C^{1,\gamma}(\mbs))\cap C^1([0,T]; C^{0,\gamma}(\mbs)) , 0<\gamma<1$
and $\starnorm{\bm{X}(t)}>0$ for $0\leq t\leq T$.
Then, $\bm{X}$ is a {\em strong solution} to the Peskin problem with initial value $\bm{X}(0)=\bm{X}_0$ 
if it satisfies equation \eqref{pvBI} for $0<t\leq T$ and $\bm{X}(t) \to \bm{X}_0$ in $C^{1,\gamma}(\mbs)$ as $t\to 0$. 
\end{definition}

Let the little H\"older spaces $h^{1,\gamma}(\mbs)$ 
be the completion of the set of smooth functions in $C^{1,\gamma}(\mbs)$ 
(see the discussion before Proposition \ref{semigroup_strong_continuity}).  
Note that $C^{1,\alpha}(\mbs)\subset h^{1,\gamma}(\mbs)$ for any $\alpha>\gamma$.

We now state our result on the local well-posedness of the Peskin problem.
\begin{theorem}\label{LWPTheorem}
Consider the Peskin problem with initial value 
$\bm{X}_0 \in h^{1, \gamma}(\mbs)$, $0<\gamma<1$ with $\starnorm{\bm{X}_0}>0$.
Then, we have the following.
\begin{enumerate}[label=(\roman*)]
\item\label{LWPext} For some time $T>0$ depending on $\bm{X}_0$, there is a mild solution $\bm{X}(t)\in C([0,T]; C^{1,\gamma}(\mbs))$.
\item\label{LWPuniq} Suppose $\bm{X}(t)\in C([0,T]; C^{1,\gamma}(\mbs))$ 
is a mild solution to the Peskin problem. Then this solution is unique within the class $C([0,T];C^{1,\gamma}(\mbs))$.
\item\label{LWPcont} Let $\bm{X}(t)\in C([0,T];C^{1,\gamma}(\mbs))$ be the mild solution to the Peskin problem with initial data $\bm{X}_0$.
Then, there is an $\epsilon>0$ such that, for all initial data $\bm{Y}_0$ satisfying $\norm{\bm{X}_0-\bm{Y}_0}_{C^{1,\gamma}}\leq \epsilon$,
there is a mild solution $\bm{X}(t;\bm{Y}_0)\in C([0,T];C^{1,\gamma}(\mbs))$. Furthermore, $\bm{X}(t;\bm{Y}_0)$ is a continuous function of 
$\bm{Y}_0\in C^{1,\gamma}(\mbs)$ with values in $C([0,T];C^{1,\gamma}(\mbs))$.
\item\label{LWPmildstrong} The function $\bm{X}(t)$ is a mild solution on $[0,T]$ if and only if it is a strong solution on $[0,T]$. 
\end{enumerate}
\end{theorem}

We prove the existence of a mild solution \eqref{mild_soln} by 
a contraction mapping argument. 
There are two ingredients to the proof of Theorem~\ref{LWPTheorem}.  
The first ingredient is a set of estimates in H\"older norms of the semigroup operator generated by $\Lambda$ in \eqref{SSD}.
The semigroup $e^{t\Lambda}$ satisfies estimates typical of linear parabolic semigroups such as the heat propagator, 
except that $\Lambda$ has the effect of taking only one spatial derivative in contrast to the Laplacian which takes two spatial derivatives.
These estimates are found by an explicit representation of $e^{t\Lambda}$ as a convolution operator with the Poisson kernel, as discussed in Section~\ref{SGestSec}.

The second ingredient is a class of smoothing estimates on the nonlinear remainder $\mathcal{R}(\bm{X}(s))$; we show that $\mathcal{R}(\bm{X}): C^{1,\gamma} \mapsto C^{2\gamma} = C^{\floor{2\gamma}, 2\gamma - \floor{2\gamma}} $ for $\gamma \in (0,1) \backslash \frac{1}{2}$.  
This shows in essence that $\mc{R}$ has the effect of taking $1+\gamma-2\gamma=1-\gamma$ derivatives. As discussed earlier, $\Lambda$ behaves like taking 
one derivative, and $\mc{R}$ is thus genuinely lower order by $\gamma$ derivatives. This allows us to view $\Lambda$ as the principal part of the evolution,
making it possible to use  Duhamel's formula \eqref{mild_soln}. 
These crucial smoothing estimates on the remainder $\mc{R}$ arise from the structure of the kernel:   (a) the components of the kernel are composed of rational functions of finite differences of $\bm{X}$ and its derivatives and (b) the kernel is a perfect derivative in $\theta'$.   Our bounds on the components of the kernel, found in Section~\ref{s:calculus},  rely on careful, albeit elementary, estimates on these rational functions.    
Finally,  since the kernel is a perfect derivative, it allows us to gain an extra $\gamma$ in our H\"older estimate, which is used to close the argument. 
We remark here that our local existence theory is close to optimal, in the sense that $\mc{R}$ takes only $\gamma$ fewer derivatives  than  
$\Lambda$, 
and $\gamma>0$ can be made arbitrarily small. We are thus at the edge of applicability of semilinear parabolic techniques; any meaningful improvement on our local 
existence theory may require fundamentally different techniques.

Once we have proven the existence of the mild solution, we show that our mild solution has the expected $C^1([0,T]; C^{0,\gamma}(\mbs))$ 
regularity. 
Since the solution satisfies the differential form of the equation pointwise, we are able to conclude the existence of a unique strong solution.  

Our next result shows that the mild solution and its time derivative are arbitrarily smooth 
for any positive time.

\begin{theorem}\label{SmoothTheorem}
Consider the mild (strong) solution $\bm{X}$ of Theorem \ref{LWPTheorem}.
The function $\bm{X}$ is in $C^1([\epsilon,T]; C^n(\mbs))$ for any $n\in\mathbb{N}$ and $\epsilon>0$.
\end{theorem}

The proof of Theorem~\ref{SmoothTheorem} is found in Section~\ref{s:smoothness}.  Since the remainder $\mathcal{R}(\bm{X})$ is a nonlinear smoothing kernel acting on $\p_\theta \bm{X}$, in order to prove higher regularity, we introduce a class of integral kernels that allow us to move derivatives in $\theta$ on the nonlinear kernel into derivatives in  $\theta'$ acting on $\p_\theta \bm{X}$.  Since the error from this operation is lower order, the regularity improvement from the semigroup lets us gain arbitrarily high regularity in space.  The corresponding smoothness in time arises from equation \eqref{SSD}.   Higher regularity in time should be achievable using similar techniques, but we do not pursue it in this paper. 

An immediate corollary of this result is that the strong solution constructed in Theorem~\ref{LWPTheorem} is  classical in the sense 
that it satisfies the jump, IB and BI formulations of the equations  pointwise.  The precise definitions and these solutions are discussed 
in Section~\ref{sect:classical}.  
\begin{corollary}\label{equiv_formulations}
The notions of classical jump, IB, BI solutions and mild (strong) solution are equivalent.  
Furthermore, the mild (strong) solution satisfies area conservation \eqref{area_conservation} and energy identity \eqref{energy_identity}.
\end{corollary}
Any classical solution, which by definition should possess sufficient smoothness, is clearly a strong solution. This then proves 
the unique existence of classical solutions and the 
equivalence of the three formulations of the Peskin problem. 

In Section \ref{sect:equilibria_stability} we study the equilibria of the Peskin problem and their stability.
 The computation of the equilibria is performed using the jump formulation 
of the equations, which is made possible by Corollary~\ref{equiv_formulations}.
The only equilibria are circles in which the material points are evenly spaced:
\begin{equation}\label{circular_equilibria}
\begin{split}
\bm{X}(\theta)&= A\bm{e}_{\rm r}+B \bm{e}_{\rm t}+ C_1\bm{e}_x+ C_2 \bm{e}_y,\; A^2+B^2>0,\\
\bm{e}_{\rm r}&=\begin{pmatrix} \cos(\theta)\\ \sin(\theta) \end{pmatrix},\;
\bm{e}_{\rm t}=\begin{pmatrix} -\sin(\theta)\\ \cos(\theta) \end{pmatrix},\;
\bm{e}_x=\begin{pmatrix} 1\\ 0\end{pmatrix},\;
\bm{e}_y=\begin{pmatrix} 0\\ 1\end{pmatrix}.
\end{split}
\end{equation}
For later reference, we let $\wh{\mc{V}}$ denote the above set of circular equilibria
and let $\mc{V}$ be the linear space in $C^{1,\gamma}(\mbs)$ spanned by the above $4$ basis vectors $\bm{e}_{{\rm r},{\rm t},x,y}$. 

We now turn to the stability of these steady states. We first study the linearization of the evolution operator 
at the above uniformly parametrized circles. By dilation, translation and rotation invariance discussed above,
the linearized operator $\mc{L}$ is the same at every circle. This makes our analysis considerably simpler than would 
be otherwise and also leads to stronger results. In Section \ref{spect_linear_stability}, we explicitly compute the spectrum of $\mc{L}$
and obtain the decay properties of the semigroup $e^{t\mc{L}}$. The operator $\mc{L}$ has a four-dimensional 
kernel that coincides with $\mc{V}$. Except for the $0$ eigenvalue corresponding to the kernel $\mc{V}$, 
all eigenvalues are negative and real, and the leading non-zero principal eigenvalue is $-1/4$. 
In fact, $\mc{L}$ is a self-adjoint operator on $L^2(\mbs;\mathbb{R}^2)$, the space of square-integrable functions 
with values in $\mathbb{R}^2$. For two functions $\bm{v},\bm{w}\in L^2(\mbs;\mathbb{R}^2)$, 
we define the standard $L^2$ inner product as:
\begin{equation}\label{inner_product}
\dual{\bm{v}}{\bm{w}}=\int_{\mbs} \bm{v}(\theta)\cdot \bm{w}(\theta)d\theta.
\end{equation}

In Section \ref{nonlinear_stability}, we establish nonlinear stability of the circular equilibria. 
To state our result we introduce some notation.
Let $\mc{P}$ be the $L^2$ projection on to $\mc{V}$ and $\Pi$ its complementary projection:
\begin{equation}\label{projPPi}
\mc{P}\bm{X}=\frac{1}{2\pi}\sum_{\ell={\rm r},{\rm t},x,y}\dual{\bm{X}}{\bm{e}_\ell}\bm{e}_\ell, \; \Pi\bm{X}=\bm{X}-\mc{P}\bm{X}.
\end{equation}
The above $L^2$ projections are clearly well-defined operators on H\"older spaces as well. Notice that 
$\mc{P}\bm{X}\in \mc{V}$ is a circle so long as it does not degenerate to a point. 
Thus, the magnitude of $\Pi\bm{X}$ measures the distance from the set of circular equilibria.  
We let the norm on $\mc{V}$, which we denote by $\norm{\cdot}_{\mc{V}}$, to be the standard Euclidean $\mathbb{R}^4$ norm
with respect to the coordinate vectors $\bm{e}_{{\rm r},{\rm t},x,y}$.
We have the following result. 
\begin{theorem}\label{StabilityTheorem}
Circles with evenly spaced material points as given in \eqref{circular_equilibria} are the only equilibria of the Peskin system.
Furthermore, there is a constant $\rho_0>0$ that depends only on $\gamma$ with the following properties.
Consider a mild solution $\bm{X}(t)$ to the Peskin problem with initial data $\bm{X}_0\in h^{1,\gamma}(\mbs)$.
Let $R>0$ be the radius of $\mc{P}\bm{X}_0$, and suppose $\norm{\Pi\bm{X}_0}_{C^{1,\gamma}}\leq \rho_0R$.
Then, the solution to the Peskin problem is defined for all positive time and converges to a circle $\bm{Z}_\infty\in \wh{\mc{V}}$. 
Furthermore, we have the following estimates.
\begin{enumerate}[label=(\roman*)]
\item\label{i:1gamma_decay}For $t\geq 0$, we have:
\begin{align}
\label{difference_to_circle_decay}
\norm{\Pi\bm{X}(t)}_{C^{1,\gamma}}&\leq C\norm{\Pi\bm{X}_0}_{C^{1,\gamma}}e^{-t/4},\\
\label{projected_dynamics_decay}
\norm{\mc{P}\bm{X}(t)-\bm{Z}_\infty}_{\mc{V}}&\leq \frac{C}{R}\norm{\Pi\bm{X}_0}_{C^{1,\gamma}}^2e^{-t/2},
\end{align}
where the above constants $C$ depend only on $\gamma$. As an immediate consequence of the above results, we have:
\begin{equation}  \label{e:XdecaytoZ}
\norm{\bm{X}(t)-\bm{Z}_\infty}_{C^{1,\gamma}}\leq C\norm{\Pi\bm{X}_0}_{C^{1,\gamma}}e^{-t/4},
\end{equation}
where $C$ depends only on $\gamma$.
\item\label{i:higher_norm_decay} For any $n\in\mathbb{N}, n\geq 2$ and $t\geq \epsilon>0$, we have:
\begin{equation}\label{higher_norm_decay}
\norm{\Pi\bm{X}(t)}_{C^n}\leq C\norm{\Pi\bm{X}_0}_{C^{1,\gamma}}e^{-t/4},
\end{equation}
where the constant $C$ depends only on $n,\epsilon$ and $\gamma$.
An immediate consequence of this and \eqref{projected_dynamics_decay} is that, for 
any $n\in\mathbb{N}, n\geq 2$ and $t\geq \epsilon>0$,
\begin{equation}\label{e:XdecaytoZhigher}
\norm{\bm{X}(t)-\bm{Z}_\infty}_{C^n}\leq C\norm{\Pi\bm{X}_0}_{C^{1,\gamma}}e^{-t/4}
\end{equation}
where the above constant $C$ depends only on $n,\epsilon$ and $\gamma$.
\end{enumerate}
\end{theorem}
To prove this theorem, we first obtain a Lipschitz estimate on the derivative of the nonlinear remainder term.
We then use a standard Lyapunov-Perron type fixed point argument on time-exponentially weighted spaces 
to obtain the exponential decay to circular equilibria. Note here that, in all of the above estimates, the 
right and left hand side of the inequalities scale proportionally with dilation, 
as they should given dilation invariance of the Peskin system.

In many results of this type, it is only possible to prove that the decay rate can be made arbitrarily close but not equal
to the value of the real part of the leading non-zero eigenvalue (in our case, $-1/4$) \cite{lunardi,sell2013dynamics}. 
Here, an explicit calculation of the kernel $e^{t\mc{L}}$
allows us to obtain a sharp linear decay estimate, which in turn leads to this sharp result. 
Inequality \eqref{projected_dynamics_decay} indicates that 
the projected dynamics on the set of equilibria given by $\mc{P}\bm{X}(t)$ is 
exponentially approaching the limiting circle $\bm{Z}_\infty$ at twice the rate of $-1/2=2\times(-1/4)$. 
This somewhat unexpected result is a consequence of the fact that the zero-eigenspace $\mc{V}$ 
and the set of equilibria $\wh{\mc{V}}$ essentially coincide, which in turn is a reflection of 
the four-dimensional group of symmetries acting on the Peskin system. Finally, exponential 
decay in higher norms given in \eqref{higher_norm_decay} follows by combining 
\eqref{difference_to_circle_decay} with the parabolic regularity estimates of Section \ref{s:smoothness}.

In Section \ref{sect:numerics}, we computationally verify the exponential decay estimates stated in 
Theorem \ref{StabilityTheorem}. The numerical scheme we develop is a boundary integral method 
based on the small scale decomposition in \eqref{SSD}
and is second order accurate in time and spectrally accurate in space.
We see that the exponential decay rate of $\Pi\bm{X}(t)$ and $\mc{P}\bm{X}(t)$ is indeed asymptotically $-1/4$ and $-1/2$, respectively.

Finally, in Section \ref{sect:global}, we address issues of global behavior.
It is convenient to define the notion of a solution on half-open time intervals.
Let the space $C^n([0,T');C^{k,\gamma}(\mbs))$ to be the union of all $C^n([0,T];C^{k,\gamma}(\mbs))$ with $0<T<T'$. 
Here, $T'>0$ is allowed to be finite or $T'=\infty$.
\begin{definition}[Solution on half-open time intervals]
If $\bm{X}(t)\in C([0,T'); C^{1,\gamma}(\mbs)), 0<\gamma<1$ and $\starnorm{\bm{X}(t)}>0$ for $0\leq t<T'$, 
$\bm{X}(t)$ is a mild solution if the restriction of $\bm{X}(t)$ to any interval $[0,T], \; 0<T<T'$ is a mild solution. 
\end{definition}

Given initial data $\bm{X}_0\in h^{1,\gamma}(\mbs)$, 
define the maximal interval of existence of a mild solution $\tau_{\rm max}(\bm{X}_0)$ as follows.
Let $\mc{S}(\bm{X}_0)$ be the set of all $T>0$ such that there exists a mild solution $\bm{X}(t)\in C([0,T];C^{1,\gamma}(\mbs))$.
We let:
\begin{equation*}
\tau_{\rm max}(\bm{X}_0)=\sup_{T\in \mc{S}(\bm{X}_0)} T.
\end{equation*}
Note that, if $\bm{X}_1(t)\in C([0,T_1];C^{1,\gamma}(\mbs))$ and $\bm{X}_2(t)\in C([0,T_2];C^{1,\gamma}(\mbs))$ for $0<T_1<T_2$
with the same initial data $\bm{X}_0$, $\bm{X}_1(t)=\bm{X}_2(t)$ up to $t=T_1$ by the uniqueness result in Theorem \ref{LWPTheorem}.
Thus, one may speak of the unique mild solution $\bm{X}(t)$ with initial data $\bm{X}_0\in h^{1,\gamma}(\mbs)$ defined up to any $t<\tau_{\max}(\bm{X}_0)$.
Therefore, $\bm{X}(t)\in C([0,\tau_{\rm max}(\bm{X}_0)); C^{1,\gamma}(\mbs))$. It is important to note here that 
$\bm{X}(t)$ cannot be in $C([0,\tau_{\rm max}(\bm{X}_0)];C^{1,\gamma}(\mbs))$. If so, we will be able to extend the solution further 
by Theorem~\ref{LWPTheorem}, 
contradicting the definition of $\tau_{\rm max}$.
If $\tau_{\rm max}(\bm{X}_0)=\infty$, we say that the solution is global.

To state our results, we introduce the $\gamma$-{\em deformation ratio}:
\begin{equation*}
\varrho_\gamma(\bm{X}):=\frac{\hnorm{\p_{\theta}\bm{X}}{0}{\gamma}}{\starnorm{\bm{X}}}.
\end{equation*}
This quantity is invariant under translation, rotation and dilation.
Note that
\begin{equation*}
\varrho_\gamma(\bm{X})=\frac{\hnorm{\p_{\theta}\bm{X}}{0}{\gamma}}{\starnorm{\bm{X}}}
\geq \frac{\sup_{\theta \in \mbs} \abs{\p_\theta \bm{X}}}{\inf_{\theta \in \mbs} \abs{\p_\theta \bm{X}}}\geq 1.
\end{equation*}
The $\gamma$-deformation ratio is thus always greater than $1$, and we may replace the last inequality with an 
equality if $\bm{X}$ is a uniformly parametrized circle. In this sense, the $\gamma$-deformation ratio measures 
the degree to which $\bm{X}$ is deformed from a uniform circle configuration.

\begin{theorem}\label{GlobalBehaviorTheorem}
Given initial data $\bm{X}_0\in h^{1,\gamma}(\mbs)$, 
consider the mild solution $\bm{X}(t)\in C([0,\tau_{\rm max}(\bm{X}_0)); C^{1,\gamma}(\mbs))$.
\begin{enumerate}[label=(\roman*)]
\item\label{GBT1}Suppose $\tau_{\rm max}(\bm{X}_0)<\infty$. Then,
\begin{equation*}
\lim_{t\to \tau_{\rm max}(\bm{X}_0)}\varrho_\alpha(\bm{X}(t))=\infty,
\end{equation*}
for any $0<\alpha<1$.
In particular, the maximal existence time $\tau_{\rm max}(\bm{X}_0)$ does not depend on $\gamma$ 
(the space $C^{1,\gamma}$ in which the mild solution is considered).
\item\label{GBT2} Suppose the solution is global, that is $\tau_{\rm max}(\bm{X}_0)=\infty$. Suppose furthermore that
\begin{equation*}
\sup_{t\geq 0} \varrho_\alpha(\bm{X}(t))<\infty
\end{equation*}
for some $0<\alpha\leq \gamma$. Then, the solution $\bm{X}(t)$ converges exponentially to a uniformly parametrized circle 
as described in Theorem \ref{StabilityTheorem}.
\end{enumerate}
\end{theorem}
In the proof of this theorem,
the energy and area conservation identities \eqref{energy_identity} and \eqref{area_conservation} play a key role.
The deformation ratio bound together with area conservation gives a lower bound on 
$\starnorm{\bm{X}}$ whereas the deformation ratio bound and energy decay give an upper bound on the norm $\norm{\bm{X}}_{C^{1,\alpha}}$.

Item \ref{GBT1} above is a consequence of these bounds on $\starnorm{\bm{X}}$ and $\norm{\bm{X}}_{C^{1,\alpha}}$ as well as 
the regularity results of Section \ref{s:smoothness}. An interesting point about item \ref{GBT1} is that {\em all} deformation ratios $\varrho_\alpha(\bm{X}), 0<\alpha<1$
must tend to $\infty$ as $t$ reaches the maximal existence time. 
In particular, this shows that the maximal existence time is independent of the value of $\gamma$ in
$C^{1,\gamma}(\mbs)$, the space in which we consider the mild solution.  This leads us to conjecture that the $0$-deformation ratio, 
$\varrho_0(\bm{X}):=\frac{\norm{\p_{\theta}\bm{X}}_{C^0}}{\starnorm{\bm{X}}}$, would blow up at the finite extinction time.

Item \ref{GBT2} states that a global solution with bounded deformation ratio converges to a circle.
If the deformation ratio is bounded, $\starnorm{\bm{X}}$ and $\norm{\bm{X}}_{C^{1,\alpha}}$ are bounded by energy decay and area conservation as discussed above.
This shows that the orbit $\bm{X}(t)$ is relatively compact in any space $C^{1,\beta}(\mbs), 0<\beta<\alpha$, meaning that $\bm{X}(t)$
has a well-defined $\omega$-limit set in $C^{1,\beta}(\mbs)$. Viewing the energy as a Lyapunov function, one can then conclude that 
the $\omega$-limit set must consist only of  stationary circles. This, together with Theorem \ref{StabilityTheorem}, allows us to establish the desired result.


\subsection{Related Results}

A recent preprint \cite{LinTong} considers the Peskin problem and establishes local well-posedness of strong solutions with initial data $\bm{X}_0 \in H^{5/2}(\mbs)$ 
and $\starnorm{\bm{X}_0}>0$ which
generates a unique solution in $C([0,T];H^{5/2}(\mbs)) \cap H^1((0,T);H^2(\mbs))$ for some $T>0$.   
Local existence follows by energy arguments, use of Fourier multiplier methods, and an application of the Schauder fixed point theorem.  
The authors also show that a solution with initial data close to a circular equilibrium converges in the $H^{5/2}(\mbs)$ norm to a circular equilibrium 
at some exponential rate. This is established with the help of the energy identity.

The Peskin problem considered here is the simplest case of a much wider class of immersed boundary problems in which a thin elastic structure 
interacts with the surrounding fluid.
One extension of the Peskin problem would be to consider different constitutive relations for the elastic force. 
Instead of $\bm{F}=\p_\theta^2\bm{X}$ in \eqref{stress_jump}, 
we may consider the following more general elasticity law:
\begin{equation*}
\bm{F}=\p_\theta(\mc{T}(\abs{\p_\theta \bm{X}})\bm{\tau}), \; \bm{\tau}=\frac{\p_\theta \bm{X}}{\abs{\p_\theta \bm{X}}}.
\end{equation*}
Here, $\bm{\tau}$ is the unit tangent vector along the curve and $\mc{T}$ is the tension, which we assumed here to 
be a function of $\abs{\p_\theta \bm{X}}$ only. This leads to the energy identity:
\begin{equation}\label{energy_identity_general}
\D{\mc{E}}{t}=-\mc{D}, \; \mc{E}=\int_{\mbs} E(\abs{\p_\theta\bm{X}}) d\theta, \; \D{E}{s}=\mc{T}(s),
\end{equation}
where $\mc{D}$ is the same as in \eqref{energy_identity}. The Peskin problem discussed in this manuscript corresponds to the case $E(s)=s^2/2$ and $T(s)=s$.
Many choices of $E(s)$ are possible; one common choice is to set $E(s)=(s-\ell)^2/2$ where $2\pi\ell$ may be considered the natural (un-forced) length of the filament.
We may also replace the Stokes equation with the Navier-Stokes equation for the fluid equations in the interior and exterior domains, 
possibly with different viscosities and mass densities. We may also consider 3D problems in which the elastic force is generated
by a 2D membrane.  All of these generalizations are important in applications, and it would therefore be more descriptive to
refer to our problem as the 2D Peskin-Stokes problem with quadratic elastic energy.

We note that the choice $E(s)=s, \mc{T}(s)=1$ leads to the 2D surface tension problem. In this case, $\mc{E}$ is simply the length of the elastic filament, and
surface tension acts to decrease the interfacial length. This energy law makes it clear that the Peskin problem
and the surface tension problem are different. Surface tension only depends on the curvature of $\Gamma$, and therefore only on the shape (or geometry) of $\Gamma$.
In contrast, in the Peskin problem, the force depends on the material parametrization. In particular, stretching the interface leads to a force in the Peskin 
problem but not in the surface tension problem, and in this sense, the interface in the surface tension problem is {\em not} elastic. 
For example, any $\theta$ parametrization of the circle will be an equilibrium configuration for 
the surface tension problem, but only the uniform Lagrangian parametrization of the circle is an equilibrium configuration for the Peskin problem.
 
The surface tension problem itself has many variants. The analytical study of the one-phase problem, in which the fluid equations (Stokes or Navier Stokes) 
are satisfied in the interior region $\Omega_{\rm i}$ only, was initiated by Solonnikov \cite{solonnikov1977solvability}, and has since been taken up by many authors. 
The two-phase surface tension problem, in which the exterior region $\Omega_{\rm e}$ is also filled with a Stokes of Navier-Stokes fluid, 
possibly of different viscosity and mass density,
has also been studied by many authors, though the results are somewhat more recent. 
We refer the reader to \cite{pruss2009analysis,pruss2016moving,shimizu2012maximal} where an extensive
list of references on these problems can be found.
We also point to several recent results on problem with structures with more complicated energies interacting 
with the surrounding fluid \cite{ambrose2017well,cheng2007navier,cheng2010interaction,liu2017well,muha2013existence,plotnikov2011modelling,plotnikov2012strain}.

There are other problems in fluid mechanics which bear similarities to ours; the closest of which is the Muskat problem. In the  simplest setup in  two dimensions, the Muskat problem features two fluids in porous media whose dynamics are governed by Darcy's law. For nearly flat interfaces, the linearization of the Muskat problem has the same symbol as the Peskin problem considered here, and one expects similar local well-posedness and stability results so long as condition $\starnorm{\bm{X}}>0$ holds, see for example \cite{AmbroseMuskat04, CGBSH2WP, CCGS13, LWPMuskatCCG, EMStab,SCH04}. This condition is also popular in other interfacial fluid problems. It is widely used in addressing well-posedness and dynamics of vortex sheets and is often referred to as the arc-chord condition. See \cite{wu2003recent} for a discussion on vortex sheets separating the same fluid.

\section{Calculus Results}\label{s:calculus}

This section contains several estimates on rational functions of difference quotients of Holder continuous functions.  These estimates will be used throughout the 
paper and will greatly simplify much of the exposition. 

\subsection{Notation}\label{fnspace_defs}
We introduce some standard function spaces.
Let $C^k(\mathbb{S}^1), k=0,1,2,\ldots$ be the space of functions on $\mathbb{S}^1$ 
with $k$ continuous derivatives. Define the norms on these spaces in the usual way:
\begin{equation*}
\norm{u}_{C^k}=\sum_{l=0}^k \snorm{u}_{C^l}, \; \snorm{u}_{C^k}=\sup_{\theta\in \mbs} \abs{ \p_\theta^k u }. 
\end{equation*}
A function $u\in C^0(\mathbb{S}^1)$ is in the H\"older space $C^{0,\gamma}(\mathbb{S}^1), \; 0<\gamma<1$ if
\begin{equation*}
\sup_{\theta,\theta'\in \mathbb{S}^1} \frac{\abs{u(\theta)-u(\theta')}}{\abs{\theta-\theta'}^\gamma}<\infty.
\end{equation*}
Given the continuity of $u$, we may also restrict the range of $\theta$ and $\theta'$ to $\abs{\theta-\theta'}<1$, for instance.
Define the $C^{0,\gamma}$ norm as:
\begin{equation*}
\norm{u}_{C^{0,\gamma}}=\norm{u}_{C^0}+\hsnorm{u}{\gamma}, \; 
\hsnorm{u}{\gamma}=\sup_{\abs{\theta-\theta'}<1} \frac{\abs{u(\theta)-u(\theta')}}{\abs{\theta-\theta'}^\gamma}.
\end{equation*}
A function $u\in C^k(\mathbb{S}^1), k=0,1,2,\ldots$ is in $C^{k,\gamma}(\mathbb{S}^1), 0<\gamma<1$ if the $k$-th derivative of $u$
is in $C^{0,\gamma}(\mathbb{S}^1)$ and we define the norm on this space by:
\begin{equation*}
\norm{u}_{C^{k,\gamma}}=\norm{u}_{C^k}+\hsnorm{ \p_\theta^k u }. 
\end{equation*}

For any function $F(\theta)$ on the circle $\mathbb{S}^1$, we define

\begin{align*}
\Delta F=F(\theta)-F(\theta').
\end{align*} 
We also let $\p_{\theta}F$ be the derivative of $F$ evaluated at $\theta$ where as $\p_{\theta'}F'$
will be the derivative evaluated at $\theta'$. We will use the same notation for vector-valued functions on the circle.
We will be considering the difference quotient of functions evaluated at $\theta$ and $\theta+h$. 
Without loss of generality, assume that $0\leq \theta-h/2<\theta<\theta+h<\theta+3h/2\leq 2\pi$. 
This can be achieved since all of our functions are periodic. 
We will often split $\mbs$ into two parts, 
\begin{equation} \label{e:Idefn}
\begin{split}
\mc{I}_s & := (\theta - h/2, \theta + 3 h/2) \\
\mc{I}_f & := \mbs \backslash \mc{I}_s
\end{split}
\end{equation}
In the following, we drop the dependence on $\theta$ in the definition of $\mc{I}_s$.

For a function $f(\theta,\theta')$, we use the notation:
\begin{equation}\label{diff_trl}
\begin{split}
(\trl_h f)(\theta,\theta')&=f(\theta+h,\theta'),\\
(\diff_h f)(\theta,\theta')&=(\trl_h f-f)(\theta,\theta')=f(\theta+h,\theta')-f(\theta,\theta').
\end{split}
\end{equation}
A product rule follows from these definitions, 
\begin{equation}  \label{e:proddiffquot}
\diff_h (f (\theta, \theta') g(\theta, \theta')) = (\diff_h f)(\theta, \theta') (\trl_h g)(\theta, \theta') + f (\theta, \theta')(\diff_h g)(\theta, \theta').
\end{equation}

\subsection{Estimates on physical space multipliers}

We begin with some straightforward estimates on rational functions of a particular form.

\begin{lemma}\label{dqests}
Suppose the functions $\bm{Z}(\theta)=(Z(\theta),W(\theta))$ and $V(\theta)$ 
belong to $C^{1,\gamma}(\mathbb{S}^1)$. Assume also that $\starnorm{\bm{Z}}>0$.
Let
\begin{equation*}
\phi(\theta,\theta')=\frac{\Delta V}{\abs{\Delta \bm{Z}}}.
\end{equation*} 
\begin{enumerate}[label=(\roman *)]
\item \label{phi'ests}
The following estimates hold for $\phi$ and its derivatives. 
\begin{align}
\label{phiest}
\abs{\phi}&\leq \frac{\hdotnorm{V}{1}{\gamma}}{\starnorm{\bm{Z}}},\\
\label{phi1est}
\abs{ \p_\theta \phi }, \abs{ \p_{\theta'} \phi }&\leq C\frac{\hdotnorm{V}{1}{\gamma}\hdotnorm{\bm{Z}}{1}{\gamma}}{\starnorm{\bm{Z}}^2}\abs{\theta-\theta'}^{\gamma-1},\\
\label{phi2est}
\abs{\p_{\theta'} \p_\theta \phi }&\leq C
\frac{\hdotnorm{V}{1}{\gamma}\hdotnorm{\bm{Z}}{1}{\gamma}^2}{\starnorm{\bm{Z}}^3}\abs{\theta-\theta'}^{\gamma-2}.
\end{align}
If, in addition, $V=Z$ or $V=W$, we have the following estimates:
\begin{align}
\label{phiest'}
\abs{\phi}&\leq 1,\\
\label{phi1est'}
\abs{\p_\theta \phi }, \abs{\p_{\theta'} \phi }&\leq C\frac{\hdotnorm{\bm{Z}}{1}{\gamma}}
{\starnorm{\bm{Z}}}\abs{\theta-\theta'}^{\gamma-1},\\
\label{phi2est'}
\abs{ \p_{\theta'} \p_\theta \phi }&\leq C
\frac{\hdotnorm{\bm{Z}}{1}{\gamma}^2}
{\starnorm{\bm{Z}}^2}\abs{\theta-\theta'}^{\gamma-2}.
\end{align}
\item \label{phi'diffests}
Suppose $0<h<\abs{\theta-\theta'+h/2}$ and $0<\theta+h<2\pi$. Then, the following estimates hold.
\begin{align}\label{diffdelphi}
\abs{\diff_h(\p_\theta\phi)}&\leq C \frac{\hdotnorm{V}{1}{\gamma}\hdotnorm{\bm{Z}}{1}{\gamma}^2}
{\starnorm{\bm{Z}}^3}\paren{h^\gamma\abs{\theta-\theta'}^{-1}+h\abs{\theta-\theta'}^{\gamma-2}}.\\
\label{diffdel2phi}
\abs{\diff_h(\p_\theta'\p_\theta\phi)}&\leq C \frac{\hdotnorm{V}{1}{\gamma}\hdotnorm{\bm{Z}}{1}{\gamma}^2}
{\starnorm{\bm{Z}}^3}\paren{h^\gamma\abs{\theta-\theta'}^{-2}+h\abs{\theta-\theta'}^{\gamma-3}}.
\end{align}
If in addition, $V=Z$ or $V=W$, we have:
\begin{align}\label{diffdelphi'}
\abs{\diff_h(\p_\theta\phi)}&\leq C \frac{\hdotnorm{\bm{Z}}{1}{\gamma}^2}
{\starnorm{\bm{Z}}^2}\paren{h^\gamma\abs{\theta-\theta'}^{-1}+h\abs{\theta-\theta'}^{\gamma-2}}.\\
\label{diffdel2phi'}
\abs{\diff_h(\p_\theta'\p_\theta\phi)}&\leq C \frac{\hdotnorm{\bm{Z}}{1}{\gamma}^2}
{\starnorm{\bm{Z}}^2}\paren{h^\gamma\abs{\theta-\theta'}^{-2}+h\abs{\theta-\theta'}^{\gamma-3}}.
\end{align}
In the above, the positive constant $C$ do not depend on $V, \bm{Z}, \theta,\theta'$ or $h$.
\end{enumerate}
\end{lemma}

\begin{poof}
Let us first prepare some elementary estimates. First, we have:
\begin{equation}\label{Vdiff}
\abs{\Delta V}\leq \snorm{V}_{C^1}\abs{\theta-\theta'}\leq \hdotnorm{V}{1}{\gamma}\abs{\theta-\theta'}. 
\end{equation}
A similar estimate holds for $\bm{Z}$. We have, by definition of $\starnorm{\bm{Z}}$, 
\begin{equation}\label{Zstar}
\abs{\Delta \bm{Z}}\geq \starnorm{\bm{Z}}\abs{\theta-\theta'}.
\end{equation}
We also have:
\begin{equation}\label{Vdiffderiv}
\begin{split}
\abs{\frac{\Delta V}{\theta -\theta'}-\p_\theta V}&=\abs{\frac{1}{\theta-\theta'}\int_0^1 \D{}{s}V(s\theta+(1-s)\theta')ds-\p_\theta V }\\
&\leq \int_0^1 \abs{\p_\theta V (\theta'+s(\theta-\theta'))-\p_\theta V(\theta)}ds\\
&\leq \snorm{V}_{C^{1,\gamma}}\int_0^1 |1 - s|^{\gamma}\abs{\theta-\theta'}^\gamma ds
\leq \hdotnorm{V}{1}{\gamma}\abs{\theta-\theta'}^\gamma.
\end{split}
\end{equation}
A similar bound holds when $\p_{\theta} V $ is replaced by $\p_{\theta'} V' $ or for $V$ replaced by $\bm{Z}$, all with the same proof. 

We first consider the first three bounds in item \ref{phi'ests}.
The bound \eqref{phiest} follows from \eqref{Vdiff} and \eqref{Zstar}.  For \eqref{phi1est}, 
we prove the bound for $\partial_\theta \phi$. The bound for $\partial_{\theta'} \phi $ can be obtained
in exactly the same way. After some calculation, we obtain:
\begin{equation}\label{dphiAB}
\begin{split}
\p_\theta \phi &=\frac{1}{\abs{\Delta \bm{Z}}^3}\paren{ \p_\theta V \abs{\Delta \bm{Z}}^2-\Delta V \Delta \bm{Z}\cdot \p_\theta \bm{Z}}=A+B,\\
A&=\frac{1}{\abs{\Delta \bm{Z}}}\paren{ \p_\theta V -\frac{\Delta V}{\theta-\theta'}},\; B=\frac{1}{\abs{\Delta \bm{Z}}^3}\paren{\Delta V\Delta \bm{Z}\cdot\paren{\frac{\Delta \bm{Z}}{\theta-\theta'}- \p_\theta \bm{Z}}}. 
\end{split}
\end{equation}
Using \eqref{Vdiffderiv} and \eqref{Zstar}, we obtain:
\begin{equation}\label{Aestphi1}
\abs{A}\leq \frac{\hdotnorm{V}{1}{\gamma}}{\starnorm{\bm{Z}}}\abs{\theta-\theta'}^{\gamma-1}.
\end{equation}
We also have:
\begin{equation}\label{Bestphi1}
\abs{B}\leq \frac{\abs{\Delta V}}{\abs{\Delta \bm{Z}}^2}\abs{\frac{\Delta \bm{Z}}{\theta-\theta'}-  \p_\theta \bm{Z} }
\leq \frac{\hdotnorm{V}{1}{\gamma}\hdotnorm{\bm{Z}}{1}{\gamma}}{\starnorm{\bm{Z}}^2}\abs{\theta-\theta'}^{\gamma-1},
\end{equation}
where we used the Cauchy-Schwarz inequality in the first inequality and \eqref{Vdiff}, \eqref{Zstar} 
and \eqref{Vdiffderiv} (as applied to $\bm{Z}$) in the last inequality. Noting that $\hdotnorm{\bm{Z}}{1}{\gamma}\geq \starnorm{\bm{Z}}$, 
we may combine \eqref{Aestphi1} and \eqref{Bestphi1} to obtain \eqref{phi1est}.

Let us now prove bound \eqref{phi2est}. We have, after some calculation:
\begin{equation}\label{mixedphiDEF}
\begin{split}
\p_{\theta'} \p_\theta \phi &=D+E+F,\\
D&=\frac{1}{\abs{\Delta \bm{Z}}^5}\paren{ \p_\theta V \Delta \bm{Z}\cdot \p_{\theta'} \bm{Z}' \abs{\Delta \bm{Z}}^2
-\Delta V(\Delta \bm{Z}\cdot \p_{\theta} \bm{Z} )(\Delta \bm{Z}\cdot \p_{\theta'} \bm{Z}' )},\\
E&=\frac{1}{\abs{\Delta \bm{Z}}^5}\paren{\p_{\theta'} V' \Delta \bm{Z}\cdot \p_{\theta} \bm{Z} \abs{\Delta \bm{Z}}^2
-\Delta V(\Delta \bm{Z}\cdot \p_{\theta} \bm{Z} )(\Delta \bm{Z}\cdot \p_{\theta'} \bm{Z}' )},\\
F&=\frac{1}{\abs{\Delta \bm{Z}}^5}\paren{\Delta V\p_{\theta} \bm{Z} \cdot \p_{\theta'} \bm{Z}' \abs{\Delta \bm{Z}}^2
-\Delta V(\Delta \bm{Z}\cdot \p_{\theta} \bm{Z} )(\Delta \bm{Z}\cdot \p_{\theta'} \bm{Z}' )}.
\end{split}
\end{equation}
We estimate $D$. We have:
\begin{equation}\label{D1D2phiests}
\begin{split}
D&=D_1+D_2,\\
D_1&=\frac{1}{\abs{\Delta \bm{Z}}^3}\paren{ \p_\theta V-\frac{\Delta V}{\theta-\theta'}}\Delta \bm{Z}\cdot\p_{\theta'} \bm{Z}' ,\\
D_2&=\frac{1}{\abs{\Delta \bm{Z}}^5}\Delta V(\Delta \bm{Z}\cdot \p_{\theta'} \bm{Z}' )\paren{\Delta \bm{Z}
\cdot\paren{\frac{\Delta \bm{Z}}{\theta-\theta'}-\p_{\theta} \bm{Z} }}.
\end{split}
\end{equation}
Using estimates \eqref{Vdiff}, \eqref{Zstar} and \eqref{Vdiffderiv}, we have:
\begin{equation}\label{D1bounddqests}
\abs{D_1}\leq \frac{1}{\abs{\Delta \bm{Z}}^2}
\abs{\p_\theta V -\frac{\Delta V}{\theta-\theta'}}\abs{\p_{\theta'} \bm{Z}' }
\leq \frac{\hdotnorm{V}{1}{\gamma}\hdotnorm{\bm{Z}}{1}{\gamma}}
{\starnorm{\bm{Z}}^2}\abs{\theta-\theta'}^{\gamma-2}.
\end{equation}
Likewise, for $D_2$, we have:
\begin{equation*}
\abs{D_2}\leq \frac{\hdotnorm{V}{1}{\gamma}\hdotnorm{\bm{Z}}{1}{\gamma}^2}
{\starnorm{\bm{Z}}^3}\abs{\theta-\theta'}^{\gamma-2}.
\end{equation*} 
Combining the above estimates and noting that $\hdotnorm{\bm{Z}}{1}{\gamma}\geq \starnorm{\bm{Z}}$, we have:
\begin{equation*}
\abs{D}\leq C\frac{\hdotnorm{V}{1}{\gamma}\hdotnorm{\bm{Z}}{1}{\gamma}^2}
{\starnorm{\bm{Z}}^3}\abs{\theta-\theta'}^{\gamma-2}
\end{equation*}
In a similar fashion, $E$ and $F$ can be shown to satisfy the same bound. This concludes the proof of \eqref{phi2est}.

To obtain the alternate estimates in item \ref{phi'ests} when $V=Z$ or $V=W$ note that:
\begin{equation*}
\abs{\Delta V}\leq \abs{\Delta \bm{Z}} \text{ if } V=Z \text{ or } V=W.
\end{equation*}
Bound \eqref{phiest'} is a direct consequence of this. We may use this to improve the bound on $B$
in \eqref{Bestphi1} to obtain \eqref{phi1est'}. We may also show that $D_2$ satisfies the bound \eqref{D1bounddqests}
and hence obtain a better bound for $D$ and similarly for $E$ and $F$ of \eqref{mixedphiDEF}. This yields \eqref{phi2est'}.

Finally, we turn to item \ref{phi'diffests}. Let us first consider \eqref{diffdelphi}.
Given expression \eqref{dphiAB}, we may estimate $\diff_h A$ and $\diff_h B$ seperately.
We have:
\begin{equation*}
\begin{split}
\diff_h A&=A_1+A_2, \\
A_1&=\frac{1}{\abs{\Delta \bm{Z}}}\diff_h\paren{\p_\theta V -\frac{\Delta V}{\theta-\theta'}},\;
A_2=\diff_h\paren{\frac{1}{\abs{\Delta \bm{Z}}}}\trl_h\paren{\p_\theta V -\frac{\Delta V}{\theta-\theta'}}.
\end{split}
\end{equation*}
We first estimate $A_1$.
\begin{equation*}
\abs{\diff_h \p_{\theta} V }\leq \hdotnorm{V}{1}{\gamma}h^\gamma.
\end{equation*}
Furthermore, similarly to the calculation in \eqref{Vdiffderiv}, we have:
\begin{equation*}
\begin{split}
\abs{\diff_h \paren{\frac{\Delta V}{\theta-\theta'}}}&\leq \int_0^1\abs{\p_\theta V ((1-s)(\theta+h)+s\theta')- \p_\theta V ((1-s)\theta+s\theta')}ds\\
&\leq \snorm{V}_{C^{1,\gamma}}\int_0^1 (1-s)^\gamma h^\gamma ds\leq \hdotnorm{V}{1}{\gamma}h^\gamma.
\end{split}
\end{equation*}
Using the above two relations and \eqref{Zstar}, we thus have:
\begin{equation}\label{A1estinphi1}
\abs{A_1}\leq C\frac{\hdotnorm{V}{1}{\gamma}}{\starnorm{\bm{Z}}}h^\gamma\abs{\theta-\theta'}^{-1}
\end{equation}
We now estimate $A_2$.
\begin{equation*}
\diff_h\paren{\frac{1}{\abs{\bm{Z}}}}=
h\int_0^1 \p_\theta \paren{\frac{1}{\abs{\bm{Z}}}}((\theta-\theta')+sh)ds
\end{equation*}
Since 
\begin{equation*}
\abs{\p_\theta \paren{\frac{1}{\abs{\bm{Z}}}}}
=\frac{\abs{\Delta\bm{Z}\cdot \p_\theta \bm{Z} }}{\abs{\Delta \bm{Z}}^3}
\leq \frac{\hdotnorm{\bm{Z}}{1}{\gamma}}{\starnorm{\bm{Z}}^2}\abs{\theta-\theta'}^{-2}
\end{equation*}
we have:
\begin{equation*}
\diff_h\paren{\frac{1}{\abs{\bm{Z}}}}
\leq \frac{\hdotnorm{\bm{Z}}{1}{\gamma}}{\starnorm{\bm{Z}}^2}\int_0^1\abs{\theta+sh-\theta'}^{-2}ds.
\end{equation*}
Note that:
\begin{equation}\label{thetadiffbound}
\frac{\abs{\theta-\theta'}}{\abs{\theta-\theta'+sh}}\leq 1+\frac{sh}{\abs{\theta-\theta'+sh}}\leq 1+\frac{sh}{\abs{\theta-\theta'+h/2}-h/2}\leq 1+2s\leq 3.
\end{equation}
Thus,
\begin{equation*}
\diff_h\paren{\frac{1}{\abs{\bm{Z}}}} \leq C\frac{\hdotnorm{\bm{Z}}{1}{\gamma}}{\starnorm{\bm{Z}}^2}h\abs{\theta-\theta'}^{-2}.
\end{equation*}
Using the above and \eqref{Vdiffderiv}, we have:
\begin{equation*}
\abs{A_2}\leq 
C\frac{\hdotnorm{V}{1}{\gamma}\hdotnorm{\bm{Z}}{1}{\gamma}h}{\starnorm{\bm{Z}}^2}\abs{\theta-\theta'}^{-2}\abs{\theta+h-\theta'}^{\gamma}.
\end{equation*}
In much the same way as in \eqref{thetadiffbound},
\begin{equation}\label{thetatheta'h}
\frac{\abs{\theta-\theta'+h}}{\abs{\theta-\theta'}}\leq 3.
\end{equation}
Thus, 
\begin{equation}\label{A2estinphi1}
\abs{A_2}\leq 
C\frac{\hdotnorm{V}{1}{\gamma}\hdotnorm{\bm{Z}}{1}{\gamma}}{\starnorm{\bm{Z}}^2}h\abs{\theta-\theta'}^{\gamma-2}.
\end{equation}
Combining \eqref{A1estinphi1} and \eqref{A2estinphi1} and using $\starnorm{\bm{Z}}\leq \norm{\bm{Z}}_{C^{1,\gamma}}$ we see that
\begin{equation}\label{diffAest}
\abs{\diff_h A}\leq 
C\frac{\hdotnorm{V}{1}{\gamma}\hdotnorm{\bm{Z}}{1}{\gamma}}{\starnorm{\bm{Z}}^2}
\paren{h^\gamma\abs{\theta-\theta'}^{-2}+h\abs{\theta-\theta'}^{\gamma-2}}.
\end{equation}

We turn to $\diff_h B$. We have:
\begin{equation*}
\begin{split}
\diff_h B&=B_1+B_2, \\
B_1&=\frac{\Delta V\Delta \bm{Z}}{\abs{\Delta \bm{Z}}^3}\cdot\diff_h\paren{ \p_\theta \bm{Z} -\frac{\Delta \bm{Z}}{\theta-\theta'}},\;
B_2=\diff_h\paren{\frac{\Delta V\Delta \bm{Z}}{\abs{\Delta \bm{Z}}^3}}\cdot \trl_h\paren{  \p_\theta \bm{Z} -\frac{\Delta \bm{Z}}{\theta-\theta'}}.
\end{split}
\end{equation*}
In much the same way we obtained the estimates for $A_1$, for $B_1$, we have:
\begin{equation}\label{B1estinphi1}
\abs{B_1}\leq C\frac{\hdotnorm{V}{1}{\gamma}\hdotnorm{\bm{Z}}{1}{\gamma}}{\starnorm{\bm{Z}}^2}h^\gamma\abs{\theta-\theta'}^{-1}.
\end{equation}
We turn to $B_2$. We have:
\begin{equation}\label{derivB2fct}
\abs{\p_\theta \paren{\frac{\Delta V\Delta \bm{Z}}{\abs{\Delta \bm{Z}}^3}}}\leq C\frac{\hdotnorm{V}{1}{\gamma}\hdotnorm{\bm{Z}}{1}{\gamma}^2}{\starnorm{\bm{Z}}^3}
\abs{\theta-\theta'}^{-2},
\end{equation}
where we used $\starnorm{\bm{Z}}\leq \hdotnorm{\bm{Z}}{1}{\gamma}$. Using the same procedure as was used for$A_2$, we have:
\begin{equation*}
\abs{\diff_h\paren{\frac{\Delta V\Delta \bm{Z}}{\abs{\Delta \bm{Z}}^3}}}\leq C\frac{\hdotnorm{V}{1}{\gamma}\hdotnorm{\bm{Z}}{1}{\gamma}^2}{\starnorm{\bm{Z}}^3}
h\abs{\theta-\theta'}^{-2}.
\end{equation*}
Combining this with \eqref{Vdiffderiv} (as applied to $\bm{Z}$), we obtain:
\begin{equation}\label{B2estinphi1}
\abs{B_2}\leq C\frac{\hdotnorm{V}{1}{\gamma}\hdotnorm{\bm{Z}}{1}{\gamma}^2}{\starnorm{\bm{Z}}^3}
h\abs{\theta-\theta'}^{\gamma-2}.
\end{equation}
Combining \eqref{B1estinphi1} and \eqref{B2estinphi1}, we obtain the estimate:
\begin{equation}\label{diffBest}
\abs{\diff_h B}\leq C\frac{\hdotnorm{V}{1}{\gamma}\hdotnorm{\bm{Z}}{1}{\gamma}^2}{\starnorm{\bm{Z}}^3}
\paren{h^\gamma\abs{\theta-\theta'}^{-2}+h\abs{\theta-\theta'}^{\gamma-2}}.
\end{equation}
Combining \eqref{diffAest} and \eqref{diffBest}, we obtain \eqref{diffdelphi}. When $V=Z$ or $V=W$, we can obtain the following bound in place of \eqref{derivB2fct}:
\begin{equation*}
\abs{\p_\theta \paren{\frac{\Delta V\Delta \bm{Z}}{\abs{\Delta \bm{Z}}^3}}}\leq C\frac{\hdotnorm{\bm{Z}}{1}{\gamma}}{\starnorm{\bm{Z}}^2}\abs{\theta-\theta'}^{-2}.
\end{equation*}
This allows us to prove \eqref{diffdelphi'}. 

Finally, we turn to \eqref{diffdel2phi} and \eqref{diffdel2phi'}. By \eqref{mixedphiDEF}, we may estimate the differences of $D,E$ and $F$
separately. To estimate $\diff_h D$, we estimate $\diff_h D_1$ and $\diff_h D_2$ where $D_{1,2}$ are given in \eqref{D1D2phiests}.
The difference $\diff_h D_1$ can be estimated in a similar way to $\diff_h A$ above and $\diff_h D_2$ similarly to $\diff_h B$ above.
The differences $\diff_h E$ and $\diff_h F$ can be estimated similarly to $\diff_h D$. We omit the details.
\end{poof}

\begin{lemma}\label{GPL1}
Suppose we have a set of functions 
$\{\bm{Z}_i\}_{i = 0}^{n}$ with $\bm{Z}_i=(Z_i,W_i) \in C^{1, \gamma}$, $\gamma \in (0, 1)$  with $\starnorm{\bm{Z}_0} > 0$.
Let $g(\theta, \theta')$ be
\begin{align*}
g(\theta, \theta')&= \prod_{i=0}^n\paren{\frac{\Delta Z_i}{\abs{\Delta \bm{Z}_0}}}^{\alpha_i}\paren{\frac{\Delta W_i}{\abs{\Delta \bm{Z}_0}}}^{\beta_i}
\end{align*}
where $\alpha_i, \beta_i \in \mathbb{N}\cup \{0\}$ and let
\begin{align*}
N=\sum_{i = 1}^{n}(\alpha_i + \beta_i), \; N_0=\sum_{i=0}^n (\alpha_i+\beta_i)=(\alpha_0+\beta_0)+N. 
\end{align*}
We have the following estimates.
\begin{align} 
\label{gest0}
\abs{g(\theta,\theta')}&\leq  \frac{1}{\starnorm{\bm{Z}_0}^N}\prod_{i = 1}^{n}\hdotnorm{\bm{Z}_i}{1}{\gamma}^{\alpha_i + \beta_i},\\
\label{e:gestmult1}
\left| \p_\theta g(\theta, \theta')\right|, \left|\p_{\theta'}g(\theta, \theta')\right| &\leq C\frac{\hdotnorm{\bm{Z}_0}{1}{\gamma}}{\starnorm{\bm{Z}_0}^{N+1}}\prod_{i = 1}^{n}\hdotnorm{\bm{Z}_i}{1}{\gamma}^{\alpha_i + \beta_i}|\theta - \theta'|^{\gamma - 1},\\
\label{mixedderiv}
\left| \p_\theta \p_{\theta'} g(\theta, \theta')\right| &\leq 
C\frac{\hdotnorm{\bm{Z}_0}{1}{\gamma}^2}{\starnorm{\bm{Z}_0}^{N+2}}\prod_{i = 1}^{n}\hdotnorm{\bm{Z}_i}{1}{\gamma}^{\alpha_i + \beta_i}\abs{\theta-\theta'}^{\gamma-2}
\end{align}
for some constant $C$ which depends only on $\gamma$ and $N_0$.
\end{lemma}

\begin{poof}
Given the assumption on the indices $\alpha_i,\beta_i$ and $N$, we have:
\begin{equation}\label{g0exp}
g=\prod_{i=0}^n \phi_i^{\alpha_i}\psi_i^{\beta_i}, \; 
\phi_i=\frac{\Delta Z_i}{\abs{\Delta \bm{Z}_0}}, \; \psi_i=\frac{\Delta W_i}{\abs{\Delta \bm{Z}_0}}.
\end{equation}
The functions $\phi_i$ and $\psi_i$ satisfy the assumptions of Lemma \ref{dqests}. Thus, using \eqref{phiest} and \eqref{phi1est}, 
for $\phi_i, i\geq 1$, we have:
\begin{equation}\label{phii01est}
\abs{\phi_i}\leq \frac{\hdotnorm{\bm{Z}_i}{1}{\gamma}}{\starnorm{\bm{Z}_0}},
\quad \abs{\p_\theta \phi_i }\leq C\frac{\hdotnorm{\bm{Z}_i}{1}{\gamma}\hdotnorm{\bm{Z}_0}{1}{\gamma}}
{\starnorm{\bm{Z}_0}^2}\abs{\theta-\theta'}^{\gamma-1},
\end{equation}
where we used the fact that $\hdotnorm{\bm{Z}_i}{1}{\gamma}$ dominates the norm of its components
The function $\psi_i, i\geq 1$ satisfy exactly the same bounds. For $\phi_0$, from \eqref{phiest'} and \eqref{phi1est'}
we have the bounds:
\begin{equation}\label{phii01est'}
\abs{\phi_0}\leq 1,
\quad \abs{ \p_\theta \phi_i  }\leq C\frac{\hdotnorm{\bm{Z}_0}{1}{\gamma}}
{\starnorm{\bm{Z}_0}}\abs{\theta-\theta'}^{\gamma-1}.
\end{equation}
The same bound holds for $\psi_0$.
Inequality \eqref{gest0} is immediate from \eqref{g0exp}, \eqref{phii01est} and \eqref{phii01est'}.
We now prove \eqref{e:gestmult1}. 
Taking the derivative of $g$ with respect to $\theta$, we obtain:
\begin{equation}\label{gderiv}
\begin{split}
\p_\theta g &=\sum_{k=0}^n (\alpha_kA_k+\beta_kB_k),\\
A_k&=(\p_\theta \phi_k ) \phi_k^{\alpha_k-1}\psi_k^{\beta_k}\prod_{i\neq k} \phi_i^{\alpha_i}\psi_i^{\beta_i},\quad 
B_k=(\p_\theta \psi_k ) \phi_k^{\alpha_k}\psi_k^{\beta_k-1}\prod_{i\neq k} \phi_i^{\alpha_i}\psi_i^{\beta_i}.
\end{split}
\end{equation}
Let us bound $A_k, k\geq 1$. 
\begin{equation*}
\abs{A_k}\leq 
\abs{\p_\theta \phi_k }\abs{\phi_k}^{\alpha_k-1}\abs{\psi_k}^{\beta_k}\prod_{i\neq k} \abs{\phi_i}^{\alpha_i}\abs{\psi_i}^{\beta_i}
\leq C\frac{\hdotnorm{\bm{Z}_0}{1}{\gamma}}{\starnorm{\bm{Z}_0}^{N+1}}
\prod_{i = i}^{n}\hdotnorm{\bm{Z}_i}{1}{\gamma}^{\alpha_i + \beta_i}|\theta - \theta'|^{\gamma - 1},
\end{equation*}
where we used \eqref{phii01est} and \eqref{phii01est'}(and the same bounds for the $\psi_i$'s). 
The terms $B_k, k\geq 0$ as well as $A_0$ satisfy exactly the same inequality.
Inequality \eqref{e:gestmult1} is now immediate.
The inequality for $\partial_{\theta'}  g$ can be obtained in exactly the same way.

We now prove \eqref{mixedderiv}. Using the notation of \eqref{gderiv}, we have:
\begin{equation}\label{g2deriv}
\p_{\theta'} \p_\theta g =\sum_{k=0}^N \paren{\alpha_k \p_{\theta'} A_k +\beta_k \p_{\theta'} B_k }
\end{equation}
Let us estimate $\partial_{\theta'} A_k$ when $k\geq 1$ (and $k$ for which $\alpha_k\neq 0)$. We see that:
\begin{equation}\label{delA_k}
\p_{\theta'} A_k = (\p_{\theta'} \p_\theta \phi_k)  g_k+ \p_{\theta} \phi_k \p_{\theta'} g_k ,\quad 
g_k=\phi_k^{\alpha_k-1}\psi_k^{\beta_k}\prod_{i\neq k} \phi_i^{\alpha_i}\psi_i^{\beta_i}.
\end{equation}
Recall from Lemma \ref{dqests} that
\begin{equation*}
\abs{ \p_{\theta'} \p_\theta \phi_k  }\leq 
C\frac{\hdotnorm{\bm{Z}_k}{1}{\gamma}\hdotnorm{\bm{Z}_0}{1}{\gamma}^2}
{\starnorm{\bm{Z}_0}^3}\abs{\theta-\theta'}^{\gamma-2}.
\end{equation*}
Thus, combining \eqref{phii01est} and the above, we have:
\begin{equation*}
\abs{ (\p_{\theta'} \p_\theta \phi_k)g_k}\leq 
C\frac{\hdotnorm{\bm{Z}_0}{1}{\gamma}^2}{\starnorm{\bm{Z}_0}^{N+2}}\prod_{i = 1}^{n}\hdotnorm{\bm{Z}_i}{1}{\gamma}^{\alpha_i + \beta_i}|\theta - \theta'|^{\gamma - 2}
\end{equation*}
Note that $g_k$ satisfies the hypothesis of the present lemma for $g$ with $\alpha_k$ replaced by $\alpha_k-1$.
We may thus use \eqref{e:gestmult1} directly to estimate $\partial g_k/\partial \theta'$. Thus, using this together with \eqref{phii01est}, 
\begin{equation*}
\abs{\p_{\theta} \phi_k \p_{\theta'} g_k }
\leq C\frac{\hdotnorm{\bm{Z}_0}{1}{\gamma}^2}{\starnorm{\bm{Z}_0}^{N+2}}\prod_{i = 1}^{n}\hdotnorm{\bm{Z}_i}{1}{\gamma}^{\alpha_i + \beta_i}|\theta - \theta'|^{2\gamma - 2}
\end{equation*}
Note that 

\begin{align*}
|\theta - \theta'|^{2\gamma - 2} \leq C|\theta - \theta'|^{\gamma - 2}
\end{align*}
for some constant $C$ depending only on $\gamma$. Combining the above two estimates, we obtain a bound on $\partial A_k/\partial \theta'$. 
Noting that $A_0$ and $B_k, k\geq 0$ satisfy the same estimates, we obtain the desired result.
\end{poof}

\begin{lemma}\label{GPL2}
Let $g(\theta, \theta')$ as in Lemma \ref{GPL1}, $0<h<\abs{\theta-\theta'+h/2}$ and $0<\theta+h<2\pi$.
\begin{align} 
\label{e:gestmult2}
\abs{\diff_h \p_{\theta'} g(\theta,\theta')} &\leq C\frac{\hdotnorm{\bm{Z}_0}{1}{\gamma}^{2}}{\starnorm{\bm{Z}_0}^{N+2}}\prod_{i = 1}^{n}\hdotnorm{\bm{Z}_i}{1}{\gamma}^{\alpha_i + \beta_i}
h\abs{\theta-\theta'}^{\gamma-2}\\
\label{diffdel2g}
\abs{\diff_h\p_\theta \p_{\theta'}g(\theta,\theta')} &\leq C\frac{\hdotnorm{\bm{Z}_0}{1}{\gamma}^3}{\starnorm{\bm{Z}_0}^{N+3}}
\prod_{i = 1}^{n}\hdotnorm{\bm{Z}_i}{1}{\gamma}^{\alpha_i + \beta_i}\paren{h^\gamma\abs{\theta-\theta'}^{-2}+h\abs{\theta-\theta'}^{\gamma-3}},
\end{align}
where the constant $C$ depends only on $\gamma$.
\end{lemma}
\begin{poof}
Let us first prove \eqref{e:gestmult2}. 
We have:
\begin{equation*}
\begin{split}
\abs{\diff_h \p_{\theta'} g}&=h\int_0^1 \abs{\p_{\theta'} \p_\theta g(\theta+sh,\theta')}ds\\
&\leq C\frac{\hdotnorm{\bm{Z}_0}{1}{\gamma}^2}{\starnorm{\bm{Z}_0}^{N+2}}\prod_{i = 1}^{n}\hdotnorm{\bm{Z}_i}{1}{\gamma}^{\alpha_i + \beta_i}
Ch\int_0^1\paren{\abs{\theta+sh-\theta'}^{\gamma-2}+\abs{\theta+sh-\theta'}^{2\gamma-2}}ds\\
&\leq C\frac{\hdotnorm{\bm{Z}_0}{1}{\gamma}^2}{\starnorm{\bm{Z}_0}^{N+2}}\prod_{i = 1}^{n}\hdotnorm{\bm{Z}_i}{1}{\gamma}^{\alpha_i + \beta_i}
Ch\paren{\abs{\theta-\theta'}^{\gamma-2}+\abs{\theta-\theta'}^{2\gamma-2}}.
\end{split}
\end{equation*}

We next prove \eqref{diffdel2g}. Recall that $\p_\theta \p_{\theta'} g$ can be written as \eqref{g2deriv}. We use the notation there.
Consider bounding $\diff_h\p_{\theta' }A_k$. From \eqref{delA_k}, we have:
\begin{equation*}
\begin{split}
\diff_h \p_\theta' A_k&=A_{k1}+A_{k2}+A_{k3}+A_{k4},\\
A_{k1}&=(\p_{\theta'}\p_\theta \phi_k)\diff_h g_k, \; 
A_{k2}=\diff_h\paren{\p_{\theta'} \p_\theta \phi_k} \trl_h g_k,\\
A_{k3}&=(\p_\theta{\phi_k})\diff_h\p_{\theta'}g_k,\;
A_{k4}=\diff_h(\p_\theta{\phi_k})\trl_h\p_{\theta'} g_k.
\end{split}
\end{equation*}
Suppose $k\geq 1$.
Let us estimate $A_{k1}$ (for $k$ such that $\alpha_k\neq 0$). Note that 
\begin{equation}\label{diffhp'g}
\begin{split}
\abs{\diff_h g_k}&\leq h\int_0^1 \abs{\p_\theta g_k(\theta+sh,\theta')}ds\\
&\leq C\frac{\hdotnorm{\bm{Z}_0}{1}{\gamma}}{\starnorm{\bm{Z}_0}^{N}}\prod_{i = 1}^{n}\hdotnorm{\bm{Z}_i}{1}{\gamma}^{\wt{\alpha}_i + \beta_i}h\int_0^1\abs{\theta+sh- \theta'}^{\gamma - 1}ds\\
&\leq C\frac{\hdotnorm{\bm{Z}_0}{1}{\gamma}}{\starnorm{\bm{Z}_0}^{N}}\prod_{i = 1}^{n}\hdotnorm{\bm{Z}_i}{1}{\gamma}^{\wt{\alpha}_i + \beta_i}h\abs{\theta-\theta'}^{\gamma - 1},
\end{split}
\end{equation}
where $\wt{\alpha}_i=\alpha_i $ if $i\neq k$ and $\wt{\alpha}_k=\alpha_k-1$ otherwise.
Combining this with \eqref{mixedderiv}, we have:
\begin{equation*}
\abs{A_{k1}}\leq 
C\frac{\hdotnorm{\bm{Z}_0}{1}{\gamma}^3}{\starnorm{\bm{Z}_0}^{N+3}}\prod_{i = 1}^{n}\hdotnorm{\bm{Z}_i}{1}{\gamma}^{\alpha_i + \beta_i}
h\paren{\abs{\theta-\theta'}^{2\gamma-3}+\abs{\theta-\theta'}^{3\gamma-3}}.
\end{equation*}
Turn next to $A_{k2}$. Using \eqref{diffdel2phi} and \eqref{gest0}, we have:
\begin{equation*}
\abs{A_{k2}}\leq 
\abs{\diff_h \p_{\theta'} \p_\theta \phi_k}\abs{\trl_h g_k}\leq 
C\frac{\hdotnorm{\bm{Z}_0}{1}{\gamma}^2}
{\starnorm{\bm{Z}}^{N+2}}\prod_{i = 1}^{n}\hdotnorm{\bm{Z}_i}{1}{\gamma}^{\alpha_i + \beta_i}\paren{h^\gamma\abs{\theta-\theta'}^{-2}+h\abs{\theta-\theta'}^{\gamma-3}}
\end{equation*}
For $A_{k3}$, using \eqref{phi1est} and \eqref{diffhp'g} we have:
\begin{equation*}
\begin{split}
\abs{A_{k3}}&\leq \abs{\p_{\theta'}\phi_k}\abs{\diff_h\p_{\theta'} g_k}\\
&\leq
C\frac{\hdotnorm{\bm{Z}_0}{1}{\gamma}^3}{\starnorm{\bm{Z}_0}^{N+3}}\prod_{i = 1}^{n}\hdotnorm{\bm{Z}_i}{1}{\gamma}^{\alpha_i + \beta_i}
h\paren{\abs{\theta-\theta'}^{2\gamma-3}+\abs{\theta-\theta'}^{3\gamma-3}}
\end{split}
\end{equation*}
Let us turn to $A_{k4}$. By \eqref{e:gestmult1} we have the bound:
\begin{equation*}
\begin{split}
\abs{\trl_h\p_{\theta'} g_k}
&\leq 
C\frac{\hdotnorm{\bm{Z}_0}{1}{\gamma}}{\starnorm{\bm{Z}_0}^{N}}\prod_{i = 1}^{n}\hdotnorm{\bm{Z}_i}{1}{\gamma}^{\wt{\alpha}_i + \beta_i}|\theta+h - \theta'|^{\gamma - 1}\\
&\leq 
C\frac{\hdotnorm{\bm{Z}_0}{1}{\gamma}}{\starnorm{\bm{Z}_0}^{N}}\prod_{i = 1}^{n}\hdotnorm{\bm{Z}_i}{1}{\gamma}^{\wt{\alpha}_i + \beta_i}|\theta- \theta'|^{\gamma - 1}
\end{split}
\end{equation*}
where we used \eqref{thetatheta'h} in the last inequality.
Thus,
\begin{equation*}
\abs{A_{k4}}\leq \abs{\diff_h(\p_\theta{\phi_k})}\abs{\trl_h\p_{\theta'} g_k}
\leq C\frac{\hdotnorm{\bm{Z}_0}{1}{\gamma}^3}{\starnorm{\bm{Z}_0}^{N+3}}\prod_{i = 1}^{n}\hdotnorm{\bm{Z}_i}{1}{\gamma}^{\alpha_i + \beta_i}
\paren{h^\gamma |\theta- \theta'|^{\gamma - 2}+h\abs{\theta-\theta'}^{2\gamma-3}}.
\end{equation*}
We may now combine our estimates on $A_{k1},\ldots,A_{k4}$ to obtain the bound on $\diff_h \p_\theta' A_k$ noting that the bound is dominated by $h^{\gamma}|\theta - \theta'|^{-2}$ and $h|\theta - \theta'|^{\gamma - 3}$. The bounds on $\diff_h \p_\theta' A_0$ as well as $\diff_h \p_\theta' B_k, k=0,\ldots, n$ can be obtained in the same fashion.
This concludes the proof. 
\end{poof}

\begin{lemma}\label{dqestslog}
Suppose the functions $\bm{Z}(\theta)=(Z(\theta),W(\theta))$ and $V(\theta)$ 
belong to $C^{1,\gamma}(\mathbb{S}^1)$. Assume also that $\starnorm{\bm{Z}}>0$.
Let
\begin{equation*}
q(\theta,\theta')=\frac{ \p_{\theta'} V' -(\Delta V/(\theta-\theta'))}{\abs{\Delta \bm{Z}}}.
\end{equation*}
\begin{enumerate}[label=(\roman *)]
\item We have the following estimates:
\begin{align}\label{q0est}
\abs{q(\theta,\theta')}&\leq \frac{\hdotnorm{V}{1}{\gamma}}{\starnorm{\bm{Z}_0}}\abs{\theta-\theta'}^{\gamma-1}.\\
\label{q1est}
\abs{\p_\theta q(\theta,\theta')}&\leq C\frac{\hdotnorm{V}{1}{\gamma}\hdotnorm{\bm{Z}}{1}{\gamma}}{\starnorm{\bm{Z}}^2}\abs{\theta-\theta'}^{\gamma-2}.
\end{align}
\item Suppose $0<h<\abs{\theta-\theta'+h/2}$ and $0<\theta+h<2\pi$. Then, the following estimate holds.
\begin{equation}
\label{diffq1est}
\abs{\diff_h \p_\theta q(\theta,\theta')}\leq C\frac{\hdotnorm{V}{1}{\gamma}\hdotnorm{\bm{Z}}{1}{\gamma}^2}{\starnorm{\bm{Z}}^3}
\paren{h^\gamma\abs{\theta-\theta'}^{-2}+h\abs{\theta-\theta'}^{\gamma-3}}
\end{equation}
where the constant $C$ depends only on $\gamma$.
\end{enumerate}
\end{lemma}
\begin{poof}
We first prove \eqref{q0est}. 
By \eqref{Zstar} and \eqref{Vdiffderiv} (applied to $\p_{\theta'} V' $) we see that 
\begin{equation*}
\abs{q}\leq \frac{\abs{(\Delta V/(\theta-\theta'))-\p_{\theta'} V' }}{\abs{\Delta \bm{Z}/(\theta-\theta')}}\abs{\theta-\theta'}^{-1}
\leq \frac{\hdotnorm{V}{1}{\gamma}}{\starnorm{\bm{Z}}}\abs{\theta-\theta'}^{\gamma-1}.
\end{equation*}
We next turn to \eqref{q1est}.
\begin{equation}\label{pthetaqexp}
\p_\theta q=\paren{\frac{\Delta V}{\theta-\theta'}-\p_{\theta} V }\frac{1}{(\theta-\theta')\abs{\Delta \bm{Z}}}-\paren{\frac{\Delta V}{\theta-\theta'}-\p_{\theta'} V' }\frac{\Delta \bm{Z}\cdot \p_{\theta} \bm{Z} }{\abs{\Delta \bm{Z}}^3}
\end{equation}
Using \eqref{Zstar} and \eqref{Vdiffderiv}, we have:
\begin{equation*}
\abs{\p_\theta q}\leq C\frac{\hdotnorm{V}{1}{\gamma}\hdotnorm{\bm{Z}}{1}{\gamma}}{\starnorm{\bm{Z}}^2}\abs{\theta-\theta'}^{\gamma-2}.
\end{equation*}
To obtain estimates in $\diff_h \p_\theta q$, take the difference of the two terms in \eqref{pthetaqexp}. This can be obtained 
in a similar fashion to the calculation in Lemma \ref{dqests} to obtain \eqref{diffdelphi}. We omit the details.
\end{poof}

\begin{lemma}\label{GPlog1}
Suppose we have a set of functions 
$\{\bm{Z}_i\}_{i = 0}^{n}$ with $\bm{Z}_i=(Z_i,W_i) \in C^{1, \gamma}$, $\gamma \in (0, 1)$  with $\starnorm{\bm{Z}_0} > 0$.
\begin{equation}\label{loglemmaf}
f(\theta, \theta')= \frac{(\p_{\theta'} V' -\Delta V/(\theta-\theta'))}{\abs{\Delta \bm{Z}_0}}\prod_{i=0}^n\paren{\frac{\Delta Z_i}{\abs{\Delta \bm{Z}_0}}}^{\alpha_i}\paren{\frac{\Delta W_i}{\abs{\Delta \bm{Z}_0}}}^{\beta_i},
\end{equation}
where $V=Z_k$ or $V=W_k$ for some $k=0,\ldots n$. Let
\begin{align*}
N=\sum_{i = 1}^{n}(\alpha_i + \beta_i), \; N_0=\sum_{i=0}^n (\alpha_i+\beta_i)=(\alpha_0+\beta_0)+N. 
\end{align*}
Then, we have:
\begin{align}\label{logest1}
|f(\theta, \theta')| &\leq \frac{\hdotnorm{\bm{Z}_k}{1}{\gamma}}{\starnorm{\bm{Z}_0}^{N+1}}
\prod_{i=1}^n \hdotnorm{\bm{Z}_i}{1}{\gamma}^{\alpha_i + \beta_i}|\theta - \theta'|^{\gamma-1},\\
\label{logest2}
\left| \partial_{\theta}f(\theta, \theta')\right| &\leq C
\frac{\hdotnorm{\bm{Z}_k}{1}{\gamma}\hdotnorm{\bm{Z}_0}{1}{\gamma}}{\starnorm{\bm{Z}_0}^{N+2}}
\prod_{i=1}^n \hdotnorm{\bm{Z}_i}{1}{\gamma}^{\alpha_i + \beta_i} |\theta - \theta'|^{\gamma-2},
\end{align}
where $C$ depends only on $\gamma$.
\end{lemma}
\begin{poof}
We can write $g$ as:
\begin{equation}\label{fqg}
f=qg, \; q=\frac{(\p_{\theta'} V' -\Delta V/(\theta-\theta'))}{\abs{\Delta \bm{Z}_0}}.
\end{equation}
Note that $g$ is exactly of the same form as the function estimated in Lemma \eqref{GPL1}.
Inequality \eqref{logest1} is a direct consequence of the \eqref{q0est} of the previous lemma and of \eqref{gest0}. 
Inequality \eqref{logest1} is now immediate. We turn to \eqref{logest2}. Note that:
\begin{equation*}
\abs{\p_\theta g }\leq \abs{\p_\theta q }\abs{g}+\abs{q}\abs{\p_\theta g }.
\end{equation*}
Each of the factors may be estimated using \eqref{q1est}, \eqref{gest0}, \eqref{q0est} and \eqref{e:gestmult2}, 
from which we obtain the desired inequality.
\end{poof}

\begin{lemma}\label{GPlog2}
Let $f(\theta,\theta')$ be as in Lemma \ref{GPlog1} and suppose $0<h<\abs{\theta-\theta'+h/2}$ and $0<\theta+h<2\pi$.
\begin{align}\label{logest3}
\abs{\diff_h f(\theta,\theta')}&\leq 
C\frac{\hdotnorm{\bm{Z}_k}{1}{\gamma}\hdotnorm{\bm{Z}_0}{1}{\gamma}}{\starnorm{\bm{Z}_0}^{N+2}}
\prod_{i=1}^n \hdotnorm{\bm{Z}_i}{1}{\gamma}^{\alpha_i + \beta_i}h |\theta - \theta'|^{\gamma-2}\\
\label{logest4}
\abs{\diff_h \p_\theta f(\theta,\theta')}&\leq
C\frac{\hdotnorm{\bm{Z}_k}{1}{\gamma}\hdotnorm{\bm{Z}_0}{1}{\gamma}^2}{\starnorm{\bm{Z}_0}^{N+3}}\prod_{i = 1}^{n}\hdotnorm{\bm{Z}_i}{1}{\gamma}^{\alpha_i + \beta_i}\paren{h^\gamma\abs{\theta-\theta'}^{-2}+h\abs{\theta-\theta'}^{\gamma-3}},
\end{align}
where the constant $C$ does not depends only on $\gamma$.
\end{lemma}
\begin{poof}
Express $f$ as in \eqref{fqg}. 
The bound \eqref{logest3} is a simple consequence of \eqref{logest2}, and can be obtained
in the same way as we obtained \eqref{e:gestmult2} from \eqref{e:gestmult1} in the proof of Lemma \ref{GPL2}.

To obtain \eqref{logest4}, we write $\diff_h \p_\theta f(\theta,\theta')$ as:
\begin{equation*}
\begin{split}
\diff_h \p_\theta f(\theta,\theta')&=B_1+B_2+B_3+B_4,\\
B_1&=(\diff_h q)(\trl_h \p_\theta g), \; B_2=q(\diff_h \p_\theta g),\\
B_3&=(\diff_h \p_\theta q)(\trl_h g), \; B_4=(\p_\theta q)(\diff_h g).
\end{split}
\end{equation*}
We estimate each term. First note that
\begin{equation*}
\begin{split}
\abs{\diff_h q} \leq h\int_0^1\abs{\p_\theta q (\theta+sh,\theta')}ds
&\leq C\frac{\hdotnorm{V}{1}{\gamma}\hdotnorm{\bm{Z}_0}{1}{\gamma}}{\starnorm{\bm{Z}_0}^2}\int_0^1\abs{\theta+sh-\theta'}^{\gamma-2}ds\\
&\leq C\frac{\hdotnorm{V}{1}{\gamma}\hdotnorm{\bm{Z}_0}{1}{\gamma}}{\starnorm{\bm{Z}_0}^2}h\abs{\theta-\theta'}^{\gamma-2}.
\end{split}
\end{equation*}
Using \eqref{e:gestmult1} and \eqref{thetatheta'h}, we have:
\begin{equation*}
\abs{\trl_h \p_\theta g}\leq C\frac{\hdotnorm{\bm{Z}_0}{1}{\gamma}}{\starnorm{\bm{Z}_0}^{N+1}}\prod_{i = 1}^{n}\hdotnorm{\bm{Z}_i}{1}{\gamma}^{\alpha_i + \beta_i}|\theta - \theta'|^{\gamma - 1}
\end{equation*}
Thus, 
\begin{equation*}
\abs{B_1}\leq C\frac{\hdotnorm{\bm{Z}_k}{1}{\gamma}\hdotnorm{\bm{Z}_0}{1}{\gamma}^2}{\starnorm{\bm{Z}_0}^{N+3}}
\prod_{i = 1}^{n}\hdotnorm{\bm{Z}_i}{1}{\gamma}^{\alpha_i + \beta_i}h|\theta - \theta'|^{2\gamma - 3}
\end{equation*}
To obtain $B_2$, we must estimate $\diff_h \p_\theta g$. Let us use the notation in \eqref{gderiv}. Consider $A_l$ and write this 
as $A_l=(\p_\theta \phi_l) g_l$ (as in \eqref{delA_k}). We have:
\begin{equation*}
\diff_h A_l=(\diff_h \p_\theta\phi_l)(\trl_h g_l)+\p_\theta\phi_l (\diff_h g_l).
\end{equation*}
Using \eqref{diffdelphi} and \eqref{gest0}, we have:
\begin{equation*}
\abs{\diff_h \p_\theta\phi_l}\abs{\trl_h g_l}\leq 
C\frac{\hdotnorm{\bm{Z}_0}{1}{\gamma}^2}{\starnorm{\bm{Z}_0}^{N+2}}
\prod_{i = 1}^{n}\hdotnorm{\bm{Z}_i}{1}{\gamma}^{\alpha_i + \beta_i}
\paren{h^\gamma\abs{\theta-\theta'}^{-1}+h\abs{\theta-\theta'}^{\gamma-2}}
\end{equation*}
Using \eqref{diffhp'g} and \eqref{phi1est}, we have:
\begin{equation*}
\abs{\p_\theta\phi_l}\abs{\diff_h g_l}\leq C\frac{\hdotnorm{\bm{Z}_0}{1}{\gamma}^2}{\starnorm{\bm{Z}_0}^{N+2}}\prod_{i = 1}^{n}\hdotnorm{\bm{Z}_i}{1}{\gamma}^{\alpha_i + \beta_i}h\abs{\theta-\theta'}^{2\gamma - 2}
\end{equation*}
Combining the above with \eqref{q0est}, we have:
\begin{equation*}
\abs{B_2}\leq 
C\frac{\hdotnorm{\bm{Z}_k}{1}{\gamma}\hdotnorm{\bm{Z}_0}{1}{\gamma}^2}{\starnorm{\bm{Z}_0}^{N+3}}\prod_{i = 1}^{n}\hdotnorm{\bm{Z}_i}{1}{\gamma}^{\alpha_i + \beta_i}
\paren{h\abs{\theta-\theta'}^{3\gamma - 3}+ h^\gamma\abs{\theta-\theta'}^{\gamma-2}+h\abs{\theta-\theta'}^{2\gamma-3}}
\end{equation*}
Let us consider $B_3$. Using \eqref{diffq1est} and \eqref{gest0}, we have:
\begin{equation*}
\abs{B_3}=\abs{\diff_h \p_\theta q}\abs{\trl_h g}\leq 
C\frac{\hdotnorm{\bm{Z}_k}{1}{\gamma}\hdotnorm{\bm{Z}_0}{1}{\gamma}^2}{\starnorm{\bm{Z}_0}^3}
\prod_{i = 1}^{n}\hdotnorm{\bm{Z}_i}{1}{\gamma}^{\alpha_i + \beta_i}
\paren{h^\gamma\abs{\theta-\theta'}^{-2}+h\abs{\theta-\theta'}^{\gamma-3}}.
\end{equation*}
Finally, using \eqref{q1est} and \eqref{diffhp'g} (as applied to $g$), we have:
\begin{equation*}
\abs{B_4}\leq 
C\frac{\hdotnorm{\bm{Z}_k}{1}{\gamma}\hdotnorm{\bm{Z}_0}{1}{\gamma}^2}{\starnorm{\bm{Z}_0}^{N+3}}\prod_{i = 1}^{n}\hdotnorm{\bm{Z}_i}{1}{\gamma}^{\alpha_i + \beta_i}h\abs{\theta-\theta'}^{2\gamma - 3}.
\end{equation*}
Combining the estimates for $B_k, k=1,\ldots,4$ and taking only the leading order terms, we obtain the desired estimate.
\end{poof}

\section{Local Existence}\label{LocExst}

\subsection{Linear Term and Semigroup Estimates}\label{SGestSec}

The operator $\Lambda u := -\frac{1}{4}\Hil{\partial_{\theta}u}$ is the principal linear part of our evolution. 
The goal of this section is essentially to prove that $\Lambda$ generates an analytic semigroup on the scale 
of H\"older spaces. See \cite{lunardi} for an in-depth discussion of analytic semigroups on H\"older spaces. 
Given that $\Lambda$ is just the square root of the one-dimensional Laplacian 
and is also the Dirichlet-to-Neumann operator on a disk (up to a multiplicative constant), 
the results we prove in this section are likely well-known to students of this operator. 
Keeping with the elementary nature of our exposition, we prove all the requisite estimates, 
mostly by elementary means.

Suppose we are given a continuous function $u$ on the unit circle $\mathbb{S}^1$. We may express $u$
in terms of Fourier series:
\begin{equation*}
u(\theta)=\sum_{k=-\infty}^\infty u_ke^{ik\theta}, \; u_k = \frac{1}{2\pi}\ds\int_{0}^{2\pi} u(\theta ')e^{-ik\theta '}d\theta '.
\end{equation*}
It is well-known that the Fourier symbol of the operator $\Lambda u=-\frac{1}{4}\Hil{\partial_{\theta}u}$ 
is given by $-\abs{k}/4$ where $k$ is the Fourier wave number. That is to say, 
\begin{equation*}
\Lambda u=-\frac{1}{4}\sum_{k=-\infty}^\infty \abs{k}u_ke^{ik\theta},
\end{equation*}
for sufficiently smooth $u$.
The operator $e^{t\Lambda}, t>0$ is therefore:
\begin{equation*}
e^{t\Lambda}u=\sum_{k=-\infty}^\infty u_ke^{ik\theta-\abs{k}t/4}.
\end{equation*}
For $t>0$, we have:
\begin{equation}\label{poisson_kernel}
\begin{split}  
e^{t\Lambda}u &= \frac{1}{2\pi} \int_{\mbs} \sum_{k = -\infty}^{\infty}e^{-|k|t/4}e^{ik(\theta - \theta ')} u(\theta') d\theta' \\
&= \frac{1}{2\pi} \int_{\mbs} P(e^{-t/4},\theta-\theta') u(\theta ')d\theta',\;
P(r,\theta)=\frac{1 - r^2}{1 - 2r\cos(\theta) + r^2}.
\end{split}
\end{equation}
Note that $P(r, \theta)$ is the Poisson kernel.

Another useful way to view $\Lambda$ is as the Dirichlet-to-Neumann map on the unit disk. 
Consider the following Laplace boundary value problem:
\begin{equation*}\label{laplaceBVP}
\Delta v=0, \text{ for } \mathbb{D}^2, \; v=u(\theta) \text{ on } \partial \mathbb{D}=\mathbb{S}^1,
\end{equation*}
where $\mathbb{D}^2$ is the open unit disk and $\theta$ is the angular coordinate. Define the 
Dirichlet-to-Neumann map $\Lambda_{\rm DN}$ as:
\begin{equation*}
\Lambda_{\rm DN}: u\mapsto \left.\PD{v}{r}\right|_{r=1}
\end{equation*}
where $r$ is the radial coordinate. Then $\Lambda=-\Lambda_{\rm DN}/4$. This explains 
the appearance of the Poisson kernel in the expression of $e^{t\Lambda}$ above.

We first state a result on the mapping properties of $\Lambda$ on H\"older spaces.
\begin{proposition}
The operator $\Lambda$ is a bounded operator from $C^{k+1,\gamma}(\mathbb{S}^1)$ to $C^{k,\gamma}(\mathbb{S}^1)$ for 
$k=0,1,2,\ldots$ and $0<\gamma<1$.
\end{proposition}
\begin{poof}
This is a consequence of the well-known fact that the Hilbert transform is a bounded map from $C^{k,\gamma}(\mbs)$ to itself.
\end{poof}

We turn our attention to the semigroup $e^{t\Lambda}$.
\begin{proposition}\label{etlamkatint}
\item For $u \in C^k(\mathbb{S}^1), \; k\in \{0\}\cup \mathbb{N}$, 
\begin{align}
\chnorm{e^{t\Lambda } u}{k} &\leq \chnorm{u}{k}, \text{ for } 0<t,\label{etlamck}\\
\chnorm{e^{t\Lambda} u}{k+1} &\leq \frac{C}{t}\chnorm{u}{k} \text{ for } 0<t<1,\label{etlamck+1}
\end{align}
where $C$ above is a constant that does not depend on $u$ (or $k$).
\end{proposition}
\begin{poof}
Let $r=e^{-t/4}$ in what follows.
Inequality \eqref{etlamck} for $k=0$ is a simple consequence of the well-known fact that
\begin{equation*}\label{IntTo1}
\frac{1}{2\pi}\int_{0}^{2\pi}P(r, \theta')d\theta' = 1, \; P(r,\theta)>0.
\end{equation*}
Indeed, 
\begin{equation*}
\abs{(e^{t\Lambda} u)(\theta)}\leq \frac{1}{2\pi} \int_{0}^{2\pi} \abs{P(r, \theta - \theta')}\abs{u(\theta ')}d\theta'\leq 
\frac{1}{2\pi} \int_{0}^{2\pi}P(r, \theta)d\theta \norm{u}_{C^0}=\norm{u}_{C^0}.
\end{equation*}
What we proved above is just the maximum principle for harmonic functions on a disk.
For $k>0$, note that
\begin{equation} \label{etlamdelkcommute}
\p_\theta^k \paren{e^{t\Lambda}u}=e^{t\Lambda}\paren{\p_\theta^k u }.
\end{equation}
Thus, 
\begin{equation}\label{snormCkbound}
\snorm{e^{t\Lambda}u}_{C^k}=\norm{e^{t\Lambda}\paren{\p_\theta^k u }}_{C^0}\leq \norm{\p_\theta^k u }_{C^0}=\snorm{u}_{C^k},
\end{equation}
where inequality \eqref{etlamck} for $k=0$ was used in the inequality above.
This concludes the proof of \eqref{etlamck}. 

Consider \eqref{etlamck+1} for $k=0$. 
\begin{equation*}
\abs{\p_\theta (e^{t\Lambda} u)}=\abs{\frac{1}{2\pi}\int_0^{2\pi} \p_\theta P(r,\theta-\theta')u(\theta')d\theta'}
\leq \frac{1}{\pi}\int_0^{\pi}\abs{\p_\theta P(r,\theta)}d\theta \norm{u}_{C^0},
\end{equation*}
where we used the symmetry of $P$ with respect to $\theta$ in the last inequality.
Note that $P(r,\theta)$ is a decreasing function of $\theta$ from $0$ to $\pi$. Thus,
\begin{equation*}
\int_0^{\pi}\abs{\p_\theta P(r,\theta)}d\theta=-\int_0^{\pi} \p_\theta P(r,\theta)d\theta=P(r,0)-P(r,\pi)=\frac{2(1+r^2)}{1-r^2}.
\end{equation*}
Thus, we have:
\begin{equation*}
\snorm{e^{t\Lambda}u}_{C^1}\leq \frac{2(1+e^{-t/2})}{\pi(1-e^{-t/2})}\norm{u}_{C^0}\leq \frac{2}{\pi}\paren{1+\frac{4}{t}}\norm{u}_{C^0}.
\end{equation*}
Using \eqref{etlamdelkcommute} and proceeding as in \eqref{snormCkbound}, we find:
\begin{equation*}
\snorm{e^{t\Lambda}u}_{C^{k+1}}\leq \frac{2}{\pi}\paren{1+\frac{4}{t}}\snorm{u}_{C^k}.
\end{equation*}
Inequality \eqref{etlamck+1} now follows with $C=1+10/\pi$.
\end{poof}

To state the next proposition, we introduce the following notation.
For $\alpha>0, \alpha\notin \mathbb{N}$, we let 
$C^\alpha(\mathbb{S}^1)=C^{\floor{\alpha},\alpha-\floor{\alpha}}(\mathbb{S}^1)$ where $\floor{\alpha}$
is the largest integer smaller than $\alpha$. This next proposition will be key for proving the local existence of mild solutions. It will also be used extensively in section \ref{s:smoothness}. 
\begin{proposition}[H\"older estimates on Semigroup]\label{holderSG}
Let $u\in C^\alpha(\mathbb{S}^1), \alpha\geq 0$ and let $\beta\geq 0$ satisfy $0\leq \beta-\alpha\leq 1$. Then, 
\begin{equation}\label{mainSGest}
\norm{e^{t\Lambda}u}_{C^\beta}\leq \frac{C}{t^{\beta-\alpha}}\norm{u}_{C^\alpha},\; 0<t<1,
\end{equation}
where the constant $C$ above depends only on $\alpha$ and $\beta$.
\end{proposition}

To prove the above proposition, we shall make use of the following standard interpolation result 
which can be found, for example, in Chapter 1 of \cite{lunardi}.
\begin{proposition}[Interpolation of Bounded Operators on H\"older Spaces]\label{holderinterpolation}
Let $0\leq \alpha_0\leq \alpha_1$ and $0\leq \beta_0\leq \beta_1$ 
and let $\mc{T}$ be a bounded operator from 
$C^{\alpha_i}(\mathbb{S}^1)$ to $C^{\beta_i}(\mathbb{S}^1), i=0,1$ so that:
\begin{equation*}
\norm{\mc{T}u}_{C^{\beta_i}}\leq K_i\norm{u}_{C^{\alpha_i}}, \; i=0,1,
\end{equation*}
where $u\in C^{\alpha_i}(\mathbb{S}^1)$ and the constant $K_i>0$ does not depend on $u$.
Let $0<\sigma<1$ and $\alpha_\sigma=(1-\sigma)\alpha_0+\sigma\alpha_1, \; \beta_\sigma=(1-\sigma)\beta_0+\sigma\beta_1$.
Suppose one of the following conditions is satisfied.
\begin{enumerate}[label=(\roman*)]
\item $\alpha_\sigma\notin\mathbb{N}$ and $\beta_\sigma\notin\mathbb{N}$.
\item $\alpha_0=\alpha_1$ and $\beta_\sigma\notin \mathbb{N}$.
\item $\alpha_\sigma\notin\mathbb{N}$ and $\beta_0=\beta_1$.
\end{enumerate}
Then, $\mc{T}$ defines a bounded operator from $C^{\alpha_\sigma}(\mathbb{S}^1)$ 
to $C^{\beta_\sigma}(\mathbb{S}^1)$ so that:
\begin{equation*}
\norm{\mc{T}u}_{C^{\beta_\sigma}}\leq CK_0^{1-\sigma}K_1^\sigma\norm{u}_{C^{\alpha_\sigma}}
\end{equation*}
where $u\in C^{\alpha_\sigma}(\mathbb{S}^1)$ and the constant $C>0$ does not depend on $\mc{T}$ (or $u$) and depends only on 
$\alpha_i, \beta_i, i=0,1$ and $\sigma$.
\end{proposition}
\begin{remark}
The restriction to non-integer values in the above proposition can be lifted if we replace the definition 
of $C^k$ spaces for integer $k$ with H\"older-Zygmund spaces. We do not need such results in this paper. 
See \cite{lunardi} for details. 
\end{remark}
\begin{poof}[Proof of Proposition \ref{holderSG}]
Consider \eqref{etlamck} at two adjacent integer levels:
\begin{equation*}
\norm{e^{t\Lambda}u}_{C^k}\leq \norm{u}_{C^k}, \; \norm{e^{t\Lambda}u}_{C^{k+1}}\leq \norm{u}_{C^{k+1}}.
\end{equation*}
Interpolating between these two inequalities using Proposition \ref{holderinterpolation}, we obtain 
\begin{equation}\label{alphaequal}
\norm{e^{t\Lambda}u}_{C^\alpha}\leq C\norm{u}_{C^\alpha}
\end{equation}
for $k<\alpha<k+1$, where $C$ depends only on $\alpha$. Since $k$ was arbitrary, this establishes \eqref{mainSGest}
for $\beta=\alpha$. Likewise, considering \eqref{etlamck+1} at adjacent integer values, we obtain:
\begin{equation}\label{alpha+1}
\norm{e^{t\Lambda}u}_{C^{\gamma+1}}\leq \frac{C}{t}\norm{u}_{C^\gamma}
\end{equation}
for any $\gamma\geq 0$. This establishes \eqref{mainSGest} for $\beta=\alpha+1$. 
Interpolating between \eqref{alphaequal} and \eqref{alpha+1} setting $\alpha=\gamma$, we obtain the rest of 
the cases of \eqref{mainSGest} so long as $\beta$ is not an integer. 
To handle the case when $\beta$ is an integer, interpolate \eqref{alphaequal} and \eqref{alpha+1} setting $\alpha=\gamma+1$.
\end{poof}
\begin{remark}
The estimates in proposition \eqref{holderSG} can also be derived directly by elementary means without recourse to 
function space interpolation. 
\end{remark}

Finally, we turn to strong continuity of the the semigroup operator $e^{t\Lambda}$. 
To this end, we introduce the little H\"older spaces $h^{k,\gamma}(\mathbb{S}^1), k=0,1,2,\cdots, 0<\gamma<1$
as the completion of $C^\infty(\mathbb{S}^1)$ (smooth functions on the unit circle) 
in $C^{k,\gamma}(\mathbb{S}^1)$. A function $u\in C^{k,\gamma}(\mathbb{S}^1)$ belongs to 
$h^{k,\gamma}(\mathbb{S}^1)$ if and only if:
\begin{equation*}
\lim_{\delta\searrow 0} \sup_{\abs{\theta-\theta'}<\delta}\frac{\abs{u^{(k)}(\theta)-u^{(k)}(\theta')}}{\abs{\theta-\theta'}^\gamma}=0,
\end{equation*}
where $u^{(k)}$ is the $k$-th derivative of $u$.
From this, it is immediate that $C^{k,\beta}(\mathbb{S}^1)$ embeds continuously into $h^{k,\alpha}(\mathbb{S}^1)$ for any $\beta>\alpha$.
We refer to Chapter 0 of \cite{lunardi} for a discussion of these issues.

Given that the the space of smooth functions is not dense in $C^{k,\gamma}(\mathbb{S}^1)$, we can only expect
strong continuity of $e^{t\Lambda}$ in the little H\"older space $h^{k,\gamma}(\mathbb{S}^1)$.
\begin{proposition}\label{semigroup_strong_continuity}
Let $u\in C^{k,\gamma}(\mathbb{S}^1), \; k=0,1,2,\ldots$ and $0<\gamma<1$. Then, 
\begin{equation}\label{hkcont}
\lim_{t\searrow 0}\norm{e^{t\Lambda}u-u}_{C^{k,\gamma}}=0
\end{equation}
if and only if $u\in h^{k,\gamma}(\mathbb{S}^1)$.
\end{proposition}
\begin{poof}
Given that $e^{t\Lambda}u$ is a smooth function for $t>0$, \eqref{hkcont} can only be true if $u\in h^{k,\gamma}(\mathbb{S}^1)$.
To prove the converse, suppose $u\in h^{k,\gamma}(\mathbb{S}^1)$. Then, by definition, for any $\epsilon>0$,
there is a function $u_\epsilon\in C^\infty(\mathbb{S}^1)$ such that $\norm{u-u_\epsilon}_{C^{k,\gamma}}\leq \epsilon$.
By the well-known properties of the Poisson kernel,
\begin{equation*}
\lim_{t\searrow 0}\norm{e^{t\Lambda}u_\epsilon-u_{\epsilon}}_{C^{k+1}}=0.
\end{equation*}
Since the $C^{k+1}$ norm dominates the $C^{k,\gamma}$ norm, the above is true in the $C^{k,\gamma}$ norm.
The desired result now follows from the boundedness of $e^{t\Lambda}$ as a bounded operator 
from $C^{k,\gamma}(\mathbb{S}^1)$ to itself, as shown in \eqref{mainSGest} (or \eqref{alphaequal}).
\end{poof}

\subsection{Nonlinear Remainder Terms}\label{NonlinSec}

Consider the term $\mc{R}$ in \eqref{SSD}. We may write 
\begin{equation}\label{RCRT}
\mc{R} = \mc{R}_C + \mc{R}_T,
\end{equation}
where
\begin{align}\label{R_Cdefn}
\mc{R}_{C}(\bm{X})(\theta) &= -\frac{1}{4\pi}\int_{\mbs}K_C(\theta,\theta')\p_{\theta'}\bm{X}' d\theta',\\\label{K_Cdefn}
K_C (\theta, \theta ' ) &= \frac{\Delta \bm{X}\cdot \p_{\theta'}\bm{X}'}{\abs{\Delta \bm{X}}^2} - \frac{1}{2}\cot\left( \frac{\theta - \theta ' }{2}\right),
\end{align}
and
\begin{align}\label{R_Tdefn}
\mc{R}_{T}(\bm{X})(\theta) &= -\frac{1}{4\pi}\int_{\mbs}K_T(\theta,\theta')\p_{\theta'}\bm{X}'d\theta ',\\
\label{K_Tdefn}
K_T(\theta,\theta')&=\p_\theta'\paren{\frac{\Delta \bm{X}\otimes \Delta \bm{X}}{\abs{\Delta \bm{X}}^2}}.
\end{align}
Suppose $\bm{X}\in C^{1,\gamma}, 0<\gamma<1$. The kernels $K_C$ and $K_T$ then
behave like $\abs{\theta-\theta'}^{\gamma-1}$ for small $\abs{\theta-\theta'}$. This follows from Lemmas \ref{GPL1}
and \ref{GPlog1} of section \ref{s:calculus}. This shows that the above integrals are well-defined for $\bm{X}\in C^{1,\gamma}(\mbs)$.
Furthermore, this suggests the following. 
Since $\bm{X}_\theta\in C^{0,\gamma}$ and the kernels behave like $\abs{\theta-\theta'}^{\gamma-1}$, both remainder 
terms $\mc{R}_C$ and $\mc{R}_T$ resemble the $\gamma$-order fractional integration operator acting on a $C^{0,\gamma}(\mbs)$ function.
This suggests $\mc{R}_C$ and $\mc{R}_T$ may lie in $C^{\floor{2\gamma},2\gamma - \floor{2\gamma}}$.
Our goal is to show that this is indeed the case. Let
\begin{align*}
F_C[\bm{u}](\theta) &= -\frac{1}{4\pi}\int_{\mbs}K_C(\theta,\theta') \bm{u}(\theta ') d\theta ',\\
F_{T}[\bm{u}](\theta) &= -\frac{1}{4\pi}\int_{\mbs}K_T(\theta,\theta') \bm{u}(\theta ')d\theta '.
\end{align*}
We will show that both $F_C$ and $F_T$ map
$C^{0, \gamma}(\mbs)$ to $C^{0, 2\gamma}(\mbs)$ for $\gamma < 1/2$ and 
$C^{0, \gamma}(\mbs)$ to $C^{1, 2\gamma-1}(\mbs)$ for $\gamma > 1/2$. In the case of $\gamma = 1/2$, we instead have that $F_C$ and $F_T$ map $C^{1, 1/2}$ to $C^{0, \alpha}$ for any $\alpha \in (0, 1)$.
The key observation is that the kernel of each operator is a perfect derivative in $\theta '$, 
which implies that each of these operators maps constants to zero.
This allows us to make use of the H\"{o}lder continuity of the function being operated on,
thus leading to a gain of regularity by $\gamma$.

We first treat the operator $F_C$.
\begin{proposition}\label{C0cotgain}
If $\bm{X} \in C^{1, \gamma}(\mbs)$ with $\gamma \in (0,1)$ and $\starnorm{\bm{X}} > 0$, $\bm{u} \in C^{0, \gamma}(\mbs)$ and

\begin{enumerate}[label=(\roman *)]
\item\label{loggamma_neq1/2} if $\gamma \neq 1/2$ then $F_C [\bm{u}] \in C^{\floor{2\gamma}, 2\gamma - \floor{2\gamma}}(\mbs)$ with 

\[
\hflnorm{F_C [ \bm{u} ]}{2\gamma} \leq C \frac{\hdotnorm{\bm{X}}{1}{\gamma}^3}{\starnorm{\bm{X}}^3}\hnorm{\bm{u}}{0}{\gamma}.
\]

\item\label{loggamma=1/2} if $\gamma = 1/2$ then $F_C [\bm{u}] \in C^{0, \alpha}(\mbs)$ for any $\alpha\in(0, 1)$ with 

\[
\hnorm{F_C [ \bm{u} ] }{0}{\alpha} \leq C \frac{\hdotnorm{\bm{X}}{1}{\gamma}^3}{\starnorm{\bm{X}}^3} \hnorm{\bm{u}}{0}{\gamma}.
\]
\end{enumerate}
In the above, the constant $C$ does not depend on $\bm{X}$ or $\bm{u}$.
\end{proposition}
In fact, we have the following somewhat stronger bound:
\begin{equation}\label{FCCdot}
\hflnorm{F_C [ \bm{u} ]}{2\gamma} \leq C \frac{\hdotnorm{\p_\theta\bm{X}}{0}{\gamma}^3}{\starnorm{\bm{X}}^3}\hnorm{\bm{u}}{0}{\gamma}
\end{equation}
This is easily seen from the proof, but can also be seen from the fact that the kernel $K_C$ is invariant under translation.
We also note that when $\gamma \leq 1/2$ we can reduce the exponents of $\hnorm{\bm{X}}{1}{\gamma}$ and $\starnorm{\bm{X}}$ by one power each.

\begin{poof}
[Proof of Proposition \ref{C0cotgain}] Let us first collect some facts about the kernel $K_C$ and the operator $F_C$. Write $K_C$ as:
\begin{equation}\label{K_CandK_Ldefn}
\begin{aligned}
K_C(\theta, \theta ') &= \left(\frac{\Delta \bm{X}\cdot \p_{\theta'}\bm{X}'}{\abs{\Delta \bm{X}}^2} - \frac{1}{\theta - \theta '} \right) + \left(\frac{1}{\theta - \theta '} - \frac{1}{2}\cot \left(\frac{\theta - \theta '}{2} \right)\right)\\
&=: K_{L}(\theta, \theta ') + R_C (\theta - \theta ').
\end{aligned}
\end{equation}
Let $\bm{X}(\theta)=(X(\theta),Y(\theta))$.  Note that $K_L$ can be written as:
\begin{equation}\label{KLXY}
K_L=(K_L^X+K_L^Y),\; 
K_L^X(\theta,\theta')=\frac{(X_\theta'-(\Delta X/(\theta-\theta')))\Delta X}{(\Delta X)^2+(\Delta Y)^2},
\end{equation}
where $\Delta X=X(\theta)-X(\theta')$ and likewise for $\Delta Y$. $K_L^Y$ is simply the expression obtained 
by swapping $X$ and $Y$ in $K_L^X$. Notice that $K_L^X$ is of the form \eqref{loglemmaf}
in Lemma \ref{GPlog1} with $V=X, \bm{Z}_0=\bm{X}$ and $n=0, \alpha_0=1, \beta_0=0$.
We may thus apply inequality \eqref{logest1} of Lemma \ref{GPlog1} to obtain the estimate:
\begin{equation*}
|K_L^X (\theta, \theta ' )| \leq \frac{\hdotnorm{\bm{X}}{1}{\gamma}}{\starnorm{\bm{X}}} \abs{\theta - \theta ' }^{\gamma-1}.
\end{equation*}
for some constant $C>0$. 
The same estimate holds for $K_L^Y$; indeed, $K_L^Y$ is of the form \eqref{loglemmaf}
with $V=Y, \bm{Z}_0=\bm{X}$ with $n=0, \alpha_0=1, \beta_0=0$. Thus, returning to \eqref{KLXY}, we obtain the estimate:
\begin{equation*}
|K_L (\theta, \theta ' )| \leq \frac{\hdotnorm{\bm{X}}{1}{\gamma}}{\starnorm{\bm{X}}} \abs{\theta - \theta ' }^{\gamma-1}.
\end{equation*}
As for $R_C$ of \eqref{K_CandK_Ldefn}, simple calculus shows that 
\begin{equation*}
\left|R_C(\theta-\theta ')\right| \leq C|\theta - \theta '|,
\end{equation*}
for some constant $C$.
Plugging this back into (\ref{K_CandK_Ldefn}) gives
\begin{equation}\label{absKCbound}
|K_C (\theta, \theta ')| \leq C\frac{\hdotnorm{\bm{X}}{1}{\gamma}}{\starnorm{\bm{X}}}|\theta - \theta '|^{\gamma - 1}.
\end{equation}

An important property of the operator $F_C$ is that it annihilates constants. 
Since $K_C(\theta, \theta ')$ is a perfect derivative, for any constant $C$ we have
\begin{equation}\label{constto0inlog}
(F_C [C])(\theta) = \frac{1}{4\pi}\int_{\mbs} \p_{\theta'}  \log\paren{\frac{\abs{\Delta \bm{X}}}{\abs{\sin((\theta-\theta')/2)}}}Cd\theta'=0
\end{equation}
by the positivity, periodicity and continuity of the argument of the logarithm. 
Note that positivity is a result of our assumption $\starnorm{\bm{X}} > 0$.
In particular, we have 
\begin{equation}\label{FCminusC}
(F_C [\bm{u}] )(\theta) = -\frac{1}{4\pi}\int_{\mbs} K_C(\theta, \theta')\bm{u}(\theta') d\theta'= \frac{1}{4\pi}\int_{\mbs} K_C (\theta, \theta') (\bm{V} - \bm{u}(\theta'))d\theta',
\end{equation}
for any constant vector $\bm{V}$. 

We turn to statement \ref{loggamma_neq1/2}. We start with the case where $\gamma < 1/2$. 
Note that
\begin{equation}\label{FCubounded}
\abs{(F_C [ \bm{u} ] )(\theta)}\leq \frac{1}{4\pi}\int_{\mbs}\abs{K_C(\theta,\theta') \bm{u}(\theta')}d\theta'\leq \frac{1}{4\pi}\norm{\bm{u}}_{C^0} \int_{\mbs}\abs{K_C(\theta,\theta')}d\theta'
\leq C\frac{\hdotnorm{\bm{X}}{1}{\gamma}}{\starnorm{\bm{X}}}\norm{\bm{u}}_{C^0},
\end{equation}
where we used \eqref{absKCbound} and $C$ does not depend on $\bm{u}$ or $\bm{X}$. We thus see that $F_C \bm{u}$ is bounded and thus well-defined.
We must bound the seminorm $\snorm{F_C [ \bm{u} ] }_{C^{0,2\gamma}}$.
Consider the difference between two points, $\theta$ and $\theta+h$. 
Then, 
\begin{align*}
\diff_h(F_C [ \bm{u} ] )&=F_C [ \bm{u} ] (\theta + h) - F_C [ \bm{u} ]  (\theta)\\
&= \frac{1}{4\pi}\int_{\mbs} \left( K_C (\theta+ h, \theta ') - K_C (\theta, \theta ')\right) (\bm{u}(\theta ) - \bm{u}(\theta')) d\theta '\\ 
& = \frac{1}{4\pi}\int_{\mc{I}_s} (\diff_h K_C)(\bm{u}(\theta ) - \bm{u}(\theta')) d\theta '
+ \frac{1}{4\pi}\int_{\mc{I}_f} (\diff_h K_C)(\bm{u}(\theta ) - \bm{u}(\theta')) d\theta '=: A + B,
\end{align*}
where we used the notation introduced in \eqref{e:Idefn} and \eqref{diff_trl}, and, in the first equality, we used \eqref{FCminusC} with 
$\bm{V}=\bm{u}(\theta)$ for both terms.
We start by bounding term $A$. Using bound \eqref{absKCbound} on $\mc{I}_s$ we have, 
\begin{equation}\label{boundAloglemma}
\begin{aligned}
|A| &\leq C\int_{\mc{I}_s}\abs{\diff_h K_C} \hnorm{\bm{u}}{0}{\gamma}|\theta - \theta '|^{\gamma} d\theta '
\leq C h^{\gamma}\hnorm{\bm{u}}{0}{\gamma}\int_{\mc{I}_s}(\abs{\trl_h K_C}+\abs{K_C})d\theta '\\
&\leq C h^{\gamma}\hnorm{\bm{u}}{0}{\gamma}\frac{\hdotnorm{\bm{X}}{1}{\gamma}}{\starnorm{\bm{X}}} \int_{\mc{I}_s} \left(\frac{1}{|\theta + h - \theta '|^{1 - \gamma}} + \frac{1}{|\theta - \theta '|^{1 - \gamma}}\right)d\theta '\\
&\leq C h^{2\gamma}\hnorm{\bm{u}}{0}{\gamma} \frac{\hdotnorm{\bm{X}}{1}{\gamma}}{\starnorm{\bm{X}}}.
\end{aligned}
\end{equation}
In the above and in the sequel, the constants $C$ may vary but do not depend on $\theta,\theta',u,\bm{X}$ or $h$.
Let us estimate $B$. From \eqref{logest3} of Lemma \eqref{GPlog2}, for $\theta'\in \mc{I}_f$, we have:
\begin{equation*}
\begin{aligned}
|\diff_h K_L(\theta,\theta')| \leq C\paren{\frac{\hdotnorm{\bm{X}}{1}{\gamma}^2}{\starnorm{\bm{X}}^2}} h|\theta - \theta '|^{\gamma - 2}.
\end{aligned}
\end{equation*}
Furthermore, it is easily seen that 
\begin{equation*}
\abs{\diff_h R_C(\theta-\theta')}\leq Ch.
\end{equation*}
We thus have, 
\begin{equation}\label{diffKCbound}
|\diff_h K_C(\theta,\theta')| \leq C\paren{\frac{\hdotnorm{\bm{X}}{1}{\gamma}^2}{\starnorm{\bm{X}}^2}} h|\theta - \theta '|^{\gamma - 2}.
\end{equation}
Hence,
\begin{equation*}
\begin{aligned}
|B| &\leq C\int_{\mc{I}_f} \abs{\diff_h K_C}\abs{\bm{u} (\theta)- \bm{u} (\theta')})d\theta'\\
&\leq \int_{\mc{I}_f}C\left( \frac{\hdotnorm{\bm{X}}{1}{\gamma}^2}{\starnorm{\bm{X}}^2} h|\theta - \theta '|^{\gamma - 2}\right) \hnorm{\bm{u}}{0}{\gamma} |\theta - \theta '|^{\gamma} d\theta '\\
&\leq 
C\frac{\hdotnorm{\bm{X}}{1}{\gamma}^2}{\starnorm{\bm{X}}^2} \hnorm{\bm{u}}{0}{\gamma}h\int_{\mc{I}_f}|\theta - \theta '|^{2\gamma - 2}  d\theta '\\
&\leq 
C\frac{\hdotnorm{\bm{X}}{1}{\gamma}^2}{\starnorm{\bm{X}}^2} \hnorm{\bm{u}}{0}{\gamma}h\int_{h/2}^{2\pi}|s|^{2\gamma - 2} ds
\leq C\frac{\hdotnorm{\bm{X}}{1}{\gamma}^2}{\starnorm{\bm{X}}^2} \hnorm{\bm{u}}{0}{\gamma}h^{2\gamma}.
\end{aligned}
\end{equation*}
Combining the bounds for $A$ and $B$, we find that:
\begin{equation*}
\abs{\diff_h (F_C [\bm{u}])}\leq C\frac{\hdotnorm{\bm{X}}{1}{\gamma}^2}{\starnorm{\bm{X}}^2}\hnorm{\bm{u}}{0}{\gamma}h^{2\gamma}.
\end{equation*}
This, together with \eqref{FCubounded} shows that $F_C u\in C^{0,2\gamma}(\mbs)$. Noting here that

\begin{align*}
\frac{\hnorm{\bm{X}}{1}{\gamma}}{\starnorm{\bm{X}}} \geq 1
\end{align*}
gives the desired bound.

We now consider the case where $\gamma > 1/2$. 
We must consider the derivative of $F_C [\bm{ u}]$. The derivative of kernel $\p_\theta K_C$ behaves like $\abs{\theta-\theta'}^{\gamma-2}$
(see \eqref{pthetaKCbound}). Thus, the derivative $\p_\theta F_C [\bm{u}]$ cannot be expressed simply 
as the integral operator against the kernel $\p_\theta K_C$ since it is too singular.
We will show that the derivative of $F_C [\bm{u}]$ exists and it is equal to:
\begin{equation}\label{Aopdef}
(\mc{A} \bm{u})(\theta)= \frac{1}{4\pi}\int_{\mbs} \p_\theta K_C(\theta,\theta')(\bm{u}(\theta)-\bm{u}(\theta'))d\theta'.
\end{equation}
We first check that this integral is well-defined. The kernel $\p_\theta K_C$ can be estimated by considering $\p_\theta K_L$ and $\p_\theta R_C$
separately. From \eqref{logest2}, we have the bound:
\begin{equation*}
\abs{\p_\theta K_L(\theta,\theta')}\leq \frac{2\hdotnorm{\bm{X}}{1}{\gamma}^2}{\starnorm{\bm{X}}^2}\abs{\theta-\theta'}^{\gamma-2}.
\end{equation*}
It is also easily seen that 
\begin{equation*}
\abs{\p_\theta R_C(\theta-\theta')}\leq C.
\end{equation*}
Combining the above, we have:
\begin{equation}\label{pthetaKCbound}
\abs{\p_\theta K_C(\theta,\theta')}\leq C\frac{\hdotnorm{\bm{X}}{1}{\gamma}^2}{\starnorm{\bm{X}}^2}\abs{\theta-\theta'}^{\gamma-2}.
\end{equation}
Thus, 
\begin{equation}
\begin{split}\label{Aopbound}
\abs{(\mc{A} \bm{u})(\theta)}&= \int_{\mbs} \abs{\p_\theta K_C(\theta,\theta')(\bm{u}(\theta')-\bm{u}(\theta))}d\theta'\\
&\leq C\frac{\hdotnorm{\bm{X}}{1}{\gamma}^2}{\starnorm{\bm{X}}^2}\norm{\bm{u}}_{C^{0,\gamma}}\int_{\mbs}\abs{\theta-\theta'}^{2\gamma-2}d\theta'
\leq C\frac{\hdotnorm{\bm{X}}{1}{\gamma}^2}{\starnorm{\bm{X}}^2}\norm{\bm{u}}_{C^{0,\gamma}}.
\end{split}
\end{equation}
where we used $2\gamma-1>0$ in the last inequality to ensure that the above integral is bounded.

We now show that the derivative of $F_C [\bm{u}]$ is indeed $\mc{A} \bm{u}$. Consider the expression:
\begin{equation*}
\begin{split}
h^{-1}\diff_h (F_C [\bm{u}])(\theta)-(\mc{A} \bm{u})(\theta)
&=\frac{1}{4\pi}\int_{\mbs} \paren{h^{-1}\diff_h K_C-\p_\theta K_C}(\theta,\theta')(\bm{u}(\theta)-\bm{u}(\theta'))d\theta'=I_1+I_2, \\
I_1&=\frac{1}{4\pi}\int_{\mc{I}_s} \paren{h^{-1}\diff_h K_C-\p_\theta K_C}(\theta,\theta')(\bm{u}(\theta)-\bm{u}(\theta'))d\theta',\\
I_2&=\frac{1}{4\pi}\int_{\mc{I}_f} \paren{h^{-1}\diff_h K_C-\p_\theta K_C}(\theta,\theta')(\bm{u}(\theta)-\bm{u}(\theta'))d\theta'.
\end{split}
\end{equation*}
In the above, we used \eqref{constto0inlog} or equivalently, \eqref{FCminusC}.
We estimate $I_1$ and $I_2$. Assume $h>0$. For $I_1$, using \eqref{absKCbound} and \eqref{pthetaKCbound} we have:
\begin{equation*}
\begin{split}
\abs{I_1}&\leq C\int_{\mc{I}_s} (h^{-1}(\abs{K_C}+\abs{\trl_h K_C})+\abs{\p_\theta K_C})\norm{\bm{u}}_{C^{0,\gamma}}\abs{\theta-\theta'}^{\gamma}d\theta'\\
&\leq C\frac{\hdotnorm{\bm{X}}{1}{\gamma}^2}{\starnorm{\bm{X}}^2}\norm{\bm{u}}_{C^{0,\gamma}}
\int_{\mc{I}_s}\paren{h^{-1}\paren{\abs{\theta-\theta'}^{2\gamma-1}+\abs{\theta+h-\theta'}^{2\gamma-1}}+\abs{\theta-\theta'}^{2\gamma-2}}d\theta'\\
&\leq C\frac{\hdotnorm{\bm{X}}{1}{\gamma}^2}{\starnorm{\bm{X}}^2}\norm{\bm{u}}_{C^{0,\gamma}}h^{2\gamma-1}.
\end{split}
\end{equation*}
To estimate $I_2$, note that, for $\theta'\in \mc{I}_f$:
\begin{equation*}
\abs{(h^{-1}\diff_h K_C-\p_\theta K_C)(\theta,\theta')}=\int_0^1 \abs{\p_\theta K_C(\theta+sh,\theta')-\p_\theta K_C(\theta,\theta')}ds
=\int_0^1 \abs{\diff_{sh} \p_\theta K_C}ds.
\end{equation*}
Using \eqref{logest4} of Lemma \ref{GPlog2}, we have:
\begin{equation*}
\abs{\diff_h \p_\theta K_L}\leq C\frac{\hdotnorm{\bm{X}}{1}{\gamma}^3}{\starnorm{\bm{X}}^3}\paren{h^\gamma\abs{\theta-\theta'}^{-2}+h\abs{\theta-\theta'}^{\gamma-3}}.
\end{equation*}
It is easily seen that this dominates $\diff_h \p_\theta R_C$ with a suitable constant. Therefore, 
\begin{equation}\label{diffpthetaKCbound}
\abs{\diff_h \p_\theta K_C}\leq C\frac{\hdotnorm{\bm{X}}{1}{\gamma}^3}{\starnorm{\bm{X}}^3}\paren{h^\gamma\abs{\theta-\theta'}^{-2}+h\abs{\theta-\theta'}^{\gamma-3}}.
\end{equation}
Thus, 
\begin{equation*}
\begin{split}
\abs{I_2}&\leq C\frac{\hdotnorm{\bm{X}}{1}{\gamma}^3}{\starnorm{\bm{X}}^3}\norm{\bm{u}}_{C^{0,\gamma}}
\int_{\mc{I}_f} \paren{h^\gamma\abs{\theta-\theta'}^{\gamma-2}+h\abs{\theta-\theta'}^{2\gamma-3}}d\theta'\\
&\leq C\frac{\hdotnorm{\bm{X}}{1}{\gamma}^3}{\starnorm{\bm{X}}^3}\norm{\bm{u}}_{C^{0,\gamma}}h^{2\gamma-1}.
\end{split}
\end{equation*}
We may obtain a similar estimate for $I_1$ and $I_2$ when $h<0$. Thus, since $2\gamma-1>0$,
\begin{equation}\label{pthetaFC=Au}
\p_\theta F_C[\bm{u}]=\lim_{h\to 0} \frac{1}{h}\diff_h (F_C [\bm{u}])=\mc{A} \bm{u}.
\end{equation}
We may now identify $\p_\theta F_C [\bm{u}]$ with the expression for $\mc{A} u$ in \eqref{Aopdef}.

We have only to bound the $C^{0,2\gamma-1}$ seminorm of $\p_\theta F_C [\bm{u}]$. Thus,
\begin{align*}
(\p_\theta F_C [\bm{u}] ) (\theta + h) - (\p_\theta F_C [\bm{u}]) (\theta) 
&= \frac{1}{4\pi}\int_{\mbs} \p_\theta K_C(\theta + h, \theta)(\bm{u} (\theta + h) - \bm{u} (\theta '))d\theta'\\
&- \frac{1}{4\pi}\int_{\mbs} \p_\theta K_C(\theta, \theta ')(\bm{u}(\theta) - \bm{u}(\theta '))d\theta'\\
&= \frac{1}{4\pi}\int_{\mc{I}_s}\p_\theta K_C(\theta, \theta ')(\bm{u}(\theta') - \bm{u}(\theta))d\theta'\\
&+ \frac{1}{4\pi}\int_{\mc{I}_s}\p_\theta K_C(\theta + h, \theta')(\bm{u}(\theta + h) - \bm{u}(\theta'))d\theta'\\
&+ \frac{1}{4\pi}\int_{\mc{I}_f}\p_\theta K_C(\theta, \theta')(\bm{u}(\theta + h) - \bm{u}(\theta))d\theta'\\
&+ \frac{1}{4\pi}\int_{\mc{I}_f}\paren{\p_\theta K_C(\theta + h, \theta ') - \p_\theta K(\theta, \theta')}(\bm{u}(\theta + h) - \bm{u}(\theta'))d\theta'\\
&\equiv I_3 + I_4 + I_5 + I_6
\end{align*}
where $\mc{I}_{s,f}$ were defined in \eqref{e:Idefn}.
We will bound each term individually and eventually recombine. 
For term $I_3$ we have, using \eqref{pthetaKCbound}, 
\begin{align*}
|I_3| &\leq C\int_{\mc{I}_s}|\p_\theta K_C(\theta, \theta ')||\bm{u}(\theta) - \bm{u}(\theta ')|d\theta'\\
& \leq C\hnorm{\bm{u}}{0}{\gamma}\frac{\hdotnorm{\bm{X}}{1}{\gamma}^{2}}{\starnorm{\bm{X}}^2}\int_{\mc{I}_s}|\theta - \theta '|^{2\gamma - 2}d\theta'
= C\hnorm{\bm{u}}{0}{\gamma}\frac{\hdotnorm{\bm{X}}{1}{\gamma}^{2}}{\starnorm{\bm{X}}^2}h^{2\gamma - 1}.
\end{align*}
Using the same argument, the same bound holds for term $I_4$. For term $I_5$, we may use \eqref{pthetaKCbound} to compute:
\begin{align*}
\abs{I_5} &\leq C\int_{\mc{I}_f}|\p_\theta K_C(\theta, \theta ')||\bm{u}(\theta)- \bm{u}(\theta + h)|d\theta '\\
&\leq Ch^{\gamma}\hnorm{\bm{u}}{0}{\gamma}\frac{\hdotnorm{\bm{X}}{1}{\gamma}^{2}}{\starnorm{\bm{X}}^2}\int_{\mc{I}_f} |\theta - \theta '|^{\gamma - 2}d\theta '
\leq C\hnorm{\bm{u}}{0}{\gamma}\frac{\hdotnorm{\bm{X}}{1}{\gamma}^{2}}{\starnorm{\bm{X}}^2}h^{2\gamma - 1}.
\end{align*}
For the term $I_6$, we may use \eqref{diffpthetaKCbound} to obtain:
\begin{equation*}
\begin{split}
\abs{I_6}&\leq C\int_{\mc{I}_f} \abs{\diff_h\p_\theta K_C}\abs{\bm{u}(\theta+h)-\bm{u}(\theta')}d\theta'\\
&\leq C\frac{\hdotnorm{\bm{X}}{1}{\gamma}^3}{\starnorm{\bm{X}}^3}\norm{\bm{u}}_{C^{0,\gamma}}
\int_{\mc{I}_f}\paren{h^\gamma\abs{\theta-\theta'}^{\gamma-2}+h\abs{\theta-\theta'}^{2\gamma-3}}\\
&\leq C\frac{\hdotnorm{\bm{X}}{1}{\gamma}^3}{\starnorm{\bm{X}}^3}\norm{\bm{u}}_{C^{0,\gamma}}h^{2\gamma-1}.
\end{split}
\end{equation*}
Collecting $I_3,\ldots, I_6$, we see that $\p_\theta F_C u$ is in $C^{0,2\gamma-1}$ with the desired estimate.

Note that the bound for statement \ref{loggamma=1/2} follows from a simple adaptation of the $\gamma < 1/2$ case.
\end{poof}

A similar result is true for the operator $F_T$.

\begin{proposition}\label{C0tengain}
Suppose $\bm{u}, \bm{X} \in C^{1, \gamma}(\mbs)$ with $\gamma \in (0, 1)$ and $\starnorm{\bm{X}} > 0$. We have the following estimates. 
\begin{enumerate}[label=(\roman*)]
\item if $\gamma \neq 1/2$ then $F_T[ \bm{u} ]\in C^{\floor{2\gamma}, 2\gamma - \floor{2\gamma}}$ with 

\[
\hflnorm{F_T [ \bm{u} ] }{2\gamma} \leq C\frac{\hdotnorm{\bm{X}}{1}{\gamma}^3}{\starnorm{\bm{X}}^3} \hnorm{\bm{u}}{0}{\gamma}.
\]

\item if $\gamma  = 1/2 $ then $F_T [ \bm{u} ] \in C^{0, \alpha}$ for any $\alpha\in (0, 1)$ with 

\[
\hnorm{F_T [ \bm{u} ] }{0}{\alpha} \leq C\frac{\hdotnorm{\bm{X}}{1}{\gamma}^3}{\starnorm{\bm{X}}^3} \hnorm{\bm{u}}{0}{\gamma}.
\]
\end{enumerate}
\end{proposition}
Like \eqref{FCCdot}, we have:
\begin{equation}\label{FTCdot}
\hflnorm{F_T [ \bm{u} ]}{2\gamma} \leq C \frac{\hdotnorm{\p_\theta\bm{X}}{0}{\gamma}^3}{\starnorm{\bm{X}}^3}\hnorm{\bm{u}}{0}{\gamma}.
\end{equation}
When $\gamma \leq 1/2$ we can reduce the exponents of $\hnorm{\bm{X}}{1}{\gamma}$ and $\starnorm{\bm{X}}$ by one power each.
 
\begin{poof}[Proof of Proposition \ref{C0tengain}]
The proof parallels that of the previous Proposition. The only difference is that we will be using Lemmas \ref{GPL1} and \ref{GPL2} to obtain the key estimates instead of Lemmas \ref{GPlog1} and \ref{GPlog2} which were used in the proof of the previous Proposition.

First, we note that for a constant vector $\bm{V}$:
\begin{equation*}
F_T[\bm{V}]= -\frac{1}{4\pi}\int_{\mbs} K_T(\theta,\theta')\bm{V}d\theta'=-\frac{1}{4\pi}\int_{\mbs} \p_\theta'\paren{\frac{\Delta \bm{X}\otimes \Delta\bm{X}}{\abs{\Delta \bm{X}}^2}}\bm{V}d\theta'=0,
\end{equation*}
where we used the fact that $K_T$ is the $\theta'$ derivative of a continuous function. Thus, like $F_C$, $F_T$ maps constant vectors to $0$.

Let $K_{T,ij}, i,j=1,2$ denote the matrix components of the kernel $K_T$. Let $\bm{X}=(X,Y)$.
We have:
\begin{equation*}
K_{T,11}=\p_{\theta'}g_{11}(\theta,\theta'),\; g_{11}(\theta,\theta')=\frac{(\Delta X)^2}{\abs{\Delta \bm{X}}^2}.
\end{equation*}
The function $g_{11}$ is in the form of $g(\theta,\theta')$ in Lemma \ref{GPL1} with $n=0, \alpha_0=2, \beta_0=0$. 
Similar considerations apply for $K_{T,12}=K_{T,21}$ and $K_{T,22}$. From Lemma \ref{GPL1}, we thus have the following estimates.
\begin{align}
\label{absKTbound}
\abs{K_{T,ij}}&\leq C\frac{\hdotnorm{\bm{X}}{1}{\gamma}}{\starnorm{\bm{X}}}\abs{\theta-\theta'}^{\gamma-1},\\
\label{pthetaKTbound}
\abs{\p_\theta K_{T,ij}}&\leq C\frac{\hdotnorm{\bm{X}}{1}{\gamma}^2}{\starnorm{\bm{X}}^2}\abs{\theta-\theta'}^{\gamma-2}.
\end{align}
From Lemma \ref{GPL2} we have the following estimates whenever $\theta'\in \mc{I}_f$:
\begin{align}
\label{diffKTbound}
\abs{\diff_hK_{T,ij}}&\leq C\frac{\hdotnorm{\bm{X}}{1}{\gamma}^2}{\starnorm{\bm{X}}^2}h\abs{\theta-\theta'}^{\gamma-2},\\
\label{diffpthetaKTbound}
\abs{\diff_h(\p_\theta K_{T,ij})}&\leq C\frac{\hdotnorm{\bm{X}}{1}{\gamma}^3}{\starnorm{\bm{X}}^3}\paren{h^\gamma\abs{\theta-\theta'}^{-2}+h\abs{\theta-\theta'}^{\gamma-3}}.
\end{align}
The rest of the proof is completely the same as the proof of the preceding Proposition. We may essentially replace $K_C$ with $K_T$ in the 
proof of Proposition \ref{C0cotgain}, replacing the use of the estimates
\eqref{absKCbound}, \eqref{pthetaKCbound}, \eqref{diffKCbound}, and \eqref{diffpthetaKCbound}
with the corresponding estimates
\eqref{absKTbound}, \eqref{pthetaKTbound}, \eqref{diffKTbound}, and \eqref{diffpthetaKTbound}.
\end{poof}
\subsection{Contraction Mapping}\label{MappingSec}

We are looking for a mild solution of the Peskin system which is a solution of

\begin{equation*}
\bm{X} (t) = e^{t\Lambda}\bm{X}_0 + \int_{0}^{t}e^{(t - s)\Lambda}\mc{R}(\bm{X}(s))ds,
\end{equation*}
where $\bm{X}: \mathbb{S}^1 \to \mathbb{R}^2$ is a closed curve and 

\begin{equation*}
\mc{R}(\bm{X}) = \mc{R}_{C}(\bm{X}) + \mc{R}_{T}(\bm{X})
\end{equation*}
as defined in (\ref{R_Cdefn}) and (\ref{R_Tdefn}). Using Propositions \ref{C0cotgain} and \ref{C0tengain} we have $\mc{R}: C^{1, \gamma}(\mathbb{S}^1)\mapsto C^{\lfloor 2\gamma \rfloor, 2\gamma - \lfloor 2\gamma \rfloor}(\mathbb{S}^1)$ if $\gamma \neq 1/2$ or $\mc{R}: C^{1, \gamma}(\mathbb{S}^1)\mapsto C^{0, \alpha}(\mathbb{S}^1)$ if $\gamma = 1/2$ for any $\alpha \in (0, 1)$. To address the local existence aspect of Theorem \ref{LWPTheorem}, we now show that the map
\begin{equation}\label{SXtX0}
S(\bm{X}, t; \bm{X}_0) := e^{\Lambda t}\bm{X}_0 + \int_{0}^{t}e^{(t - s)\Lambda}\mc{R}(\bm{X}(s))ds
\end{equation}
has a fixed point in a certain subset of $C([0,T];C^{1, \gamma}(\mathbb{S}^1))$ for suitably chosen values of $T$. We will use 

\begin{theorem}[Banach Fixed Point Theorem]
Let $(U, d)$ be a non-empty complete metric space. If $S : U \to U$ is an operator with $d(Su, Sv) \leq qd(u, v)$ for $q\in [0, 1)$, then $S$ has a unique fixed point $u^{*} \in U$. 
\end{theorem}
Given our previous estimates on $\mc{R}(\bm{X})$ and their dependence on $\starnorm{\bm{X}}$, we take a subset of $C^{1, \gamma}(\mathbb{S}^1)$ which includes only $\bm{Y} \in C^{1, \gamma}$ with $\starnorm{\bm{Y}} \geq m > 0$. We define our set as follows.

\begin{proposition}[Adapted from \cite{majdabertozzi} proposition 8.7]\label{closedsetO}
For any $M > m > 0$, the set $O^{M, m} := \{\bm{X}\in C^{1, \gamma}(\mc{D}) : \hnorm{\bm{X}}{1}{\gamma} \leq M \text{ and } \starnorm{\bm{X}} \geq m \}$ is closed in $C^{1,\gamma}(\mathbb{S}^1)$ . By extension, for any $T \geq t \geq 0$ the set $O^{M, m}_{t} = \{\bm{X} \in C([0, t]; C^{1, \gamma}(\mathbb{S}^1)): \bm{X}(s) \in O^{M, m}\text{ for all } s\in [0, t] \}$ is closed in $C([0, T]; C^{1, \gamma}(\mathbb{S}^1))$.
\end{proposition}

\begin{poof}
For $\bm{X}, \bm{Y} \in C^{1, \gamma}(\mathbb{S}^1)$, we have by the reverse triangle inequality

\begin{align*}
|\starnorm{\bm{X}} - \starnorm{\bm{Y}}| &= \left| \inf_{\theta \neq \theta '}\frac{|\bm{X} - \bm{X}'|}{|\theta - \theta '|} - \inf_{\theta \neq \theta '} \frac{|\bm{Y} - \bm{Y}'|}{|\theta - \theta '|}\right| \leq \sup_{\theta\neq \theta '} \left| \frac{|\bm{X} - \bm{X}'|}{|\theta - \theta '|} - \frac{|\bm{Y} - \bm{Y}'|}{|\theta - \theta '|}\right|\\
&\leq \sup_{\theta\neq \theta '} \frac{|\bm{X} - \bm{Y} - (\bm{X}' - \bm{Y}')|}{|\theta - \theta '|} \leq \hnorm{\bm{X} - \bm{Y}}{1}{\gamma}
\end{align*}

Also by the reverse triangle inequality,

\begin{equation*}
|\hnorm{\bm{X}}{1}{\gamma} - \hnorm{\bm{Y}}{1}{\gamma}| \leq \hnorm{\bm{X} - \bm{Y}}{1}{\gamma}.
\end{equation*}

Therefore the maps $\starnorm{\cdot} : C^{1, \gamma} \to [0, \infty)$ and $\hnorm{\cdot}{1}{\gamma}: C^{1, \gamma} \to [0, \infty)$ are continuous. Since $O^{M, m}$ is the intersection of preimages of two closed sets under continuous maps, it is closed in $C^{1, \gamma}$.  
Since $[0, t]$ is closed in $[0, T]$, the second statement follows. 
\end{poof}

To show that the map $S(\bm{X}, t; \bm{X}_0)$ is a contraction over the above set, we first prove the following:
\begin{proposition}\label{LipRemainder}
For any $M > m > 0$ and $\gamma \in(0, 1)$, the remainder term $\mc{R}(\bm{X})$ is Lipschitz on any convex set $\mc{B}\subset O^{M, m}$ and
\begin{enumerate}[label=(\roman *)]
\item if $\gamma \neq 1/2$, $\mc{R} : \mc{B} \to C^{\floor{2\gamma}, 2\gamma - \floor{2\gamma}}(\mathbb{S}^1)$ with
\begin{align*}
\hnorm{\mc{R}(\bm{X}) - \mc{R}(\bm{Y})}{\floor{2\gamma}}{2\gamma - \floor{2\gamma}} \leq C\frac{M^{4}}{m^{4}}\hnorm{\bm{X} - \bm{Y}}{1}{\gamma}.
\end{align*}
\item if $\gamma = 1/2$, then for any $\alpha \in (0, 1)$, $\mc{R} : \mc{B} \to C^{0, \alpha}(\mathbb{S}^1)$ with
\begin{align*}
\hnorm{\mc{R}(\bm{X}) - \mc{R}(\bm{Y})}{0}{\alpha} \leq C\frac{M^{4}}{m^{4}}\hnorm{\bm{X} - \bm{Y}}{1}{\gamma}.
\end{align*}
\end{enumerate}
\end{proposition}

\begin{poof}
We show only the $\gamma \neq 1/2$ case as the $\gamma = 1/2$ case follows from the arguments when $\gamma < 1/2$. Let $M > m > 0$ and $\mc{B} \subset O^{M, m}$ be convex. It suffices to show that the linearization of $\mc{R}$, $\p_{\bm{X}}\mc{R}(\bm{X})$, is bounded on $C^{\floor{2\gamma}, 2\gamma - \floor{2\gamma}}(\mathbb{S}^1)$ for any $\bm{X} \in \mc{B}$ since

\begin{align*}
\hnorm{\mc{R}(\bm{X}) - \mc{R}(\bm{Y})}{\floor{2\gamma}}{2\gamma - \floor{2\gamma}} & = \hnorm{\int_{0}^{1}\frac{d}{ds}\mc{R}((1 - s)\bm{X} + s\bm{Y})}{\floor{2\gamma}}{2\gamma - \floor{2\gamma}}ds\\
 &\leq \int_{0}^{1}\hnorm{\p_{\bm{X}}\mc{R}((1 - s)\bm{X} + s\bm{Y})}{\floor{2\gamma}}{2\gamma - \floor{2\gamma}}ds \hnorm{\bm{X} - \bm{Y}}{1}{\gamma}.
\end{align*}

We can prove that $\p_{\bm{X}}\mc{R} = \p_{\bm{X}}\mc{R}_{C} + \p_{\bm{X}}\mc{R}_{T}$ is bounded on $C^{\floor{2\gamma}, 2\gamma - \floor{2\gamma}}$ by showing that both $\p_{\bm{X}}\mc{R}_C$ and  $\p_{\bm{X}}\mc{R}_T$ are. Let $\bm{X}, \bm{Z}\in \mc{B}$. After some brief calculations we find

\begin{align*}\label{RemainderCot}
\p_{\bm{X}}\mc{R}_C(\bm{X})\bm{Z} &= \frac{1}{4\pi}\int_{\mbs}\left(\frac{\Delta\bm{X}\cdot\p_{\theta'}\bm{X}'}{|\Delta\bm{X}|^2} - \frac{1}{2}\cot\left(\frac{\theta - \theta'}{2}\right) \right)\left(\bm{V}_1 - \p_{\theta'}\bm{Z}' \right)d\theta'\\
&+ \frac{1}{4\pi}\int_{\mbs}\left(\frac{|\Delta\bm{X}|^2\Delta\bm{X}\cdot\p_{\theta'}\bm{Z}' + |\Delta\bm{X}|^2\Delta\bm{Z}\cdot\p_{\theta'}\bm{X}' - 2\Delta\bm{X}\cdot\p_{\theta'}\bm{X}'\Delta\bm{X}\cdot\Delta\bm{Z}}{|\Delta\bm{X}|^4} \right)\left(\bm{V}_2 - \p_{\theta'}\bm{X}' \right)d\theta'\\
\end{align*}
for any constant vectors $\bm{V}_1$ and $\bm{V}_2$. We can rewrite this as

\begin{align*}
\p_{\bm{X}}\mc{R}_C(\bm{X})\bm{Z} &= \frac{1}{4\pi}\int_{\mbs}\left(\frac{\Delta\bm{X}\cdot\p_{\theta'}\bm{X}'}{|\Delta\bm{X}|^2} - \frac{1}{2}\cot\left(\frac{\theta - \theta'}{2}\right) \right)\left(\bm{V}_1 - \p_{\theta'}\bm{Z}' \right)d\theta'\\
&+ \frac{1}{4\pi}\int_{\mbs}\left(\frac{|\Delta\bm{X}|^2\Delta\bm{X}\cdot\p_{\theta'}\bm{Z}' - \Delta\bm{X}\cdot\p_{\theta'}\bm{X}'\Delta\bm{X}\cdot\Delta\bm{Z}}{|\Delta\bm{X}|^4} \right)\left(\bm{V}_2 - \p_{\theta'}\bm{X}' \right)d\theta'\\
&+ \frac{1}{4\pi}\int_{\mbs}\left(\frac{|\Delta\bm{X}|^2\Delta\bm{Z}\cdot\p_{\theta'}\bm{X}' - \Delta\bm{X}\cdot\p_{\theta'}\bm{X}'\Delta\bm{X}\cdot\Delta\bm{Z}}{|\Delta\bm{X}|^4} \right)\left(\bm{V}_2 - \p_{\theta'}\bm{X}' \right)d\theta'\\
&=: R_1 + R_2^1 + R_2^2.
\end{align*}
From Proposition \ref{C0cotgain} we know that 

\begin{equation*}
\hnorm{R_1}{\floor{2\gamma}}{2\gamma - \floor{2\gamma}} \leq C \frac{\hdotnorm{\bm{X}}{1}{\gamma}^3}{\starnorm{\bm{X}}^3} \hdotnorm{\bm{Z}}{1}{\gamma},
\end{equation*}
So we turn to $R_2^i$ terms. Since the choice of $\bm{V}_2$ must be the same for both, we must treat them simultaneously. However, the bounds for $R_2^1$ also hold for $R_2^2$ and are found using the exact same methods and choice of $\bm{V}_2$. With this in mind, we calculate bounds explicitly for $R_2^1$ only. Note that we can rewrite the kernel of $R_2^1$ as

\begin{equation}\label{R_11_lemmaform}
\begin{aligned}
\frac{|\Delta\bm{X}|^2\Delta\bm{X}\cdot\p_{\theta'}\bm{Z}' - \Delta\bm{X}\cdot\p_{\theta'}\bm{X}'\Delta\bm{X}\cdot\Delta\bm{Z}}{|\Delta\bm{X}|^4} &= \frac{|\Delta\bm{X}|^2\Delta\bm{X}\cdot\left(\p_{\theta'}\bm{Z}' - \frac{\Delta\bm{Z}}{\theta - \theta'} \right)}{|\Delta\bm{X}|^4}\\
&+ \frac{\Delta\bm{X}\cdot\left(\frac{\Delta\bm{X}}{\theta - \theta'} - \p_{\theta'}\bm{X}'\right)\Delta\bm{X}\cdot\Delta\bm{Z}}{|\Delta\bm{X}|^4},
\end{aligned}
\end{equation}
which is of the form needed to use Lemmas \ref{GPlog1} and \ref{GPlog2}. Using them, we find

\begin{align*}
\left|\frac{|\Delta\bm{X}|^2\Delta\bm{X}\cdot\left(\p_{\theta'}\bm{Z}' - \frac{\Delta\bm{Z}}{\theta - \theta'} \right)}{|\Delta\bm{X}|^4} \right| &\leq C\frac{\hdotnorm{\bm{Z}}{1}{\gamma}}{\starnorm{\bm{X}}}|\theta - \theta'|^{\gamma}\\
\left| \frac{\Delta\bm{X}\cdot\left(\frac{\Delta\bm{X}}{\theta - \theta'} - \p_{\theta'}\bm{X}'\right)\Delta\bm{X}\cdot\Delta\bm{Z}}{|\Delta\bm{X}|^4}\right| &\leq C\frac{\hdotnorm{\bm{Z}}{1}{\gamma}}{\starnorm{\bm{X}}}|\theta - \theta'|^{\gamma}.
\end{align*}
Combining them yields

\begin{align}\label{R_11_bound}
\left|\frac{|\Delta\bm{X}|^2\Delta\bm{X}\cdot\p_{\theta'}\bm{Z}' - \Delta\bm{X}\cdot\p_{\theta'}\bm{X}'\Delta\bm{X}\cdot\Delta\bm{Z}}{|\Delta\bm{X}|^4} \right| &\leq C\frac{\hdotnorm{\bm{X}}{1}{\gamma}\hdotnorm{\bm{Z}}{1}{\gamma}}{\starnorm{\bm{X}}^2}|\theta - \theta'|^{\gamma}.
\end{align}
Taking $\bm{V}_2 = \p_{\theta}\bm{X}$ gives

\begin{align*}
R_{2}^1 (\theta) &= \frac{1}{4\pi}\int_{\mbs}  \left(\frac{|\Delta\bm{X}|^2\Delta\bm{X}\cdot\left(\p_{\theta'}\bm{Z}' - \frac{\Delta\bm{Z}}{\theta - \theta'} \right)}{|\Delta\bm{X}|^4}\right) (\p_{\theta}\bm{X} - \p_{\theta'}\bm{X}')d\theta'\\
&+ \frac{1}{4\pi}\int_{\mbs}  \left(\frac{\Delta\bm{X}\cdot\left(\frac{\Delta\bm{X}}{\theta - \theta'} - \p_{\theta'}\bm{X}'\right)\Delta\bm{X}\cdot\Delta\bm{Z}}{|\Delta\bm{X}|^4}\right) (\p_{\theta}\bm{X} - \p_{\theta'}\bm{X}')d\theta'.
\end{align*}
Using (\ref{R_11_bound}),

\begin{align*}
|R_{2}^{1}(\theta)| &\leq C \frac{\hdotnorm{\bm{X}}{1}{\gamma}\hdotnorm{\bm{Z}}{1}{\gamma}}{\starnorm{\bm{X}}^2} \int_{\mbs}|\theta - \theta'|^{\gamma - 1}|\p_{\theta}\bm{X} - \p_{\theta'}\bm{X}'|d\theta'\\
&\leq C\frac{\hdotnorm{\bm{X}}{1}{\gamma}^2\hdotnorm{\bm{Z}}{1}{\gamma}}{\starnorm{\bm{X}}^2}\int_{\mbs}|\theta - \theta'|^{2\gamma - 1}d\theta' \leq C\frac{\hdotnorm{\bm{X}}{1}{\gamma}^2\hdotnorm{\bm{Z}}{1}{\gamma}}{\starnorm{\bm{X}}^2}.
\end{align*}

Let us now work on the case when $\gamma < 1/2$. Let $0 < h < 2\pi$. Using our usual sets $\mc{I}_s$ and $\mc{I}_f$ and letting $\bm{V}_1 = \p_{\theta}\bm{X}$,

\begin{align*}
\diff_{h}R_{2}^{1} &= \frac{1}{4\pi}\int_{\mbs} \trl_{h}\left( \frac{|\Delta\bm{X}|^2\Delta\bm{X}\cdot\p_{\theta'}\bm{Z}' - \Delta\bm{X}\cdot\p_{\theta'}\bm{X}'\Delta\bm{X}\cdot\Delta\bm{Z}}{|\Delta\bm{X}|^4}\right)(\p_{\theta}\bm{X} - \p_{\theta'}\bm{X}')d\theta'\\
&+ \frac{1}{4\pi}\int_{\mbs} \left( \frac{|\Delta\bm{X}|^2\Delta\bm{X}\cdot\p_{\theta'}\bm{Z}' - \Delta\bm{X}\cdot\p_{\theta'}\bm{X}'\Delta\bm{X}\cdot\Delta\bm{Z}}{|\Delta\bm{X}|^4}\right)(\p_{\theta}\bm{X} - \p_{\theta'}\bm{X}')d\theta'\\
&= \frac{1}{4\pi}\int_{\mc{I}_s} \trl_{h}\left( \frac{|\Delta\bm{X}|^2\Delta\bm{X}\cdot\p_{\theta'}\bm{Z}' - \Delta\bm{X}\cdot\p_{\theta'}\bm{X}'\Delta\bm{X}\cdot\Delta\bm{Z}}{|\Delta\bm{X}|^4}\right)(\p_{\theta}\bm{X} - \p_{\theta'}\bm{X}')d\theta'\\
&+ \frac{1}{4\pi}\int_{\mc{I}_s} \left( \frac{|\Delta\bm{X}|^2\Delta\bm{X}\cdot\p_{\theta'}\bm{Z}' - \Delta\bm{X}\cdot\p_{\theta'}\bm{X}'\Delta\bm{X}\cdot\Delta\bm{Z}}{|\Delta\bm{X}|^4}\right)(\p_{\theta}\bm{X} - \p_{\theta'}\bm{X}')d\theta'\\
&+ \frac{1}{4\pi}\int_{\mc{I}_f}\diff_{h}\left( \frac{|\Delta\bm{X}|^2\Delta\bm{X}\cdot\p_{\theta'}\bm{Z}' - \Delta\bm{X}\cdot\p_{\theta'}\bm{X}'\Delta\bm{X}\cdot\Delta\bm{Z}}{|\Delta\bm{X}|^4}\right)(\p_{\theta}\bm{X} - \p_{\theta'}\bm{X}')d\theta'\\
&=: I_1 + I_2 + I_3.
\end{align*}
For term $I_1$, we may reuse (\ref{R_11_bound}) and find

\begin{align*}
\left|I_1 \right| &\leq C \frac{\hdotnorm{\bm{X}}{1}{\gamma}\hdotnorm{\bm{Z}}{1}{\gamma}}{\starnorm{\bm{X}}^2} \int_{\mc{I}_s}|\theta + h - \theta'|^{1 - \gamma}|\p_{\theta}\bm{X} - \p_{\theta'}\bm{X}'|d\theta'\\
&\leq \frac{\hdotnorm{\bm{X}}{1}{\gamma}^2 \hdotnorm{\bm{Z}}{1}{\gamma}}{\starnorm{\bm{X}}^2} \int_{\mc{I}_s}|\theta + h - \theta'|^{1 - \gamma}|\theta - \theta'|^{\gamma}d\theta'\\
&\leq C h^{\gamma} \frac{\hdotnorm{\bm{X}}{1}{\gamma}^2 \hdotnorm{\bm{Z}}{1}{\gamma}}{\starnorm{\bm{X}}^2} \int_{\mc{I}_s}|\theta + h - \theta'|^{1 - \gamma}d\theta' \leq Ch^{2\gamma} \frac{\hdotnorm{\bm{X}}{1}{\gamma}^2 \hdotnorm{\bm{Z}}{1}{\gamma}}{\starnorm{\bm{X}}^2} .
\end{align*}
The same bound holds for term $I_2$. For term $I_3$, since we are on $\mc{I}_f$ we may apply estimate (\ref{logest3}) to (\ref{R_11_lemmaform}) and get

\begin{align*}
\left|\diff_h \left(\frac{|\Delta\bm{X}|^2\Delta\bm{X}\cdot\p_{\theta'}\bm{Z}' - \Delta\bm{X}\cdot\p_{\theta'}\bm{X}'\Delta\bm{X}\cdot\Delta\bm{Z}}{|\Delta\bm{X}|^4}\right) \right| &\leq Ch\frac{\hdotnorm{\bm{Z}}{1}{\gamma}\hdotnorm{\bm{X}}{1}{\gamma}^2}{\starnorm{\bm{X}}^3}|\theta - \theta'|^{\gamma - 2}.
\end{align*}
Thus,

\begin{align*}
\left| I_3\right| &\leq Ch\frac{\hdotnorm{\bm{Z}}{1}{\gamma}\hdotnorm{\bm{X}}{1}{\gamma}^2}{\starnorm{\bm{X}}^3}\int_{\mc{I}_f}|\theta - \theta'|^{\gamma - 2}|\p_{\theta}\bm{X}-\p_{\theta'}\bm{X}'|d\theta'\\
&\leq Ch\frac{\hdotnorm{\bm{Z}}{1}{\gamma}\hdotnorm{\bm{X}}{1}{\gamma}^3}{\starnorm{\bm{X}}^3}\int_{\mc{I}_f}|\theta - \theta'|^{2\gamma - 2}d\theta' \leq Ch^{2\gamma}\frac{\hdotnorm{\bm{Z}}{1}{\gamma}\hdotnorm{\bm{X}}{1}{\gamma}^3}{\starnorm{\bm{X}}^3},
\end{align*}
as desired. 

For the case of $\gamma > 1/2$, we have

\begin{align*}
\p_{\theta}R_{2}^{1}(\theta) &= \frac{1}{4\pi}\int_{\mbs}\p_{\theta}\left( \frac{|\Delta\bm{X}|^2\Delta\bm{X}\cdot\p_{\theta'}\bm{Z}' - \Delta\bm{X}\cdot\p_{\theta'}\bm{X}'\Delta\bm{X}\cdot\Delta\bm{Z}}{|\Delta\bm{X}|^4}\right)(\p_{\theta}\bm{X} - \p_{\theta'}\bm{X}')d\theta', 
\end{align*}
where now we are forced to set vector $\bm{V}_2 = \p_{\theta}\bm{X}(\theta)$. See the Proposition \ref{C0cotgain} following equation \eqref{Aopdef} for details. We now need only show that $\p_{\theta}R_{2}^1 \in C^{0, 2\gamma - 1}$. Note that by linearity, 

\begin{align*}
\p_{\theta}\left( \frac{|\Delta\bm{X}|^2\Delta\bm{X}\cdot\p_{\theta'}\bm{Z}' - \Delta\bm{X}\cdot\p_{\theta'}\bm{X}'\Delta\bm{X}\cdot\Delta\bm{Z}}{|\Delta\bm{X}|^4}\right) &= \p_{\theta}\left(\frac{|\Delta\bm{X}|^2\Delta\bm{X}\cdot\left(\p_{\theta'}\bm{Z}' - \frac{\Delta\bm{Z}}{\theta - \theta'} \right)}{|\Delta\bm{X}|^4}\right)\\
&+ \p_{\theta}\left(\frac{\Delta\bm{X}\cdot\left(\frac{\Delta\bm{X}}{\theta - \theta'} - \p_{\theta'}\bm{X}'\right)\Delta\bm{X}\cdot\Delta\bm{Z}}{|\Delta\bm{X}|^4}\right).
\end{align*}
Thus, (\ref{logest2}) implies

\begin{align}\label{R_21_deriv_bound}
\left| \p_{\theta}\left( \frac{|\Delta\bm{X}|^2\Delta\bm{X}\cdot\p_{\theta'}\bm{Z}' - \Delta\bm{X}\cdot\p_{\theta'}\bm{X}'\Delta\bm{X}\cdot\Delta\bm{Z}}{|\Delta\bm{X}|^4}\right)\right| &\leq C\frac{\hdotnorm{\bm{X}}{1}{\gamma}^2 \hdotnorm{\bm{Z}}{1}{\gamma}}{\starnorm{\bm{X}}^3}|\theta - \theta'|^{\gamma - 2}.
\end{align}
By periodicity,

\begin{align*}
|\p_{\theta}R_{2}^{1}(\theta)| &\leq C\frac{\hdotnorm{\bm{X}}{1}{\gamma}^2 \hdotnorm{\bm{Z}}{1}{\gamma}}{\starnorm{\bm{X}}^3}\int_{\mbs}|\theta - \theta'|^{\gamma - 2}|\p_{\theta}\bm{X} - \p_{\theta'}\bm{X}'|d\theta'\\
&\leq C\frac{\hdotnorm{\bm{X}}{1}{\gamma}^3 \hdotnorm{\bm{Z}}{1}{\gamma}}{\starnorm{\bm{X}}^3}\int_{\mbs}|\theta - \theta'|^{2\gamma - 2}d\theta'\\
&\leq C\frac{\hdotnorm{\bm{X}}{1}{\gamma}^3 \hdotnorm{\bm{Z}}{1}{\gamma}}{\starnorm{\bm{X}}^3}
\end{align*}
since $2\gamma - 1 > 0$. To finish the proof, without loss of generality, let $0 < h < 2\pi$. As usual, we break our domain of integration into $\mc{I}_s$ and $\mc{I}_f$. Then,

\begin{align*}
\diff_h R_{2}^{1} &= \frac{1}{4\pi}\int_{\mbs}\trl_{h}\left(\p_{\theta}\left( \frac{|\Delta\bm{X}|^2\Delta\bm{X}\cdot\p_{\theta'}\bm{Z}' - \Delta\bm{X}\cdot\p_{\theta'}\bm{X}'\Delta\bm{X}\cdot\Delta\bm{Z}}{|\Delta\bm{X}|^4}\right)\right)(\p_{\theta}\bm{X}(\theta + h) - \p_{\theta'}\bm{X}')d\theta'\\
&- \frac{1}{4\pi}\int_{\mbs}\left(\p_{\theta}\left( \frac{|\Delta\bm{X}|^2\Delta\bm{X}\cdot\p_{\theta'}\bm{Z}' - \Delta\bm{X}\cdot\p_{\theta'}\bm{X}'\Delta\bm{X}\cdot\Delta\bm{Z}}{|\Delta\bm{X}|^4}\right)\right)(\p_{\theta}\bm{X}(\theta) - \p_{\theta'}\bm{X}')d\theta'\\
&= \frac{1}{4\pi}\int_{\mc{I}_s}\trl_{h}\left(\p_{\theta}\left( \frac{|\Delta\bm{X}|^2\Delta\bm{X}\cdot\p_{\theta'}\bm{Z}' - \Delta\bm{X}\cdot\p_{\theta'}\bm{X}'\Delta\bm{X}\cdot\Delta\bm{Z}}{|\Delta\bm{X}|^4}\right)\right)(\p_{\theta}\bm{X}(\theta + h) - \p_{\theta'}\bm{X}')d\theta'\\
&- \frac{1}{4\pi}\int_{\mc{I}_s}\left(\p_{\theta}\left( \frac{|\Delta\bm{X}|^2\Delta\bm{X}\cdot\p_{\theta'}\bm{Z}' - \Delta\bm{X}\cdot\p_{\theta'}\bm{X}'\Delta\bm{X}\cdot\Delta\bm{Z}}{|\Delta\bm{X}|^4}\right)\right)(\p_{\theta}\bm{X}(\theta) - \p_{\theta'}\bm{X}')d\theta'\\
&+ \frac{1}{4\pi}\int_{\mc{I}_f}\diff_{h}\left(\p_{\theta}\left( \frac{|\Delta\bm{X}|^2\Delta\bm{X}\cdot\p_{\theta'}\bm{Z}' - \Delta\bm{X}\cdot\p_{\theta'}\bm{X}'\Delta\bm{X}\cdot\Delta\bm{Z}}{|\Delta\bm{X}|^4}\right)\right)(\p_{\theta}\bm{X}(\theta + h) - \p_{\theta'}\bm{X}')d\theta'\\
&+ \frac{1}{4\pi}\int_{\mc{I}_f}\left(\p_{\theta}\left( \frac{|\Delta\bm{X}|^2\Delta\bm{X}\cdot\p_{\theta'}\bm{Z}' - \Delta\bm{X}\cdot\p_{\theta'}\bm{X}'\Delta\bm{X}\cdot\Delta\bm{Z}}{|\Delta\bm{X}|^4}\right)\right)(\p_{\theta}\bm{X}(\theta) - \p_{\theta}\bm{X}(\theta))d\theta'\\
&=: II_1 + II_2 + II_3 + II_4.
\end{align*}
For term $II_1$, we make use of (\ref{R_21_deriv_bound}),

\begin{align*}
|II_1| &\leq C\frac{\hdotnorm{\bm{X}}{1}{\gamma}^2 \hdotnorm{\bm{Z}}{1}{\gamma}}{\starnorm{\bm{X}}^3}\int_{\mc{I}_s}|\theta + h - \theta'|^{\gamma - 2}|\p_{\theta}\bm{X}(\theta + h) - \p_{\theta'}\bm{X}'|d\theta'\\
&\leq C\frac{\hdotnorm{\bm{X}}{1}{\gamma}^3 \hdotnorm{\bm{Z}}{1}{\gamma}}{\starnorm{\bm{X}}^3}\int_{\mc{I}_s}|\theta + h - \theta'|^{2\gamma - 2}d\theta'\\
&\leq Ch^{2\gamma - 1}\frac{\hdotnorm{\bm{X}}{1}{\gamma}^3 \hdotnorm{\bm{Z}}{1}{\gamma}}{\starnorm{\bm{X}}^3},
\end{align*}
where we have used $2\gamma - 1 > 0$. The same bound holds for term $II_2$. For term $II_3$ we make use of linearity and estimate (\ref{logest4}) to find

\begin{align}
\left|\diff_{h}\p_{\theta}\left( \frac{|\Delta\bm{X}|^2\Delta\bm{X}\cdot\p_{\theta'}\bm{Z}' - \Delta\bm{X}\cdot\p_{\theta'}\bm{X}'\Delta\bm{X}\cdot\Delta\bm{Z}}{|\Delta\bm{X}|^4}\right) \right| &\leq C\frac{\hdotnorm{\bm{X}}{1}{\gamma}^3 \hdotnorm{\bm{Z}}{1}{\gamma}}{\starnorm{\bm{X}}^4}\left(h^{\gamma}|\theta - \theta'|^{- 2} + h|\theta - \theta'|^{\gamma - 3}\right).
\end{align}
Using this,

\begin{align*}
|II_3| &\leq C\frac{\hdotnorm{\bm{X}}{1}{\gamma}^3 \hdotnorm{\bm{Z}}{1}{\gamma}}{\starnorm{\bm{X}}^4}\int_{\mc{I}_f}\left(h^{\gamma}|\theta - \theta'|^{- 2} + h|\theta - \theta'|^{\gamma - 3}\right)|\p_{\theta}\bm{X}(\theta + h) - \p_{\theta'}\bm{X}'|d\theta'\\
&\leq C\frac{\hdotnorm{\bm{X}}{1}{\gamma}^4 \hdotnorm{\bm{Z}}{1}{\gamma}}{\starnorm{\bm{X}}^4}\int_{\mc{I}_f}\left(h^{\gamma}|\theta - \theta'|^{- 2} + h|\theta - \theta'|^{\gamma - 3}\right)|\theta + h - \theta'|^{\gamma}d\theta'\\
&\leq C\frac{\hdotnorm{\bm{X}}{1}{\gamma}^4 \hdotnorm{\bm{Z}}{1}{\gamma}}{\starnorm{\bm{X}}^4}\int_{\mc{I}_f}\left(h^{\gamma}|\theta - \theta'|^{\gamma - 2} + h|\theta - \theta'|^{2\gamma - 3}\right)d\theta'\\
&\leq C\frac{\hdotnorm{\bm{X}}{1}{\gamma}^4 \hdotnorm{\bm{Z}}{1}{\gamma}}{\starnorm{\bm{X}}^4}\left(h^{\gamma}(3/2h)^{2\gamma - 1} + h(3/2h)^{2\gamma - 2}\right) = Ch^{2\gamma - 1}\frac{\hdotnorm{\bm{X}}{1}{\gamma}^4 \hdotnorm{\bm{Z}}{1}{\gamma}}{\starnorm{\bm{X}}^4}.
\end{align*}
Finally, via (\ref{R_21_deriv_bound}), for term $II_4$ 

\begin{align*}
|II_4| &\leq C\frac{\hdotnorm{\bm{X}}{1}{\gamma}^2 \hdotnorm{\bm{Z}}{1}{\gamma}}{\starnorm{\bm{X}}^3}\int_{\mc{I}_f}|\theta - \theta'|^{\gamma - 2}|\p_{\theta}\bm{X}(\theta + h) - \p_{\theta}\bm{X}(\theta)|d\theta'\\
&\leq Ch^{\gamma}\frac{\hdotnorm{\bm{X}}{1}{\gamma}^3 \hdotnorm{\bm{Z}}{1}{\gamma}}{\starnorm{\bm{X}}^3}\int_{\mc{I}_f}|\theta - \theta'|^{\gamma - 2}d\theta'\\
&\leq Ch^{\gamma}\frac{\hdotnorm{\bm{X}}{1}{\gamma}^3 \hdotnorm{\bm{Z}}{1}{\gamma}}{\starnorm{\bm{X}}^3}\left((3/2h)^{\gamma - 1} - (2\pi - 1/2h)^{\gamma - 1}\right)\\
&\leq Ch^{2\gamma - 1}\frac{\hdotnorm{\bm{X}}{1}{\gamma}^3 \hdotnorm{\bm{Z}}{1}{\gamma}}{\starnorm{\bm{X}}^3}.
\end{align*}
Thus, if $\gamma > 1/2$, $R_{2}^1 \in C^{1, 2\gamma - 1}$. 

Linearizing $\mc{R}_{T}$ gives

\begin{align*}
\p_{\bm{X}}\mc{R}_T (\bm{X})\bm{Z} &= \frac{1}{4\pi}\int_{\mbs} \frac{(|\Delta\bm{X}|^2\p_{\theta'}\bm{X}' - \Delta\bm{X}\cdot\p_{\theta'}\bm{X}'\Delta\bm{X})\otimes\Delta\bm{X}}{|\Delta\bm{X}|^4}(\p_{\theta'}\bm{Z}' - \bm{V}_3)d\theta'\\
&+ \frac{1}{4\pi}\int_{\mbs} \frac{\Delta\bm{X}\otimes(|\Delta\bm{X}|^2\p_{\theta'}\bm{X}' - \Delta\bm{X}\cdot\p_{\theta'}\bm{X}'\Delta\bm{X})}{|\Delta\bm{X}|^4}(\p_{\theta'}\bm{Z}' - \bm{V}_3)d\theta'\\
&+ \frac{1}{4\pi}\int_{\mbs}\frac{\left(\Delta\bm{X}\cdot\Delta\bm{Z}\p_{\theta'}\bm{X}' - \Delta\bm{Z}\cdot\p_{\theta'}\bm{X}'\Delta\bm{X}\right)\otimes\Delta\bm{X}}{|\Delta\bm{X}|^4}(\p_{\theta'}\bm{X}' - \bm{V}_4)d\theta'\\
&+ \frac{1}{4\pi}\int_{\mbs}\frac{\left(\Delta\bm{X}\cdot\Delta\bm{Z}\p_{\theta'}\bm{X}' - \Delta\bm{X}\cdot\p_{\theta'}\bm{Z}'\Delta\bm{X}\right)\otimes\Delta\bm{X}}{|\Delta\bm{X}|^4}(\p_{\theta'}\bm{X}' - \bm{V}_4)d\theta'\\
&+ \frac{1}{4\pi}\int_{\mbs}\frac{\left(|\Delta\bm{X}|^2\p_{\theta'}\bm{Z}' - \Delta\bm{X}\cdot\p_{\theta'}\bm{X}'\Delta\bm{Z}\right)\otimes\Delta\bm{X}}{|\Delta\bm{X}|^4}(\p_{\theta'}\bm{X}' - \bm{V}_4)d\theta'\\
&+ \frac{1}{4\pi}\int_{\mbs}\frac{\Delta\bm{X}\cdot\Delta\bm{Z}\left(|\Delta\bm{X}|^2\p_{\theta'}\bm{X}' - \Delta\bm{X}\cdot\p_{\theta'}\bm{X}'\Delta\bm{X}\right)\otimes\Delta\bm{X}}{|\Delta\bm{X}|^6}(\p_{\theta'}\bm{X}' - \bm{V}_4)d\theta'\\
&+ \frac{1}{4\pi}\int_{\mbs}\frac{\Delta\bm{X}\otimes\left(\Delta\bm{X}\cdot\Delta\bm{Z}\p_{\theta'}\bm{X}' - \Delta\bm{Z}\cdot\p_{\theta'}\bm{X}'\Delta\bm{X}\right)}{|\Delta\bm{X}|^4}(\p_{\theta'}\bm{X}' - \bm{V}_4)d\theta'\\
&+ \frac{1}{4\pi}\int_{\mbs}\frac{\Delta\bm{X}\otimes\left(\Delta\bm{X}\cdot\Delta\bm{Z}\p_{\theta'}\bm{X}' - \Delta\bm{X}\cdot\p_{\theta'}\bm{Z}'\Delta\bm{X}\right)}{|\Delta\bm{X}|^4}(\p_{\theta'}\bm{X}' - \bm{V}_4)d\theta'\\
&+ \frac{1}{4\pi}\int_{\mbs}\frac{\Delta\bm{X}\otimes\left(|\Delta\bm{X}|^2\p_{\theta'}\bm{Z}' - \Delta\bm{X}\cdot\p_{\theta'}\bm{X}'\Delta\bm{Z}\right)}{|\Delta\bm{X}|^4}(\p_{\theta'}\bm{X}' - \bm{V}_4)d\theta'\\
&+\int_{\mbs}\frac{\Delta\bm{X}\otimes\Delta\bm{X}\cdot\Delta\bm{Z}\left(|\Delta\bm{X}|^2\p_{\theta'}\bm{X}' - \Delta\bm{X}\cdot\p_{\theta'}\bm{X}'\Delta\bm{X}\right)}{|\Delta\bm{X}|^6}(\p_{\theta'}\bm{X}' - \bm{V}_4)d\theta'\\
&+ \frac{1}{4\pi}\int_{\mbs} \frac{(|\Delta\bm{X}|^2\p_{\theta'}\bm{X}' - \Delta\bm{X}\cdot\p_{\theta'}\bm{X}'\Delta\bm{X})\otimes\Delta\bm{Z}}{|\Delta\bm{X}|^4}(\p_{\theta'}\bm{X}' - \bm{V}_4)d\theta'\\
&+ \frac{1}{4\pi}\int_{\mbs} \frac{\Delta\bm{Z}\otimes(|\Delta\bm{X}|^2\p_{\theta'}\bm{X}' - \Delta\bm{X}\cdot\p_{\theta'}\bm{X}'\Delta\bm{X})}{|\Delta\bm{X}|^4}(\p_{\theta'}\bm{X}' - \bm{V}_4)d\theta'.
\end{align*}
For any constant vectors $\bm{V}_3, \bm{V}_4$. We have rearranged some terms so that it is easy to see that we can use the same arguments as we did for $\mc{R}_C$ and get the same result. For the sake of brevity, we will not repeat them here. We can combine all of the terms and bounds and find

\begin{align*}
\hnorm{\p_{\bm{X}}\mc{R}(\bm{X})\bm{Z}}{\floor{2\gamma}}{2\gamma - \floor{2\gamma}} &\leq C\frac{\hdotnorm{\bm{X}}{1}{\gamma}^4\hdotnorm{\bm{Z}}{1}{\gamma}}{\starnorm{\bm{X}}^4}.
\end{align*}

\end{poof}

Consider a set of functions with
\begin{align*}
\hnorm{\bm{X} - \bm{X}_0}{1}{\gamma} \leq \frac{1}{2}\starnorm{\bm{X}_0}.
\end{align*}
In this set, we have
\begin{align*}
\starnorm{\bm{X}} &\geq \starnorm{\bm{X}_0} - \hnorm{\bm{X} - \bm{X}_0}{1}{\gamma} \geq \frac{1}{2}\starnorm{\bm{X}_0}=m,\\
\chnorm{\bm{X}}{1,\gamma}&\leq \chnorm{\bm{X}_0}{1,\gamma}+\frac{1}{2}\starnorm{\bm{X}_0}=M.
\end{align*}
With this in mind, define the set
\begin{align}\label{setBT}
\mc{B}_{T} = \left\lbrace \bm{X} \in C([0, T]; C^{1, \gamma}(\mbs)) : \norm{\bm{X} -\bm{X}_0}_{C([0, T]; C^{1, \gamma})} \leq \frac{1}{2}\starnorm{\bm{X}_0}  \right\rbrace,
\end{align}
where we have abused notation slightly to write $\bm{X}_0$ as the function that is constant in time taking value $\bm{X}_0$. 
Note that $\mc{B}_T$ is a convex set in $O^{M,m}_T$.
We are now in position to prove the existence of a mild solution to the Peskin system. 

\begin{proposition}\label{LocalExist}
There exists some time $T > 0$ such that $S(\bm{X}, t; \bm{X}_0)$ forms a contraction on $\mc{B}_T$.
\end{proposition}

\begin{poof}
We first show that there is some $T_1>0$ for which $S$ maps $B_{T_1}$ to itself. For $\bm{X}\in \mc{B}_{T_1}$ and $\gamma\neq 1/2$,
\begin{equation}\label{SXtX0-X0}
\begin{split}
\hnorm{S(\bm{X}, t; \bm{X}_0)-\bm{X}_0}{1}{\gamma} &\leq \hnorm{e^{t\Lambda}\bm{X}_0-\bm{X}_0}{1}{\gamma} + \int_{0}^{t}\hnorm{e^{(t - s)\Lambda}\mc{R}(\bm{X}(s))}{1}{\gamma}ds\\
&\leq \hnorm{(e^{t\Lambda} - 1)\bm{X}_0}{1}{\gamma} + C\int_{0}^{t}(t - s)^{\gamma - 1}\hflnorm{\mc{R}(\bm{X}(s))}{2\gamma}ds\\
&\leq \hnorm{(e^{t\Lambda} - 1)\bm{X}_0}{1}{\gamma} + C\frac{M^5}{m^4}\int_{0}^{t}(t - s)^{\gamma - 1}ds \leq \chnorm{(e^{t\Lambda} - 1)\bm{X}_0}{1,\gamma} + Ct^{\gamma}\frac{M^5}{m^4},
\end{split}
\end{equation}
where we used equation \eqref{mainSGest} between the first and second line and Propostions \ref{C0cotgain} and \ref{C0tengain} between lines two and three. 
Since $\bm{X}_0 \in h^{1, \gamma}(\mbs)$, by Proposition \ref{semigroup_strong_continuity}, we may take $T_1$ small enough so that
\begin{align*}
\norm{(e^{t\Lambda} - 1)\bm{X}_0}_{C^{1, \gamma}}+ Ct^{\gamma}\frac{M^5}{m^4} \leq \frac{1}{2}\starnorm{\bm{X}_0}
\end{align*}
for all $0\leq t \leq T_1$. 
We now show that $S$ forms a contraction on $\mc{B}_{T_2}$ for some $T_2>0$. Let $\bm{X}$, $\bm{Y}\in \mc{B}_t$. Then,
\begin{align*}
\hnorm{S(\bm{X}, t; \bm{X}_0) - S(\bm{Y}, t; \bm{X}_0 )}{1}{\gamma} &\leq \int_{0}^{t}\hnorm{e^{(t - s)\Lambda}[\mc{R}(\bm{X}(s)) - \mc{R}(\bm{Y}(s))]}{1}{\gamma}\\
&\leq C\int_{0}^{t}(t - s)^{\gamma - 1}\hflnorm{\mc{R}(\bm{X}(s)) - \mc{R}(\bm{Y}(s))}{2\gamma}ds\\
&\leq C\frac{M^4}{m^4}\norm{\bm{X} - \bm{Y}}_{C([0, t]; C^{1, \gamma})}\int_{0}^{t}(t - s)^{\gamma - 1} ds\\
&\leq Ct^{\gamma}\frac{M^4}{m^4}\norm{\bm{X} - \bm{Y}}_{C([0, t]; C^{1, \gamma})},
\end{align*}
where Proposition \ref{LipRemainder} was used between lines two and three. There exists some time $T_2>0$ such that 
\begin{align*}
Ct^{\gamma}\frac{M^4}{m^4} \leq \frac{1}{2}
\end{align*}
for all $0\leq t \leq T_2$. Taking $T = \min\{T_1, T_2 \}$ gives the desired result. In the case of $\gamma = 1/2$, we use the second statements in Propositions \ref{C0cotgain}, \ref{C0tengain} and \ref{LipRemainder} with the choice of $\alpha = 3\gamma/2$ and the result follows from the same arguments.

\end{poof}

Proposition \ref{LocalExist} gives the local existence of a solution to the Peskin problem for any initial data $\bm{X}_0 \in h^{1, \gamma}(\mbs)$ with $\starnorm{\bm{X}_0} > 0$. Our next result shows that our mild solution is, in fact, a solution to the differential form of the PDE away from the initial time.  In particular, we have 

\begin{lemma} \label{mild_is_strong}
Let  $\bm{X}(t) \in C([0,T]; C^{1,\gamma}(\mbs))$ be a mild solution with initial data $\bm{X}_0$ as in the assumptions in Theorem~\ref{LWPTheorem}, then 
\begin{equation} \label{e:classical}
\p_t \bm{X} = \Lambda \bm{X} + \mathcal{R}(\bm{X})
\end{equation}
for the time interval $(0,t)$.  Furthermore, $\p_t \bm{X} \in C([0,T]; C^{0,\gamma}(\mbs))$.  
\end{lemma} 
\begin{poof}
We use some of the ideas of the proof of Lemma 4.1.6 in \cite{lunardi}. First, we claim that $\int_0^t \bm{X} (s) ds \in C^{1, \gamma}$ and 
\begin{equation} \label{e:avmild}
\bm{X}(t) = \bm{X}_0 + \Lambda \int_0^t \bm{X}(s) ds + \int_0^t \mathcal{R}(\bm{X}(s)) ds
\end{equation}
for all $[0,t]$.  The first fact follows from $\chnorm{\int_0^t \bm{X}(s) ds }{1,\gamma} \leq t \sup_{0\leq s \leq t} \chnorm{\bm{X}(s)}{1,\gamma} \leq  Ct$.

To show \eqref{e:avmild} we note from \eqref{etlamck} that $\chnorm{e^{(s-\sigma) \Lambda} \mathcal{R}(\bm{X}(\sigma))}{1, \gamma} \leq
\chnorm{\mathcal{R}(\bm{X}(\sigma)}{1, \gamma} \leq C(\chnorm{\bm{X}(\sigma)}{1, \gamma}, \starnorm{\bm{X}(\sigma)})$ for $0 \leq \sigma \leq s \leq t$, and so by Fubini,
\begin{equation} \label{e:intmild}
\begin{split}
\int_0^t \bm{X}(s) ds 
& = \int_0^t e^{s \Lambda} \bm{X}_0 ds + \int_0^t \int_0^s e^{(s-\sigma) \Lambda} \mathcal{R}(\bm{X}(\sigma)) d \sigma ds \\
& = \int_0^t e^{s \Lambda} \bm{X}_0 ds + \int_0^t  \int^t_\sigma e^{(s-\sigma) \Lambda} \mathcal{R}(\bm{X}(\sigma)) d s d\sigma.
\end{split}
\end{equation}

From the construction of the $\Lambda$ operator in subsection~\ref{SGestSec}, we can write $ \Lambda e^{t \Lambda}= \p_t (e^{t \Lambda})$.  Therefore
\begin{align*}
 \int_0^t \Lambda \bm{X}(s) ds 
 &  = 
\int_0^t \Lambda e^{s \Lambda} \bm{X}_0 ds + \int_0^t   \int_0^s  \Lambda e^{(s - \sigma ) \Lambda} \mathcal{R}(\bm{X}(\sigma)) d \sigma d s  \\
 &  = 
\int_0^t \Lambda e^{s \Lambda} \bm{X}_0 ds + \int_0^t  \int^t_\sigma \Lambda e^{(s-\sigma) \Lambda} \mathcal{R}(\bm{X}(\sigma)) d s d\sigma  \\
& = 
(e^{t \Lambda} - 1 ) \bm{X}_0 + \int_0^t \left( e^{(t-\sigma) \Lambda} - 1 \right) \mathcal{R}(\bm{X}(\sigma)) d \sigma \\
& = \bm{X}(t) - \bm{X}_0 - \int_0^t \mathcal{R}(\bm{X}(\sigma)) d \sigma 
\end{align*}
and hence \eqref{e:avmild}.

2. Since $\int_0^t \bm{X}(s) ds \in C^{1,\gamma}$ then $\Lambda \int_0^t \bm{X}(s) ds \in C^{0,\gamma}$.  Furthermore, from the local well-posedness theory, $\mathcal{R}(\bm{X}(t)) \in C^{\floor{2\gamma}, 2\gamma - \floor{2\gamma}}$.   We can now define the finite difference, 
\[
{\bm{X}(t + h) - \bm{X}(t) \over h} = { 1\over h} \int_{t}^{t+ h} \Lambda  \bm{X}(s) ds + {1\over h} \int_t^{t+h} \mathcal{R}(\bm{X}(s)) ds.
\]
Thanks to Lipschitz continuity of $\mc{R}$ proved in Proposition \ref{LipRemainder}, $\mathcal{R}(\bm{X})$ is continuous at $t$. Then,
\[
\lim_{t \to 0} {1\over h} \int_t^{t+h} \mathcal{R}(\bm{X}(s)) ds = \mathcal{R}(\bm{X}(t)).
\]
Likewise, $ \lim_{t \to 0} { 1\over h} \int_{t}^{t+ h} \Lambda  \bm{X}(s) ds = \Lambda \bm{X}(t)$.  
Combining the limits yields \eqref{e:classical} with $\p_t \bm{X} \in C([0,T], C^{0,\gamma}(\mbs))$.  
\end{poof}

We also show that a strong solution satisfies the mild form of the equation.  This fact will also be used in the stability analysis.
\begin{lemma} \label{strong_is_mild}
Let  $\bm{X}(t) \in C([0,t]; C^{1,\gamma}(\mbs))\cap C^1([0,T];C^{0,\gamma}(\mbs))$ solve \eqref{SSD} with initial data $\bm{X}_0$, then $\bm{X}(t)$ satisfies
\eqref{mild_soln} on the time interval $[0,T]$. 
\end{lemma} 
\begin{poof}
If we fix $t> 0$ and define $\bm{Z}(s) = e^{(t-s)\Lambda} \bm{X}(s)$ then for any $0 \leq s < t$ we have 
\begin{align*}
\p_s \bm{Z}(s) & = e^{(t-s)\Lambda} \p_s \bm{X}(s) - e^{(t-s)\Lambda} \Lambda \bm{X}(s) \\
& = e^{(t-s)\Lambda} \mathcal{R}(\bm{X}(s)), 
\end{align*}
where $\p_s \bm{X}(s), \Lambda \bm{X}(s) \in C^{0,\gamma}(\mbs)$. Integrating on $[0,t]$ yields 
\[
\bm{X}(t) - e^{t \Lambda} \bm{X}_0 = \int_0^t e^{(t-s)\Lambda} \mathcal{R}(\bm{X}(s)) ds
\]
and \eqref{mild_soln}.
\end{poof}

We are now able to prove Theorem \ref{LWPTheorem}.

\begin{poofof}{Theorem \ref{LWPTheorem}}
We have already proved item \ref{LWPext} in Proposition \ref{LocalExist} and item \ref{LWPmildstrong} in Lemmas \ref{mild_is_strong} and \ref{strong_is_mild}.
We prove item \ref{LWPuniq}.
Suppose we have two solutions $\bm{Y}$ and $\bm{Z}$ in $C([0,T];C^{1,\gamma}(\mbs))$ with the same initial value $\bm{X}_0$.
Define $T_*$ to be:
\begin{equation*}
T_*=\sup \lbrace \tau: 0\leq \tau\leq T, \; \bm{Y}(t)=\bm{Z}(t) \text{ for all } 0\leq t\leq \tau\rbrace
\end{equation*}
We show that $T_*=T$. Suppose otherwise. Then, $T_*<T$. First note that $\bm{Y}(T_*)=\bm{Z}(T_*)\in h^{1,\gamma}(\mbs)$.
If $T_*=0$, this is true by assumption. For $T_*>0$, this follows from:
\begin{equation*}
\norm{\bm{Y}(T_*)}_{C^{1,\gamma'}}\leq \frac{C}{T_*^{\gamma'-\gamma}}\norm{\bm{X}_{0}}_{C^{1,\gamma}}+\int_0^{T_*}\frac{C}{(T_*-s)^{\gamma'-\gamma}}\norm{\bm{Y}(s)}_{C^{1,\gamma}}ds<\infty,
\end{equation*}
where $\gamma<\gamma'<1$. Thus, $\bm{Y}(T_*)\in C^{1,\gamma'}(\mbs)\subset h^{1,\gamma}(\mbs)$. We also have $\starnorm{\bm{Y}(T_*)}>0$ by
the definition of the mild solution.
We may thus consider a mild solution to the Peskin problem starting at $t=T_*$ with initial value $\bm{X}_*=\bm{Y}(T_*)=\bm{Z}(T_*)$.
By the contraction mapping argument of Proposition \ref{LocalExist}, there is a unique mild solution $\bm{W}(t)$ with initial data $\bm{X}_*$ 
for some time $0\leq t\leq \wt{T}_*\leq T-T_*, \wt{T}_*>0$. By the uniqueness of the fixed point of 
the contraction map, we must have $\bm{Y}(t+T_*)=\bm{W}(t)=\bm{Z}(t+T_*)$ for $0\leq t\leq \wt{T}_*$. 
This is a contradiction. 

%
%
%

We next prove item \ref{LWPcont}. Let $\bm{X}(t)$ be a mild solution in $C([0,T];C^{1,\gamma}(\mbs))$ with initial data $\bm{X}_0$.
Take some time $T_1>0$ and consider the set $\wt{\mc{B}}_{T_1}$:
\begin{equation*}
\wt{\mc{B}}_{X(t_*),T_1}= \left\lbrace \bm{X} \in C([0, T_1]; C^{1, \gamma}(\mbs)) : \norm{\bm{X} -\bm{X}(t_*)}_{C([0, T]; C^{1, \gamma})} 
\leq \frac{1}{2}\inf_{0\leq t\leq T}\starnorm{\bm{X}(t)}  \right\rbrace,
\end{equation*}
Set 
\begin{equation*}
m=\frac{1}{2}\inf_{0\leq t\leq T} \starnorm{\bm{X}(t)}, \; 
M=m+\sup_{0\leq t\leq T} \norm{\bm{X}(t)}_{C^{1,\gamma}}.
\end{equation*}
Consider the map $\mc{S}_{\bm{Y}_0}$ taking 
$\bm{Y}(t)$ to $\mc{S}(\bm{Y},t;\bm{Y}_0)$ (see \eqref{SXtX0}). We may estimate, in the same way as in \eqref{SXtX0-X0},
\begin{equation}
\hnorm{S(\bm{X}, t; \bm{X}_0)-\bm{Y}_0}{1}{\gamma}\leq  \chnorm{(e^{t\Lambda} - 1)\bm{X}_0}{1,\gamma} +\chnorm{\bm{X}_0-\bm{Y}_0}{1,\gamma}
+ Ct^{\gamma}\frac{M^5}{m^4},
\end{equation}
Thus, by taking $\norm{\bm{X}_0-\bm{Y}_0}{C^{1,\gamma}}$ and $T_1$ small enough, we see that $\mc{S}_{\bm{Y}_0}$ maps $\mc{B}_{\bm{X}_0,T_1}$ to itself.
That this is a contraction follows in the same way as Proposition \ref{LocalExist}. Thus, there is an $\epsilon>0$ and $T_1>0$ depending on $\bm{X}_0$
such that, for all initial data $\bm{Y}_0$ satisfying $\norm{\bm{Y}_0-\bm{X}_0}_{C^{1,\gamma}}\leq \epsilon$, 
there is a mild solution $\bm{Y}(t)$ to the Peskin problem up to time $T_1$. We also see that 
\begin{align*}
(\bm{X} - \bm{Y})(t) &= e^{\Lambda t}(\bm{X}_0 - \bm{Y}_0) + \int_{0}^{t}e^{\Lambda (t - s)}\left[ \mc{R}(\bm{X})(s) - \mc{R}(\bm{Y})(s)\right]ds\\
\end{align*}
for all $t \in [0,T_1]$. Taking the $C^{1, \gamma}$ norm of this equation and using Proposition \ref{LipRemainder} and semigroup estimate \eqref{mainSGest} gives
\begin{align*}
\hnorm{(\bm{X} - \bm{Y})(t)}{1}{\gamma} &\leq C\hnorm{\bm{X}_0 - \bm{Y}_0}{1}{\gamma} + C\frac{M^{4}}{m^{4}}\int_{0}^{t}(t - s)^{\gamma}\hnorm{(\bm{X} - \bm{Y})(s)}{1}{\gamma}ds.
\end{align*} 
Making use of a generalized Gr\"{o}nwall's lemma from Lemma~7.0.3 of \cite{lunardi} gives
\begin{align*}
\hnorm{(\bm{X} - \bm{Y})(t)}{1}{\gamma} &\leq C\hnorm{\bm{X}_0 - \bm{Y}_0}{1}{\gamma}\leq C\epsilon.
\end{align*}
Since $t\in[0, T]$ was arbitrary, taking the supremum over all $t \in [0, T_1]$ this shows that we have continuity in $C([0,T_1];C^{1,\gamma}(\mbs))$.
Selecting $t_*=T_1$ and repeating this argument with initial data $\bm{Y}_1$ close to $\bm{X}(T_1)$, by possibly reducing the value of $\epsilon$, 
we can extend our result up to $t=T_1$. 
This process can be repeated until we reach $t=T$.
\end{poofof}

\section{Smoothness and equivalence of solutions} \label{s:smoothness}

In this section we first prove higher regularity of our solution in space and time in Subsection~\ref{sect:higher}.  Once the regularity of the boundary integral formulation of our problem has been established, we show that our solution is equivalent to the other formulations of the Peskin problem in Subsection~\ref{sect:classical}. For $n\in \mathbb{N}$ and $\alpha > 0$, we occasionally write $C^{n + \alpha}$ to mean $C^{n + \floor{\alpha}, \alpha - \floor{\alpha}}$ for convenience. 

\subsection{Higher regularity} 
\label{sect:higher}

Given that the symbol of our leading order operator is parabolic, it is natural to ask whether the contours will become smooth as time evolves. In order to establish higher regularity, our method involves converting spatial derivatives in $\theta$ on the integral operator into derivatives in $\theta '$, similar to a method developed in \cite{GuoHallstromSpirn}. Our goal is to show that the nonlinear remainder terms carry the regularity of $\bm{X}$, and improved regularity on $\bm{X}$ follows from an application of  semigroup estimate \eqref{mainSGest}. 

We introduce the notation
\begin{align*}
\bm{\chi}_j&:=\Delta \p_\theta^j \bm{X}=\p_\theta^j \bm{X}-\p_{\theta'}^j\bm{X}', \; \bm{\chi}_j=(\chi_{j,1},\chi_{j,2}),
\end{align*}
for $\bm{X}\in C^{j,\gamma}(\mbs)$ and $\starnorm{\bm{X}} > 0$ and the following sets
\begin{align*}  
\mc{S}_{k} &=\left\{ \sum_{l = 0}^{M} a_l g_l : a_l \in \mathbb{R}, g_l \in \mc{S}_{k}^{1}, M\in \mathbb{N} \right\} , \\
\mc{S}_{k}^1 &= \left\{ \prod_{l=0}^k\paren{\frac{\chi_{l, 1}}{\abs{\bm{\chi}_0}}}^{\alpha_l} \paren{\frac{\chi_{l, 2}}{\abs{\bm{\chi}_0}}}^{\beta_l}, 
\alpha_l, \beta_l \in \mathbb{N}\cup \{0\},\; \sum_{l = 0}^{k}l(\alpha_l + \beta_l)=k \right\}.
\end{align*}
In this section, a term in $\mc{S}_k^1$ will be called a $k$-monomial. 
Note that the sum:
\begin{equation}\label{total_deriv_number}
\sum_{l=0}^k l(\alpha_l+\beta_l)
\end{equation}
is the total number of derivatives in one $k$-monomial. One may thus say that a $k$-monomial 
is a monomial that has exactly $k$ derivatives. Since all $\alpha_l, \beta_l$ are positive, this implies that the 
largest $l$ for which $\alpha_l, \beta_l$ can be non-zero is $l=k$.
Also, note that a $k$-monomial satisfies the assumptions for the function $g$ of 
Lemma \ref{GPL1} and \ref{GPL2} so long as $\bm{X}\in C^{k + 1,\gamma}(\mbs), \; \gamma\in (0,1)$. The number $N$ that appears in the estimates of Lemmas \ref{GPL1} and \ref{GPL2} are given by:
\begin{equation}\label{Ninkmonomial}
N=\sum_{l=1}^k (\alpha_k+\beta_k)\leq k.
\end{equation}

The relevance of the class of $k$-monomials and their linear combinations to the problem at hand is contained in the following lemma.
\begin{lemma}\label{Sk1bounds}
If $k \geq 0$, $f \in \mc{S}_k$, and $\bm{X} \in C^{k +1,\gamma}(\mbs)$ for some  $\gamma \in (0, 1)$ with $\starnorm{\bm{X}} > 0$ then 
 \begin{equation} \label{e:multintbyparts}
 \p_\theta f +\p_{\theta'} f \in \mathcal{S}_{k+1}.
 \end{equation}
\end{lemma}

\begin{poof}
It is clear that we have only to prove the assertion for one $k$-monomial $f\in S_k^1$. 
Let us adopt the notation of Lemma \ref{GPL1}. We set:
\begin{equation*}
f=\prod_{l=0}^k\phi_l^{\alpha_l}\psi_l^{\beta_l}, \; \phi_l=\frac{\chi_{l, 1}}{\abs{\bm{\chi}_0}}, \;  \psi_l=\frac{\chi_{l, 2}}{\abs{\bm{\chi}_0}}.
\end{equation*}
Let us compute $\p_\theta f +\p_{\theta'} f$. We have:
\begin{equation*}
\begin{split}
\p_\theta f+\p_{\theta'} f &=\sum_{l=0}^k (\alpha_lA_l+\beta_lB_l),\\
A_l&=(\p_\theta \phi_l+\p_{\theta'} \phi_l) \phi_l^{\alpha_l-1}\psi_l^{\beta_l}\prod_{i\neq l} \phi_i^{\alpha_i}\psi_i^{\beta_i},\quad 
B_l=(\p_\theta \psi_l+\p_{\theta'}\psi_l) \phi_l^{\alpha_l}\psi_l^{\beta_l-1}\prod_{i\neq l} \phi_i^{\alpha_i}\psi_i^{\beta_i}.
\end{split}
\end{equation*}
We may compute:
\begin{equation*}
\p_\theta \phi_l+\p_{\theta'} \phi_l=\frac{\chi_{l+1,1}}{\abs{\bm{\chi}_0}}-\frac{\chi_{l,1}\bm{\chi}_0\cdot \bm{\chi}_1}{\abs{\bm{\chi}_0}^3}.
\end{equation*}
It is clear that each $A_l$ gives rise to a linear combination of monomials
and that each monomial has $k+1$ total derivatives.
Therefore, each monomial in $A_l$ is a $(k+1)$-monomial 
belonging to $\mc{S}_{k+1}^1$. The same considerations apply for $B_l$. This concludes the proof.
\end{poof}

The next lemma provides an explicit representation of high order derivatives on the nonlinear remainder $\mathcal{R}(\bm{X})$. 
\begin{lemma}\label{l:movederiv}
Assume $\bm{X}\in C^{k+1, \alpha}(\mbs)$  with $k\geq 1$  and $\alpha \in (0, 1)$, and assume $\starnorm{\bm{X}} > 0$. Then 
\begin{equation}  \label{e:highdervSIO}
\p_{\theta}^{k}\mc{R}\left( \bm{X}  \right)
= \sum_{j=1}^{k} \int_{\mbs} \p_{\theta'}P^k_{k+1-j} \p_{\theta'}^{j}\bm{X}' d\theta' 
+   F_C [\p_{\theta}^{k+1}\bm{X}] +  F_T [\p_{\theta}^{k+1}\bm{X}],
\end{equation}
where $P^k_{k+1-j} \in \mc{S}_{k+1-j}$ for  $j \in \{1,\ldots, k\}$. 
\end{lemma}
\begin{poof}
We will prove  \eqref{e:highdervSIO} inductively using the structure of the nonlinear remainder and the form of the kernels in the class, $\mc{S}_j$.  
First recall   
\begin{align*}
\mc{R}(\bm{X}) &= F_C [ \p_\theta \bm{X} ]+ F_T [ \p_\theta \bm{X} ]
\end{align*}
where
\begin{align*}
F_C [ \bm{u} ]
&= {1\over 4\pi}  \int_{\mbs} (\p_{\theta'}M_0)\bm{u}'d\theta', \quad M_0(\theta,\theta')=\log\paren{\frac{\abs{\Delta \bm{X}}}{2\abs{\sin((\theta-\theta')/2)}}}
=\log\paren{\frac{\abs{\bm{\chi}_0}}{2\abs{\sin((\theta-\theta')/2)}}},
\end{align*}
and
\begin{align*}
F_T [ \bm{u} ] &= {1\over 4\pi}  \int_{\mbs}(\p_\theta'N_0)\bm{u}'d\theta', \quad
N_0(\theta,\theta')=-\paren{\frac{\Delta \bm{X}\otimes \Delta \bm{X}}{\abs{\Delta\bm{X}}^2}}
=-\paren{\frac{\bm{\chi}_0\otimes \bm{\chi}_0}{\abs{\bm{\chi}_0}^2}}\in \mc{S}_0.
\end{align*}
We first note that we can convert derivatives in $\theta$ into $\theta'$ derivatives so long as $\bm{u} \in C^1(\mbs)$ and $\bm{X} \in C^{2}(\mbs)$ 
as follows. Using \eqref{pthetaFC=Au} and \eqref{Aopdef}, we have:
\begin{equation}\label{generateCot0}
\begin{split}
\p_\theta F_C[\bm{u}]&=\frac{1}{4\pi}\int_{\mbs}\p_{\theta'}\p_{\theta}M_0(\bm{u}'-\bm{u})d\theta'\\
&=\frac{1}{4\pi}\int_{\mbs}\paren{-\p_{\theta'}^2M_0+\p_\theta'\paren{\p_{\theta}M_0+\p_{\theta'}M_0}(\bm{u}'-\bm{u})}d\theta'\\
&=\frac{1}{4\pi}\int_{\mbs}\paren{\p_{\theta'}M_0(\p_{\theta'}\bm{u}')+\p_{\theta'}\paren{(\p_{\theta}+\p_{\theta'})M_0}(\bm{u}'-\bm{u})}d\theta'
\end{split}
\end{equation}
Let
\begin{equation*}
M_1=\p_\theta M_0+\p_{\theta'}M_0=\frac{\bm{\chi}_0\cdot \bm{\chi}_1}{\abs{\bm{\chi}_0}^2}\in \mc{S}_1.
\end{equation*}
From \eqref{generateCot0}, we have:
\begin{equation}\label{generateCot}
\p_\theta F_C[\bm{u}]= F_C[\p_\theta\bm{u}]+\frac{1}{4\pi}\int_{\mbs} (\p_\theta'M_1)\bm{u}'d\theta'
\end{equation}
where we used the fact that $M_1$ is continuous so that the integral $\p_\theta'M_1 \bm{u}$ vanishes.

In much the same way, we have:
\begin{equation}\label{generateTensor}
\p_\theta F_T[\bm{u}]= F_T[\p_\theta\bm{u}]+\frac{1}{4\pi}\int_{\mbs} (\p_\theta'N_1)\bm{u}'d\theta'
\end{equation}
where, by Lemma \ref{Sk1bounds}, 
\begin{equation*}
N_1=\p_\theta N_0+\p_{\theta'}N_0\in \mc{S}_1.
\end{equation*}
Combining \eqref{generateCot} and \eqref{generateTensor}, we have:
\begin{equation} \label{e:F_comm}
 \p_{\theta}  F_C [ \bm{u} ]    + \p_{\theta}  F_T [ \bm{u} ]    =  \int_{\mbs}  \p_{\theta'} P_1^1  
 \bm{u}'d \theta' +  F_C [\p_\theta  \bm{u} ]    + F_T [ \p_\theta \bm{u} ], \; P^1_{1} :=  {1\over 4\pi}(M_1 + N_1)\in \mc{S}_1
\end{equation}
In particular, we can write,
\begin{equation}  \label{e:Rdiffcase1}
\p_{\theta}\mc{R}(\bm{X}) = \int_{\mbs} \p_{\theta'} P^1_{1} \p_{\theta'} \bm{X}' d \theta' 
+  F_C \left[ \p_{\theta}^2\bm{X} \right] + F_T \left[ \p_{\theta}^2\bm{X} \right] 
\end{equation}
which satisfies \eqref{e:highdervSIO} with $k=1$.  

We proceed by mathematical induction. Suppose the statement of the Proposition is true for some $k\geq 1$.
We show that it is true for $k+1$.
Assume $\bm{X}\in C^{k+2,\alpha }(\mbs)$. Since $\bm{X}\in C^{k+2,\alpha}(\mbs)\subset C^{k+1,\alpha}(\mbs)$,
\eqref{e:highdervSIO} is true by the induction hypothesis. 
If we differentiate \eqref{e:highdervSIO} in $\theta$ and write out the resulting nonlinear kernel, we find 
\begin{align}\label{e:indctRform2}
\p_{\theta}^{k+1}\mc{R}(\bm{X}) 
& =
\sum_{j=1}^{k} \int_{\mbs} \p_{\theta'} \p_\theta \left( P^k_{k+1-j} \right) \left(  \p_{\theta'}^{j}\bm{X}'   - \p_{\theta}^{j}\bm{X}    \right) d \theta'
+  \p_\theta F_C [ \p_\theta^{k+1} \bm{X} ] +  \p_\theta F_T[ \p_\theta^{k+1} \bm{X} ] .
\end{align}
A rigorous justification of this calculation will proceed in the same way as we obtained \eqref{pthetaFC=Au}
in the proof of Proposition \ref{C0cotgain}. We omit this proof.
Since the contour $\bm{X}$ is smooth enough, we can use \eqref{e:multintbyparts} to move derivatives off of the kernels.  In particular we write:
\begin{align*}
\p_\theta P^k_{k+1-j} = - \p_{\theta'} P^k_{k+1 - j} + \widetilde{P}^k_{k+2-j}
\end{align*}
%
for some new multipliers $\widetilde{P}^k_{k+2-j} \in \mc{S}_{k+2-j}$.  Furthermore, from 
\eqref{e:F_comm},  
\[
\p_\theta F_C [ \p_\theta^{k+1} \bm{X} ] + \p_\theta F_T [ \p_\theta^{k+1} \bm{X} ]
= \int_{\mbs} \p_{\theta'} P_1^1 \p_{\theta'}^{k+1} \bm {X}'d \theta' 
+ F_C [ \p_\theta^{k+2} \bm{X}]  +F_T [ \p_\theta^{k+2} \bm{X} ] .
\]
We can now rewrite the terms of \eqref{e:indctRform2} as follows:
\begin{equation*}
\begin{split}
\p_{\theta}^{k+1}\mc{R}(\bm{X}) 
& = 
\sum_{j=1}^{k} \int_{\mbs} \p_{\theta'} \paren{- \p_{\theta'} P^k_{k+1 - j} + \widetilde{P}^k_{k+2-j}}\left(  \p_{\theta'}^{j}\bm{X}' - \p_\theta^j \bm{X} \right)   d\theta'   \\
& \quad +\int_{\mbs} \p_{\theta'} P^1_1 \p_{\theta'}^{k+1} \bm {X}' d \theta' 
+  F_C [ \p_\theta^{k+2} \bm{X} ] + F_T[ \p_\theta^{k+2} \bm{X} ] \\
& = 
\sum_{j=1}^{k} \paren{\int_{\mbs} \p_{\theta'} P^k_{k+1-j} \p_{\theta'}^{j+1}\bm{X}'  d\theta'  
+ \int_{\mbs} \p_{\theta'} \widetilde P^k_{k+2-j} \p_{\theta'}^{j}\bm{X}'  d\theta' }\\
& \quad +\int_{\mbs} \p_{\theta'} P^1_1 \p_{\theta'}^{k+1} \bm {X}' d \theta' 
+  F_C [ \p_\theta^{k+2} \bm{X} ] + F_T[ \p_\theta^{k+2} \bm{X} ]  \\
& = \sum_{j=1}^{k+1} \int_{\mbs}  \p_{\theta'}  P^{k+1}_{k+2-j} \p_{\theta'}^{j}\bm{X}' d\theta' 
+ F_C [ \p^{k+2}_\theta \bm{X} ]  +  F_T [ \p_\theta^{k+2} \bm{X} ],
\end{split}
\end{equation*}
where we defined
\begin{align*}
 P^{k+1}_{k+2-j} = \left\{ 
 \begin{array}{ll}
\widetilde P^{k}_{k+1}  & \hbox{ if } j = 1 \\
P^k_{k+2 - j} + \widetilde P^k_{k+2-j} & \hbox{ if } 2 \leq j \leq k  \\
P_1^k + P_1^1 & \hbox{ if }  j = k + 1.
\end{array}
\right.
\end{align*}  
It may be easily checked that $P^{k+1}_{k+2 - j} \in S_{k+2-j}$ for $1\leq j \leq k+1$.  
\end{poof}

For the next step, we use our calculus lemmas to show that the remainder terms satisfy the following smoothing estimate. 


\begin{lemma}\label{l:CarryRegularity}
Assume $\bm{X}\in C^{n, \alpha}(\mbs)$  with $n\geq 2$  where $\alpha \in (1/2, 1)$, and assume $\starnorm{\bm{X}} > 0$.  Then
\begin{equation} \label{e:Rregimprove}
\chnorm{\mc{R}(\bm{X}) }{n , 2 \alpha-1 } \leq C\paren{\frac{\chnorm{\bm{X}}{n,\alpha}}{\starnorm{\bm{X}}}}^{n+2}\chnorm{\bm{X}}{n,\alpha},
\end{equation}
where the constant $C$ depends only on $n$ and $\alpha$.
\end{lemma}
As the reader will see from the proof, the above bound is suboptimal. It will, however, be sufficient for our purposes.

\begin{poof}[Proof of Lemma \ref{l:CarryRegularity}]

In order to prove Lemma~\ref{l:CarryRegularity}, we express the terms of \eqref{e:highdervSIO} in three parts,
\[
\p_\theta^{n-1} \mc{R}(\bm{X}) = G_{n-1} + F_C[\p_\theta^n \bm{X}] + F_T [\p_\theta^n \bm{X}], 
\] 
where
\begin{align*}
G_{n-1} = \sum_{j=1}^{n-1}  G^{n-1}_j &:= \sum_{j=1}^{n-1} \int_{\mbs} \p_{\theta'}  P^{n-1}_{n-j}\p_{\theta'}^j \bm{X}' d\theta'  .
\end{align*}

From Propositions~\ref{C0cotgain} and \ref{C0tengain} we have the bounds 
\begin{align}\label{highFCFTest}
\chnorm{\p_\theta F_C [ \p^{n-1}_\theta{ \bm{X} } ] +  \p_\theta F_T [ \p^{n-1}_\theta{ \bm{X} } ] }{0, 2\alpha - 1} 
&\leq C\frac{\chnorm{\bm{X}}{1,\alpha}^3}{\starnorm{\bm{X}}^3} \chnorm{\bm{X}}{n,\alpha}.
\end{align}

We must now show that $\p_\theta G^{n-1}_j$ are in $C^{0,2\alpha-1}$. The important point here is that $P^{n-1}_{n-j}\in \mc{S}_{n-j}$, 
and thus Lemmas \ref{GPL1} and \ref{GPL2} are applicable. 
Suppose $P^{n-1}_{n-j}$ can be written as:
\begin{equation*}
P^{n-1}_{n-j}=\sum_{\ell=1}^M a_\ell g_\ell, \; g_\ell \in \mc{S}_{n-j}^1.
\end{equation*}
Any of the terms $g_\ell$ has the form:
\begin{equation*}
g_\ell=\prod_{l=0}^{n-j}\paren{\frac{\chi_{l, 1}}{\abs{\bm{\chi}_0}}}^{\alpha_l} \paren{\frac{\chi_{l, 2}}{\abs{\bm{\chi}_0}}}^{\beta_l},\quad
\sum_{l=1}^{n-j}l(\alpha_l+\beta_l)=n-j.
\end{equation*}
Then, using a procedure that is exactly the same as the proof of Proposition \ref{C0tengain} (or Proposition \ref{C0cotgain}), we have:
\begin{align*}
\chnorm{\p_\theta\int_{\mbs}\p_{\theta'}g_\ell \p_{\theta'}^j\bm{X}'d\theta'}{0,2\alpha-1}=C\frac{\chnorm{\bm{X}}{1,\alpha}^3\prod_{l=1}^{n-j}\chnorm{\p_{\theta}^l\bm{X}}{1,\alpha}^{\alpha_l+\beta_l}}{\starnorm{\bm{X}}^{N+3}}
\chnorm{\p_\theta^j\bm{X}}{0,\alpha}
\end{align*}
where $N$ is given as in \eqref{Ninkmonomial}. A rather crude over-estimation yields:
\begin{equation}\label{Gjn-1est}
\chnorm{\p_\theta G_j^{n-1}}{0,2\alpha-1}\leq C\paren{\frac{\chnorm{\bm{X}}{n-j+1,\alpha}}{\starnorm{\bm{X}}}}^{n-j+3}\chnorm{\bm{X}}{j,\alpha}.
\end{equation}
Combining \eqref{highFCFTest} and \eqref{Gjn-1est}, we obtain the following estimate:
\begin{equation*}
\chnorm{\p_\theta^n\mc{R}(\bm{X})}{0,2\alpha-1}\leq C\paren{\frac{\chnorm{\bm{X}}{n,\alpha}}{\starnorm{\bm{X}}}}^{n+2}\chnorm{\bm{X}}{n,\alpha}.
\end{equation*}
Recall from Proposition \ref{C0cotgain} and \ref{C0tengain} that 
\begin{equation*}
\chnorm{\mc{R}(\bm{X})}{1,2\alpha-1}\leq C\paren{\frac{\chnorm{\bm{X}}{1,\alpha}}{\starnorm{\bm{X}}}}^{3}\chnorm{\bm{X}}{1,\alpha}.
\end{equation*}
We thus have the desired estimate.
\end{poof}



Finally, using the above proposition, we have the following proposition.
\begin{proposition}\label{p:spatialsmooth}
Suppose $\bm{X}(t)$ is a mild solution to the Peskin problem in the interval $[0,T]$ with initial $\bm{X}_0\in h^{1,\gamma}(\mbs)$.
For any $n \in \mathbb{N}$, $\gamma \in (0, 1)$, $\bm{X}(t)$
is in $C([\epsilon, T]; C^{n,\alpha}(\mbs))$  for any $\epsilon \in (0, T)$ and any $\alpha \in (0, 1)$.
\end{proposition}

\begin{poof}

1. To prove higher regularity, we will be estimating a differentiated form of the mild solution, expressed in a more general form,
\begin{align} \label{e:generalYmild}
\bm{Y}(t) = e^{ (t-t_0) \Lambda} \bm{Y}(t_0) + \int_{t_0}^t e^{(t -s) \Lambda} F(\bm{Y}(s))ds. 
\end{align}
If we assume that $\bm{Y}(t_0) \in C^\beta(\mbs)$ and $F(\bm{Y}(s)) \in C([t_0,T]; C^\alpha(\mbs))$ for all $t_0 \leq s \leq t$
then the semigroup estimate \eqref{mainSGest}, implies that for any $\alpha \geq 0$ and $\delta \in (0,1]$ with
$1 + \alpha - \delta$ non-integer valued,   
\begin{equation} \label{e:mildgenform}
\chnorm{\bm{Y}(t) }{1 + \alpha -\delta} \leq  {C \over (t-t_0)^{1+\alpha-\beta-\delta}}\chnorm{\bm{Y}(t_0)}{\beta} + C \sup_{t_0\leq s \leq t} \chnorm{F(\bm{Y}(s))}{\alpha }.
\end{equation}
This follows from \eqref{mainSGest} and 
\begin{align*}
\chnorm{\bm{Y}(t) }{1 + \alpha - \delta}
& \leq \chnorm{e^{t \Lambda} \bm{Y}(t_0) }{1 +  \alpha - \delta}
+ \int_{t_0}^t \chnorm{ e^{(t-s) \Lambda} F(\bm{Y}(s))}{1 +  \alpha - \delta } \\
& \leq {C \over (t - t_0)^{1 + \alpha-\delta - \beta}}\chnorm{ \bm{Y}(t_0) }{\beta} + \int_{t_0}^t  (t -s)^{\delta - 1} \chnorm{ F(\bm{Y}(s))}{\alpha } ds .
\end{align*}
We will be moving $t_0$ away from zero during the iteration procedure.

2.  Let $\bm{X} \in C([0,T]; C^{1, \gamma}(\mbs))$ be the mild solution generated in Theorem~\ref{LWPTheorem}.  We first claim that for any $\epsilon_2 \in (0,\epsilon)$ our solution 
$\bm{X} \in L^\infty ([\epsilon_2, T]; C^{2, \alpha}(\mbs))$, for all $0 < \alpha <1$. Here, $L^\infty([\epsilon,T]; C^{k,\beta}(\mbs))$ refers to the set of functions of $t$ with values in $C^{k,\beta}(\mbs)$ with bounded $C^{k,\beta}(\mbs)$ norm for $\epsilon\leq 0\leq T$. To show this we define 
an iteration $\gamma_{(1,k)}$ starting with $\gamma_{(1,0)} = \gamma$ and updates $\gamma_{(1, j+1)} = {3\over2} \gamma_{(1,j)}$.  We stop the iteration once $1 < \gamma_{(1,k+1)} < 2$, i.e. $k = \lfloor{ \ln {2 / \gamma} \over \ln {3 / 2} } \rfloor$. Without loss of generality, we may assume that none of these $\gamma_{1, j}$'s fall on an integer value.
Subdividing $0 < \epsilon_{(1,1)} < \epsilon_{(1,2)} < \cdots < \epsilon_{(1,k)} < \epsilon_2$, we use \eqref{e:mildgenform} 
and if $ \bm{X} \in L^\infty ([\epsilon_{(1,j)}, T]; C^{1+ \gamma_{(1,j)}}(\mbs))$  then 
\[
\chnorm{\bm{X}}{1+{3\over2} {\gamma_{(1,j)} }} 
\leq C {1\over (t - \epsilon_{(1,j)})^{\gamma_{(1,j)}\over2}} \chnorm{\bm{X}(\epsilon_{(1,j)}) }{ 1 + \gamma_{(1,j)} } 
+ C \sup_{\epsilon_{1,j} \leq s \leq t} \chnorm{\mathcal{R}(\bm{X}(s))} {{7 \over 4} \gamma_{(1,j)}},
\]
where we have chosen $\alpha = {7  \over 4 } \gamma_{(1,j)}$ and $\delta = {1 \over 4}\gamma_{(1,j)}$ in \eqref{e:mildgenform}.
Since 
\[
\chnorm{\mathcal{R}(\bm{X}(s))} {{7 \over 4} \gamma_{(1,j)}} \leq C ( \chnorm{\bm{X}}{1 + {7\over 8} \gamma_{(1,j)}}, \starnorm{\bm{X}})
\leq C ( \chnorm{\bm{X}}{1 + \gamma_{(1,j)}}, \starnorm{\bm{X}}), 
\]
by Propositions~\ref{C0cotgain}, \ref{C0tengain} and the induction hypothesis, 
then $\bm{X} \in L^\infty ([\epsilon_{(1,j+1)}, T]; C^{1+\gamma_{(1,j+1)}}(\mbs))$.

Therefore, we have $\bm{X} \in L^\infty ([\epsilon_{(1,k)}, T]; C^{1+ \gamma_{(1,k+1)}}(\mbs))$ for some $\gamma_{(1,k+1)} \in (1,2)$, which implies
that $\bm{X} \in L^\infty ([\epsilon_{(1,k)}, T]; C^{2,  \tilde\gamma}(\mbs))$ for some $\tilde \gamma \in (0,1)$.  Finally we choose $\alpha = 2-\delta$ for 
 any small $\delta > 0$,  
then  \eqref{mainSGest} implies 
\[
\chnorm{\bm{X}}{2, 1 - 2\delta } 
\leq C {1\over (t - \epsilon_{(1,k)})^{1 - 2 \delta - \tilde\gamma}} \chnorm{\bm{X}(\epsilon_{(1,j)}) }{2 + \tilde \gamma } 
+ C \sup_{\epsilon_{(1,k)} \leq s \leq t} \chnorm{\mathcal{R}(\bm{X}(s))} {1, 1-\delta},
\]
and since $\chnorm{\mathcal{R}(\bm{X}(s))} {1, 1-\delta} \leq C ( \chnorm{\bm{X}(s)} {1, 1- {\delta\over 2}}, \starnorm{\bm{X}}) \leq C ( \chnorm{\bm{X}(s)} {2 , \tilde \gamma}, \starnorm{\bm{X}})$, by  Propositions~\ref{C0cotgain}, \ref{C0tengain}, 
and since $\delta > 0$ is arbitrarily small, we achieve regularity for  the full range of H\"older norms.   

3. We now show that $\bm{X} \in L^\infty ([\epsilon, T], C^{n , \alpha}(\mbs))$ for any $\alpha \in (0,1)$.   To show this we subdivide $0 < \epsilon_2 < \epsilon_3 < \ldots < \epsilon_n = \epsilon$, and we suppose  that $\bm{X} \in L^\infty([\epsilon_k, T], C^{k , \alpha}(\mbs))$ for some $k \geq 2$ and all $\alpha \in (0,1)$.   
By Lemma~\ref{l:CarryRegularity} we find that $\p_\theta^k \mathcal{R}(\bm{X}) \in C^{0,\tilde \alpha}(\mbs)$ for any $\tilde \alpha \in (0,1)$.  We can then use estimate  \eqref{e:mildgenform} and find that 
\[
\chnorm{ \p^k_{\theta} \bm{X} (t) }{1 , \tilde \alpha - \delta}  \leq {C \over (t - t_0)^{1+\tilde \alpha -2 \delta}} \chnorm{ \p_\theta^k \bm{X}(\epsilon_k)   }{0, \delta} + 
\sup_{\epsilon_{k} \leq s \leq t} \chnorm{ \p_\theta^k \mc{R}(\bm{X}(s))}{ 0, \tilde \alpha} ,
\]
and since $\tilde \alpha$ can be made arbitrary in $(0,1)$, then 
 $\bm{X} \in L^\infty ([\epsilon_{k+1}, T]; C^{k +1, 1  - 2 \delta}(\mbs))$ by suitable choices of $\alpha, \tilde \alpha, \delta$, 
 and so $\bm{X} \in L^\infty([\epsilon_{k+1}, T]; C^{k +1, \alpha}(\mbs))$ for any $\alpha \in (0,1)$.  
This process can be repeated $n - 1$ times to conclude that  $\bm{X}\in L^\infty([\epsilon_n, T]; C^{n , \alpha}(\mbs))$ for any $\alpha \in (0, 1)$. 
In particular we see that $\bm{X}\in L^\infty([\epsilon, T]; C^{n , \alpha}(\mbs))$ for any  $\alpha \in (0,1)$. 

4. We must show that $\bm{X} \in C([\epsilon,T];C^{n,\alpha}(\mbs))$.
Using interpolation on H\"older norms (see, for example, Chapter 1 of \cite{lunardi}), we have, for any $\epsilon\leq t_1\leq t_2\leq T$, 
\begin{align*}
\norm{\bm{X}(t_1)-\bm{X}(t_2)}_{C^{n,\alpha}}
&\leq C\norm{\bm{X}(t_1)-\bm{X}(t_2)}_{C^{k,\beta}}^{1-\delta}\norm{\bm{X}(t_1)-\bm{X}(t_2)}_{C^{1,\gamma}}^\delta\\
&\leq C\paren{\norm{\bm{X}(t_1)}_{C^{k,\beta}}+\norm{\bm{X}(t_2)}_{C^{k,\beta}}}^{1-\delta}\norm{\bm{X}(t_1)-\bm{X}(t_2)}_{C^{1,\gamma}}^\delta\\
&\text{ where } n+\alpha=\delta(1+\gamma)+(1-\delta)(k+\beta), \delta>0.
\end{align*}
Note that $\bm{X}(t)\in L^\infty([\epsilon,T];C^{k,\beta}(\mbs)$ for any $k\in \mathbb{N}$ and $\beta\in (0,1)$.
The above inequality can therefore be made to hold for any $n\in \mathbb{N}$ and $\alpha\in (0,1)$ by 
choosing $k+\beta$ large enough.
We thus see that $\bm{X}\in C([\epsilon,T]; C^{n,\alpha}(\mbs))$ for any $n\in \mathbb{N}$, $\alpha\in(0,1)$.
\end{poof}

\begin{poof}[Proof of Theorem \ref{SmoothTheorem}]
We show that the solution is in $C^1([\epsilon,T]; C^{k,\beta}(\mbs))$ for $\epsilon>0$ for any $n$ and $\alpha\in (0,1)$.
Recall from Theorem \ref{LWPTheorem} that the mild solution $\bm{X}(t)$ is a strong solution and thus, 
\begin{equation}\label{XinC1C0gamma}
\p_t\bm{X}(t)\in C^0([\epsilon,T];C^{0,\gamma}(\mbs)).
\end{equation}
Since $\bm{X}(t)$ is a strong solution, it satisfies:
\begin{equation*}
\p_t\bm{X}=\Lambda \bm{X}+\mc{R}(\bm{X}).
\end{equation*}
By the previous Proposition, $\bm{X}\in C([\epsilon,T];C^{n+1,\alpha}(\mbs))$ for any $n\in \mathbb{N}, \alpha\in (0,1)$. Thus,
$\Lambda\bm{X}\in C([\epsilon,T];C^{n,\alpha}(\mbs))$. Furthermore, by Lemma \ref{l:CarryRegularity}, $\mc{R}(\bm{X})\in L^\infty([\epsilon,T];C^{n+1,2\alpha-1}(\mbs))$
so long as $\alpha\in(1/2,1)$. Thus, we see that $\p_t\bm{X}\in L^\infty([\epsilon,T];C^{n,\alpha}(\mbs))$ for any $n\in \mathbb{N}$ and $\alpha\in(0,1)$.
That $\p_t\bm{X}\in C([\epsilon,T];C^{n,\alpha}(\mbs))$ now follows from the same argument 
as in Step 4 of the proof of Proposition \ref{p:spatialsmooth}, interpolating \eqref{XinC1C0gamma}
with the bound $\p_t\bm{X}\in L^\infty([\epsilon,T];C^{k,\beta}(\mbs))$ for large enough $k+\beta$.
\end{poof}

Note that Theorem \ref{SmoothTheorem} was proved using continuity of $\bm{X}(t)$ in time at low order H\"older norms and 
bounds for $\bm{X}(t)$ in higher order H\"older norms.
The former is the result of Lipschitz continuity of $\mc{R}$ in $C^{1,\gamma}(\mbs)$ given in Proposition \ref{LipRemainder} 
and the latter of bounds on $\mc{R}$ in higher order H\"older norms given in Lemma \ref{l:CarryRegularity}. 
To obtain higher order regularity in time, it would thus suffice to establish continuity of derivatives of $\mc{R}$ in $C^{1,\gamma}(\mbs)$
and bounds on such derivatives in higher order H\"older norms. One should be able to obtain the former by a generalization of the proof
in Proposition \ref{LipRemainder}. Indeed, we will establish Lipschitz continuity of the derivative of $\mc{R}$ in Proposition \ref{doublePertN},
which will be used in the study of stability of steady states. Bounds on the derivative of $\mc{R}$ in higher order H\"older norms should 
follow by arguments similar to those that led to Lemma \ref{l:CarryRegularity}. We will not pursue this here; our paper is already quite long.

\subsection{Classical Solutions and the Equivalence of Formulations}\label{sect:classical}

In this subsection, we discuss classical solutions and the equivalence of formulations.
Given the three formulations we have of the Peskin problem introduced in Section \ref{sect:intro},
we will introduce three notions of classical solutions corresponding to each formulation.

In the statement below, we let $C^k(\Omega_{\rm i,e})$ denote $C^k$ functions in $\Omega_{\rm i,e}$ respectively and 
$C^k(\overline{\Omega}_{\rm i,e})$ denote the $C^k$ functions in $\Omega_{\rm i,e}$ whose
derivatives up to order $k$ are uniformly continuous (so that limiting values of the derivatives 
are well-defined at the interface $\Gamma$).
\begin{definition}[Classical jump and IB solutions]
Let $\bm{X}(t)\in C((0,T];C^2(\mbs))\cap C^1((0,T];C(\mbs))$
and $\bm{X}_0\in C^{1,\gamma}(\mbs), \; 0<\gamma<1$ and $\starnorm{\bm{X}_0}>0$.
The function $\bm{X}(t)$ is a {\em classical jump solution} with initial data $\bm{X}_0$ if:
\begin{enumerate}[label=(\roman*)]
\item \label{condition_for_up} For each $t$ satisfying $0<t\leq T$, there exist functions $\bm{u}(\bm{x},t)$ and $p(\bm{x},t)$ satisfying 
equations \eqref{stokes_eqn}-\eqref{stress_jump}, 
\eqref{Xt=u} and \eqref{upinfty}. The velocity field $\bm{u}$ restricted to $\Omega_{\rm i,e}$ is in 
$C^2(\Omega_{\rm i,e})\cap C^1(\overline{\Omega}_{\rm i,e})$ respectively and pressure $p$
restricted to $\Omega_{\rm i,e}$ is in $C^1(\Omega_{\rm i,e})\cap C(\overline{\Omega}_{\rm i,e})$ respectively.
\item $\bm{X}(t)\to \bm{X}_0$ in the $C^{1,\gamma}$ norm as $t \to 0$.
\end{enumerate}
In item \ref{condition_for_up} above, if we require that $\bm{u}, p$ satisfy \eqref{IB_stokes} in the distributional sense
instead of \eqref{stokes_eqn}-\eqref{stress_jump}, we say that $\bm{X}(t)$ is a {\em classical IB solution}.
\end{definition}

\begin{definition}[Classical BI solution]
Let $\bm{X}(t)\in C((0,T];C^2(\mbs))\cap C^1((0,T];C(\mbs))$
and $\bm{X}_0\in C^{1,\gamma}(\mbs), \; 0<\gamma<1$ and $\starnorm{\bm{X}_0}>0$.
The function $\bm{X}(t)$ is a {\em classical BI solution} with initial data $\bm{X}_0$
if it satisfies \eqref{BIform} for $0<t\leq T$ and 
$\bm{X}(t)\to \bm{X}_0$ in $C^{1,\gamma}$ as $t\to 0$.
\end{definition}

\begin{poof}[Corollary~\ref{equiv_formulations}]
That the classical jump and IB solutions are equivalent follows from standard arguments. We refer the reader, for example, 
to \cite{lai2001remark,li2006immersed} for details. 

Suppose we have a classical IB solution. We now show that it is a classical BI solution.
Let $\bm{X}(t)$ be a classical IB solution and let $\bm{u}(\bm{x},t)$ and $p(\bm{x},t)$ be the corresponding 
velocity and pressure fields. Define, as in \eqref{bndry_integral_u} and \eqref{bndry_integral_p}, 
\begin{align*}
\wt{\bm{u}}(\bm{x},t)&=\int_{\mbs} G(\bm{x}-\bm{X}(\theta',t))\p^2_{\theta'}\bm{X}(\theta',t)d\theta',\\
\wt{p}(\bm{x},t)&=\int_{\mbs} \bm{P}_{\rm st}(\bm{x}-\bm{X}(\theta',t))\cdot \p^2_{\theta'}\bm{X}(\theta',t)d\theta'.
\end{align*}
That $\wt{\bm{u}}$ and $\wt{p}$ satisfy \eqref{IB_stokes} in a distributional sense, again, follows from standard arguments.
We have thus only to show that $\wt{\bm{u}}=\bm{u}$. Once we know that this is true, we may use \eqref{Xt=u}
to see that $\bm{X}$ is indeed a classical BI solution. Consider the functions $\bm{w}=\bm{u}-\wt{\bm{u}}$
and $q=p-\wt{p}$. The functions $\bm{w}$ and $q$ satisfy the following equations in the sense of distributions.
\begin{equation}\label{wqeqn}
-\Delta \bm{w}+\nabla q=0, \; \nabla \cdot \bm{w}=0.
\end{equation}
From this, we see that, in the sense of distributions,
\begin{equation*}
\Delta q=0.
\end{equation*}
By Weyl's theorem 
$q$ is  smooth, and thus $q$ is a harmonic function. By \eqref{upinfty}, $q$ is bounded and thus 
by Liouville's theorem, reduces to a constant. Substituting this back into \eqref{wqeqn}, we see that 
each component of $\bm{w}$ is also harmonic. Now, note 
\begin{equation}\label{udecay}
\wt{\bm{u}}(\bm{x},t)=-\int_{\mbs} \p_{\theta'}G(\bm{x}-\bm{X}(\theta',t))\p_{\theta'}\bm{X}(\theta',t)d\theta'
\end{equation}
The kernel $\p_\theta'G(\bm{x}-\bm{X}')$ behaves like $\abs{\bm{x}-\bm{X}'}^{-1}$ as $\bm{x}\to \infty$, 
and thus, $\wt{\bm{u}}(\bm{x})\to 0$ as $\bm{x}\to \infty$.
By \eqref{upinfty}, $\bm{w}$ is thus bounded and 
tends to $0$ as $\bm{x}\to \infty$. Thus, $\bm{w}=0$.

The steps going from \eqref{BIform} to \eqref{SSD} are easily justified since $\bm{X}(t)\in C^2(\mbs)$ for each fixed $t>0$.
A classical BI solution is thus a strong solution. We know from Lemma \ref{strong_is_mild} that a strong solution is a mild solution.

As established in Theorem \ref{SmoothTheorem}, the mild solution $\bm{X}(t)$ satisfies $\bm{X}, \p_t\bm{X}\in C^\infty(\mbs)$ for every fixed $t>0$.
By Lemma \ref{mild_is_strong} the mild solutions are strong, and we may follow our inference backwards from \eqref{SSD} to \eqref{BIform} 
to see that this solution is in fact a classical BI solution. Standard arguments again show that a classical BI solution 
with smooth $\bm{X}$ and $\p_t\bm{X}$ for fixed $t$ is a classical jump solution (and therefore, also a classical IB solution).
%

Finally, we establish the conservation of area and the energy identity.  Our strong solution satisfies the jump formulation of the Peskin problem, and from this, we see that \eqref{area_conservation} is a direct consequence of  
\eqref{incomp} and \eqref{Xt=u}. To obtain the energy identity \eqref{energy_identity}, take the inner product of \eqref{stokes_eqn} with $\bm{u}$
and integrate by parts. The boundary term at $\bm{x}\to \infty$ vanish thanks to the fact that $\nabla \bm{u}$ behaves like $\abs{\bm{x}}^{-2}$
as $\bm{x}\to \infty$. This decay follows directly from the boundary integral representation \eqref{bndry_integral_u}, which we may 
rewrite as in \eqref{udecay}.
\end{poof}

\section{Equilibria and Stability}\label{sect:equilibria_stability}

\subsection{Equilibria}

The goal of this Section is to calculate the equilibria of the Peskin system 
and determine the stability of the equilibria.
We first calculate the equilibria, or stationary states, of the problem. 
Consider any solution $\bm{X}\in C^{1,\alpha}$ that is stationary.
Note first that stationary solutions are classical jump solutions, as 
established in Corollary \ref{equiv_formulations}. Recall by Proposition \ref{equiv_formulations} that $\bm{X}$
must satisfy the energy identity \eqref{energy_identity}. 
Since the $\bm{X}$ is stationary, the energy $\mc{E}$ is not changing in time. Thus, the dissipation $\mc{D}=0$
and thus $\nabla \bm{u}=0$. Given \eqref{upinfty}, $\bm{u}=0$ everywhere.
From the Stokes equation \eqref{stokes_eqn}, we thus conclude that:
\begin{equation*}
\nabla p=0 \text{ in } \Omega_{\rm i,e}.
\end{equation*}
Thus, the pressure $p$ is constant within each domain. Let these pressure values be $p_{\rm i,e}$ respectively in $\Omega_{\rm i,e}$.
Then, by the stress boundary condition \eqref{stress_jump}, we have:
\begin{equation*}
-\triangle p\bm{n}\abs{\p_\theta \bm{X}}=\p_\theta^2 \bm{X}, \; \triangle p=p_{\rm i}-p_{\rm e}, \; \bm{n}=\abs{\p_\theta \bm{X}}^{-1}R_{\pi/2}^{-1}\p_\theta \bm{X},
\end{equation*}
where $\bm{n}$ is the outward normal (pointing from $\Omega_{\rm i}$ to $\Omega_{\rm e}$) 
and $R_{\pi/2}$ is the $2\times 2$ matrix of rotation by $\pi/2$ in the counter-clockwise direction.
For definiteness, we have assumed that the parametrization of the curve is the counter-clockwise direction.
This immediately yields:
\begin{equation*}
\p_\theta^2\bm{X}=\triangle p R_{\pi/2}\p_\theta \bm{X}.
\end{equation*}
This is an easily solved differential equation. Noting that $\bm{X}(2\pi)=\bm{X}(0)$, we have $\triangle p=1$ and 
\begin{equation*}
\begin{split}
\bm{X}(\theta)= A\bm{e}_{\rm r}+B \bm{e}_{\rm t}+ C_1\bm{e}_x+ C_2 \bm{e}_y,
\end{split}
\end{equation*}
where $A,B,C_1,C_2$ are real numbers and $e_{{\rm t},{\rm t},x,y}$ were defined in \eqref{circular_equilibria}
The only equilibrium configurations, therefore, are circles with equidistant material points. 
To make sure that the curve does not degenerate to a point, we impose the condition $A^2+B^2>0$.
We have thus proved the following:
\begin{proposition}\label{equilibria_are_circles}
The only stationary mild solutions of the Peskin problem are circles with equidistant material points.
We may parametrize this set $\wh{\mc{V}}$ as in \eqref{circular_equilibria}.
\end{proposition}
The above proposition is part of Theorem~\ref{StabilityTheorem}, but we have restated this here for easier reference.

\subsection{Spectral and Linear Stability}\label{spect_linear_stability}

To study the stability of these equilibria, we linearize our equation around the stationary solutions.
We study the spectrum of the resulting linear operator $\mc{L}$ and the decay properties of $e^{t\mc{L}}$, 
the semigroup operator generated by $\mc{L}$.
The linearized operator of equation \eqref{SSD} at $\bm{X}$ is given by:
\begin{equation*}
\mc{L}_{\bm{X}}\bm{Y}=\left.\D{}{\epsilon}\paren{\Lambda (\bm{X}+\epsilon \bm{Y})+\mc{R}(\bm{X}+\epsilon \bm{Y})}\right|_{\epsilon=0}=\Lambda \bm{Y}+\p_{\bm{X}}\mc{R}(\bm{X})\bm{Y}.
\end{equation*}
The derivative $\p_{\bm{X}}\mc{R}$ was already introduced in Proposition \ref{LipRemainder}.
Let us linearize around the unit circle:
\begin{equation*}
\bm{X}_{\star}(\theta)=\begin{pmatrix}\cos(\theta)\\ \sin(\theta)\end{pmatrix}.
\end{equation*}
We may compute the linearization around $\bm{X}_\star$ using expressions \eqref{pvBI} and \eqref{pthetaG}. 
This linearization $\mc{L}=\mc{L}_{\bm{X}_{\star}}$ acting on $\bm{Y}$ is given by:
\begin{equation}
\begin{split}\label{Ldefn}
\mc{L}\bm{Y}&=\Lambda \bm{Y}+\p_{\bm{X}}\mc{R}(\bm{X}_\star)\bm{Y}\\
&=\int_{\mbs} \paren{G(\Delta \bm{X}_{\star})\p_{\theta'}^2 \bm{Y}'+(\p_{\bm{X}} G(\Delta \bm{X}_{\star})\Delta \bm{Y})\p_{\theta'}^2\bm{X}_{\star}'}d\theta',\\
\p_{\bm{X}}G(\Delta \bm{X}_{\star})\Delta \bm{Y}&=
-\frac{\Delta \bm{X}_{\star}\cdot \Delta \bm{Y}}{\abs{\Delta \bm{X}_{\star}}^2}I+\frac{\Delta \bm{Y}\otimes\Delta \bm{X}_{\star}}{\abs{\Delta \bm{X}_{\star}}^2}\\
&\quad +\frac{\Delta \bm{X}_{\star}\otimes\Delta \bm{Y}}{\abs{\Delta \bm{X}_{\star}}^2}-2\frac{\Delta \bm{X}_{\star}\cdot \Delta \bm{Y}}{\abs{\Delta \bm{X}_{\star}}^2}\frac{\Delta \bm{X}_{\star}\otimes\Delta \bm{X}_{\star}}{\abs{\Delta \bm{X}_{\star}}^2}.
\end{split}
\end{equation}
In fact, given the translation, rotation and dilation symmetries discussed in Section \ref{sect:intro} 
(see \eqref{dilation_invariance} and surrounding discussion)
linearization around any circular equilibrium will produce the same linearization. Thus:
\begin{equation}\label{L_same_for_circles}
\mc{L}_{\bm Z}=\mc{L}_{\bm{X}_\star}=\mc{L} \text{ for any } \bm{Z} \text{ that is a circular stationary state}.
\end{equation}
We introduce some notation. Define the projection operators:
\begin{equation}\label{proj_trl}
\mc{P}_{\rm trl}\bm{w}=\frac{1}{2\pi}\paren{\dual{\bm{w}}{\bm{e}_x}\bm{e}_x+\dual{\bm{w}}{\bm{e}_y}\bm{e}_y}, \; \Pi_{\rm trl}\bm{w}=\bm{w}-\mc{P}_{\rm trl}\bm{w},
\end{equation}
where the inner product $\dual{\cdot}{\cdot}$ was introduced in \eqref{inner_product}.
The projection $\mc{P}_{\rm trl}$ extracts the translation component of the function (or curve) $\bm{w}$.
After a rather long calculation, we find that \eqref{Ldefn} can be expressed as follows:
\begin{equation}\label{mcLexp}
\begin{split}
\mc{L}\bm{Y}&=R_\theta \Lambda R_\theta^{-1}\bm{Y}+\frac{1}{4}\mc{P}_{\rm trl}\bm{Y}=R_\theta \Lambda R_\theta^{-1}\Pi_{\rm trl}\bm{Y},\\
(R_\theta \Lambda R_\theta^{-1}\bm{w})(\theta)&=\begin{pmatrix}\cos(\theta) &-\sin(\theta)\\ \sin(\theta) &\cos(\theta) \end{pmatrix} 
\begin{pmatrix} \Lambda & 0 \\ 0 & \Lambda \end{pmatrix} \begin{pmatrix}\cos(\theta) &\sin(\theta)\\ -\sin(\theta) &\cos(\theta) \end{pmatrix}
\begin{pmatrix} u(\theta) \\ v(\theta)\end{pmatrix}, \quad \bm{w}=\begin{pmatrix} u\\ v\end{pmatrix}.
\end{split}
\end{equation}
where $\Lambda$ is the familiar operator that first appeared in \eqref{SSD}. 
The $2\times 2$ rotation matrices $R_\theta$ and $R_\theta^{-1}$ act simply by matrix vector multiplication.
From this expression, the eigenvalues of $\mc{L}$ are easily obtained, which are given by:
\begin{equation*}
\lambda_k=-\frac{k}{4}, \; k=\lbrace 0\rbrace\cup \mathbb{N}.
\end{equation*}
For each eigenvalue, the eigenvectors for are given as follows. 
The eigenspace for $\lambda_0=0$ is the span generated by:
\begin{equation}\label{zero_eigenspace}
\lambda_0=0: \quad \bm{e}_{\rm r},\;  \bm{e}_{\rm t}, \; \bm{e}_x, \; \bm{e}_y,
\end{equation}
where the vectors $\bm{e}_{\rm r,t}, \bm{e}_{x,y}$ were defined in \eqref{circular_equilibria}.
The eigenspace $\mc{V}$ for $\lambda_0=0$ coincide with the set of circular equilibria $\wh{\mc{V}}$ (except for the non-degeneracy condition) 
and reflect the group symmetries of our system.
The vector $\bm{e}_{\rm r,t}$ correspond to dilation and rotational symmetry respectively and $\bm{e}_{x,y}$
with translational symmetry. 
For $\lambda_1=-1/4$, the two-dimensional eigenspace is spanned by:
\begin{equation}\label{primary_decay_mode}
\lambda_1=-\frac{1}{4}: \quad \begin{pmatrix} \cos(2\theta)\\ \sin(2\theta) \end{pmatrix}, \; \begin{pmatrix} -\sin(2\theta)\\ \cos(2\theta) \end{pmatrix}.
\end{equation}
For $\lambda_k=-k/4, \; k\geq 2$, the four-dimensional eigenspace is spanned by:
\begin{equation*}
\lambda_k=-\frac{k}{4},\; k\geq 2:\quad \cos(k\theta)\bm{e}_{\rm r},\; \sin(k\theta)\bm{e}_{\rm r},\; \cos(k\theta)\bm{e}_{\rm t},\; \sin(k\theta)\bm{e}_{\rm t}.
\end{equation*}
Of the above four eigenmodes, 
the two eigenmodes proportional to $\bm{e}_{\rm r}$ correspond to radial deformations with a change in shape from the circular configuration, 
whereas the other two eigenmodes proportional to $\bm{e}_{\rm t}$ correspond to tangential deformations without change in circular shape.
The eigenvectors corresponding to the above eigenvalues 
can be easily checked to form a complete orthogonal set in $L^2(\mbs;\mathbb{R}^2)$, and thus the above 
list exhausts the spectrum of $\mc{L}$ as an operator on $L^2(\mbs;\mathbb{R}^2)$. 
	We have only obtained the spectrum of $\mc{L}$ as an operator defined in $L^2(\mbs)$, not as a closed operator on $C^{1,\gamma}(\mbs)$.It is not difficult to show that the $\mc{L}$ has the same spectrum in $C^{1,\gamma}(\mbs)$. We will, however, not need this information here since the above spectral structure will allow us to explicitly compute the semigroup operator $e^{t\mc{L}}$, from which mapping properties of the semigroup operator $e^{t\mc{L}}$ in $C^{1,\gamma}(\mbs)$ will be determined directly.

It is interesting that the spectral structure of $\mc{L}$ and $\Lambda$ (acting component-wise on $\mathbb{R}^2$ valued functions) are very similar.
The difference between the two can be written in terms of the Hilbert transform:
\begin{equation}\label{LLambdaB}
\mc{L}\bm{w}=\Lambda \bm{w}+\frac{1}{4}\begin{pmatrix} 0 &-\mc{H}\\ \mc{H} & 0 \end{pmatrix}\bm{w}.
\end{equation}

The calculation of the above spectrum around a stationary circle is essentially due to \cite{cortez2004parametric}, where the authors 
treat the problem in which the Navier-Stokes equation replaces the Stokes equation in our problem.
For a Navier-Stokes fluid,  a boundary integral reduction is not available and they instead base their calculation 
on the linearization of the jump and IB formulation of the equations (see Section \ref{sect:intro} and Section \ref{sect:classical}).
Our results can be obtained by setting the Reynolds number (or mass density of the fluid) to $0$ in their results, 
except that they seemed to have missed the principal non-zero eigenspace/eigenvector corresponding to $\lambda_1$ above.

Except for the eigenvalues corresponding to the four-dimensional group of symmetries, all eigenvalues are negative.
In this sense, we may say that the uniform circular configurations are spectrally stable.
Our next result shows that the linear evolution semigroup $e^{t\mc{L}}$ has decay properties expected from this spectral structure of $\mc{L}$.
This may be understood as a linear stability result.
Recall the projection operators introduced in \eqref{projPPi}. It is easily checked that:
\begin{equation*}
\mc{L}= R_\theta \Lambda R_\theta^{-1}\Pi=\Pi R_\theta \Lambda R_\theta^{-1}.
\end{equation*}
Let $e^{t\mc{L}}$ be the evolution operator generated by $\mc{L}$. 
From the above, we see that
\begin{equation*}
e^{t\mc{L}}\bm{w}=R_\theta e^{t\Lambda} R_\theta^{-1}\Pi\bm{w}+\mc{P}\bm{w}.
\end{equation*}
We have the following exponential decay estimate for $e^{t\mc{L}}$.
\begin{proposition}\label{SGDecay}
Let $\bm{w}\in C^\alpha(\mathbb{S}^1), \alpha\geq 0$ and let $\beta\geq 0$ satisfy $0\leq \beta-\alpha\leq 1$. Then, 
\begin{equation*}
\norm{e^{t\mc{L}}\Pi \bm{w}}_{C^\beta}\leq Ce^{-t/4}\paren{\frac{1}{t^{\beta-\alpha}}+1}\norm{\Pi \bm{w}}_{C^\alpha},\; t>0,
\end{equation*}
where the constant $C$ above depends only on $\alpha$ and $\beta$.
\end{proposition}

\begin{poof}
Note first that, clearly $R_\theta,R_\theta^{-1}, \Pi$ are bounded operators from $C^\alpha$ to itself for any $\alpha\geq 0$.
Thus, for $t\leq 4\ln 2$, say, 
\begin{equation}\label{expmcLPiw1}
\norm{e^{t\mc{L}}\Pi \bm{w}}_{C^\beta}=\norm{R_\theta e^{t\Lambda}R_\theta^{-1}\Pi \bm{w}}_{C^\beta}
\leq C\norm{e^{t\Lambda}R_\theta^{-1}\Pi \bm{w}}_{C^\beta}
\leq \frac{C}{t^{\beta-\alpha}}\norm{R_\theta^{-1}\Pi \bm{w}}_{C^\alpha}
\leq \frac{C}{t^{\beta-\alpha}}\norm{\Pi \bm{w}}_{C^\alpha}.
\end{equation}
where we used Proposition \ref{holderSG} in the second inequality.

We turn to the estimate for $t\geq 4\ln 2$.
Introduce the following projection acting on scalar valued functions on $\mbs$:
\begin{equation*}
\wh{\mc{P}}w=\frac{1}{2\pi}\dual{w}{1}, \; \wh{\Pi} w=w-\wh{\mc{P}}w.
\end{equation*}
where $\dual{\cdot}{\cdot}$ is the standard $L^2$ inner product. 
It can then be easily checked that:
\begin{equation}\label{expmcLPiw20}
e^{t\mc{L}}\Pi \bm{w}=R_\theta e^{t\Lambda}\wh{\Pi} R_\theta^{-1}\Pi \bm{w}.
\end{equation}
In the above, the operator $e^{t\Lambda}\wh{\Pi}$ acts component-wise. Let us examine the action of this operator. 
The projection $\wh{P}$ is the spectral projection for $e^{t\Lambda}$ for the eigenvalue $0$. 
Performing a calculation similar to \eqref{poisson_kernel} we have:
\begin{equation*}
(e^{t\Lambda}\wh{\Pi}w)(\theta)=\frac{1}{2\pi}\int_{\mbs} e^{-t/4}\wh{P}(e^{-t/4},\theta-\theta')w(\theta')d\theta', \; \wh{P}(r,\theta)=\frac{2(\cos(\theta)-r)}{1-2r\cos\theta+r^2}.
\end{equation*}
Using the fact that $\exp(-t/4)\leq 1/2$ for $t\geq 4\ln 2$, it is easily seen that:
\begin{equation*}
\int_{\mbs} \abs{\wh{P}(e^{-t/4},\theta)}d\theta \leq C, \; \int_{\mbs} \abs{\p_\theta \wh{P}(e^{-t/4},\theta)}d\theta \leq C,
\end{equation*}
where the above $C$ are constants independent of $t$, so long as $t\geq 4\ln 2$. This then immediately shows, 
in a manner similar to Proposition \ref{etlamkatint} that:
\begin{equation*}
\norm{e^{t\Lambda}\wh{\Pi}w}_{C^k}\leq Ce^{-t/4}\norm{w}_{C^k}, \quad \norm{e^{t\Lambda}\wh{\Pi}w}_{C^{k+1}}\leq Ce^{-t/4}\norm{w}_{C^k}, \quad t\geq 4\ln 2.
\end{equation*}
for $k=\lbrace 0\rbrace \cup \mathbb{N}$. We may now use Proposition \ref{holderinterpolation} as in the proof of Proposition \ref{holderSG} to 
obtain:
\begin{equation*}
\norm{e^{t\Lambda}\wh{\Pi}w}_{C^\beta}\leq Ce^{-t/4}\norm{w}_{C^\alpha}, \; 0\leq \beta-\alpha\leq 1, \; t\geq 4\ln 2.
\end{equation*}
Using this, we may estimate \eqref{expmcLPiw20} as follows.
\begin{equation}\label{expmcLPiw2}
\norm{e^{t\mc{L}}\Pi \bm{w}}_{C^\beta}=\norm{R_\theta e^{t\Lambda}\wh{\Pi} R_\theta^{-1}\Pi \bm{w}}_{C^\beta}
\leq C\norm{e^{t\Lambda}\wh{\Pi} R_\theta^{-1}\Pi \bm{w}}_{C^\beta}
\leq Ce^{-t/4}\norm{R_\theta^{-1}\Pi \bm{w}}_{C^\alpha}\leq Ce^{-t/4}\norm{\Pi \bm{w}}_{C^\alpha}.
\end{equation}
Combining estimate \eqref{expmcLPiw1} valid for $t\leq 4\ln 2$ and the above estimate valid for $t\geq 4\ln 2$, we obtain the desired estimate.
\end{poof}

\subsection{Nonlinear Stability}\label{nonlinear_stability}

The goal of this subsection is to prove Theorem \ref{StabilityTheorem}. 
Note that equation \eqref{SSD} can be written as follows:
\begin{equation}\label{LNeqn}
\p_t\bm{X}=\mc{L}\bm{X}+\mc{N}(\bm{X}), \; \mc{N}(\bm{X})=\Lambda\bm{X}+\mc{R}(\bm{X})-\mc{L}\bm{X}
\end{equation}
Our strategy is similar to the one we used to establish the existence/uniqueness of mild solutions of the Peskin problem; 
we turn the above into an integral equation as in \eqref{mild_soln} and use the linear estimate in Proposition~\ref{SGDecay} to 
prove exponential decay to circular equilibria.
This is a standard technique used to study stability of equilibria in ODEs and in parabolic problems, and is sometimes referred 
to as the Lyapunov-Perron method \cite{sell2013dynamics}.
In order to establish nonlinear stability, we must obtain estimates on the remainder term $\mc{N}$.

Recall from Section \ref{sect:intro} that $\mc{V}$ is the kernel of $\mc{L}$, or the four-dimensional eigenspace corresponding to 
the eigenvalue $0$ spanned by \eqref{zero_eigenspace}. According to Proposition \ref{equilibria_are_circles}, 
the set $\wh{\mc{V}}\subset \mc{V}$ corresponds to the set of equilibria. Take any point $\bm{Z}\in \wh{\mc{V}}$. 
Since $\bm{Z}$ is a stationary mild solution (and is, consequently, a strong solution)
\begin{equation*}
0=\mc{L}\bm{Z}+\mc{N}(\bm{Z})=\mc{N}(\bm{Z}),
\end{equation*}
where we used the fact that $\mc{L}\bm{Z}=0$ since $\bm{Z}\in \wh{\mc{V}}\subset\mc{V}$. Furthermore, 
\begin{equation*}
\p_{\bm{X}}\mc{N}[\bm{Z}]\bm{Y}=\Lambda\bm{Y}+\p_{\bm{X}}\mc{R}[\bm{Z}]\bm{Y}-\mc{L}\bm{Y}=0.
\end{equation*}
This follows from the definition of $\mc{L}$ in \eqref{Ldefn} and the fact that the linearization at all stationary circles 
are the same, see \eqref{L_same_for_circles}. For later convenience, let us restate these observations:
\begin{equation}\label{N0pXN0}
\mc{N}(\bm{Z})=0 \text{ and } \p_{\bm X}\mc{N}[\bm{Z}]\bm{Y}=0 \text{ for } \bm{Z}\in \wh{\mc{V}}, \; \bm{Y}\in C^{1,\gamma}(\mbs).
\end{equation}

Our next step is to establish the Lipschitz continuity of $\p_{\bm{X}}\mc{R}$.
\begin{proposition}\label{doublePertN}
Given any $M \geq m > 0$ and convex set $\mc{B}\subset O^{M, m}=\{\bm{Y}\in C^{1, \gamma}: \hnorm{\bm{Y}}{1}{\gamma} \leq M\text{ and } \starnorm{\bm{Y}}\geq m\}$, with $0<\gamma<1$, and $\bm{V}, \bm{W}$ and $\bm{Z}\in \mc{B}$, if $\gamma \neq 1/2$ we have:
\begin{align*}
\hnorm{\p_{\bm{X}}\mc{R}(\bm{V})\bm{Z} - \p_{\bm{X}}\mc{R}(\bm{W})\bm{Z}}{\floor{2\gamma}}{2\gamma-\floor{2\gamma}} &\leq C\frac{M^{4}}{m^{5}}\hnorm{\bm{V} - \bm{W}}{1}{\gamma}\hnorm{\bm{Z}}{1}{\gamma}.
\end{align*}
If $\gamma = 1/2$, we have for any $\alpha \in (0, 1)$
\begin{align*}
\hnorm{\p_{\bm{X}}\mc{R}(\bm{V})\bm{Z} - \p_{\bm{X}}\mc{R}(\bm{W})\bm{Z}}{0}{\alpha} &\leq C\frac{M^{4}}{m^{5}}\hnorm{\bm{V} - \bm{W}}{1}{\gamma}\hnorm{\bm{Z}}{1}{\gamma}.
\end{align*}
\end{proposition}

\begin{poof}
We prove only the $\gamma \neq 1/2$ case as the $\gamma = 1/2$ case is easily derived by adapting the  $\gamma < 1/2$ case. Define:
\begin{align*}
\p^2_{\bm{X}}\mc{R}[\bm{Y}](\bm{V})(\bm{Z}) := \D{}{\epsilon}\left.(\p_{\bm{X}}\mc{R}(\bm{Y}+\epsilon \bm{V})\bm{Z})\right|_{\epsilon=0}.
\end{align*}
We show that, if $\bm{Y}\in C^{1,\gamma}$ and $\starnorm{\bm{Y}}>0$, the following estimate holds:
\begin{equation}\label{p2XRestimate}
\norm{\p^2_{\bm{X}}\mc{R}[\bm{Y}](\bm{V})(\bm{Z})}_{C^{\floor{2\gamma}, 2\gamma - \floor{2\gamma}}}\leq C\frac{\norm{\bm{Y}}_{C^{1,\gamma}}^4}{\starnorm{\bm{Y}}^5}\norm{\bm{V}}_{C^{1,\gamma}}\norm{\bm{Z}}_{C^{1,\gamma}}.
\end{equation}
Once we have this estimate, the desired bound is immediate, for:
\begin{align*}
\hnorm{\p_{\bm{X}}\mc{R}(\bm{V})\bm{Z} - \p_{\bm{X}}\mc{R}(\bm{W})\bm{Z}}{\floor{2\gamma}}{2\gamma - \floor{2\gamma}} &= \hnorm{\int_{0}^{1}\frac{d}{ds}\p_{\bm{X}}\mc{R}((1 - s)\bm{V} + s\bm{W})\bm{Z}ds}{\floor{2\gamma}}{2\gamma - \floor{2\gamma}}\\
&\leq \int_0^1 \norm{\p^2_{\bm{X}}\mc{R}[(1-s)\bm{V}+s\bm{W}](\bm{V}-\bm{W})(\bm{Z})}_{C^{\floor{2\gamma}, 2\gamma - \floor{2\gamma}}}ds\\
&\leq C\frac{M^4}{m^5}\hnorm{\bm{V} - \bm{W}}{1}{\gamma}\hnorm{\bm{Z}}{1}{\gamma},
\end{align*}
where we used the assumption that $\mc{B}$ is convex in the last inequality.

Recall from \eqref{RCRT} that 
\begin{align}\label{FR_FT_defn}
\mc{R}(\bm{X}) = \mc{R}_C(\bm{X}) + \mc{R}_T(\bm{X}).
\end{align}
We thus have:
\begin{align*}
\p_{\bm{X}}\mc{R}(\bm{X})\bm{Z} &= \p_{\bm{X}}\mc{R}_C(\bm{X})\bm{Z} + \p_{\bm{X}}\mc{R}_T (\bm{X})\bm{Z}.
\end{align*}
We shall focus on estimating the $\mc{R}_C$ term as estimates on $\mc{R}_T$ can be obtained in a similar fashion. 
In the proof of Proposition \ref{LipRemainder} we calculated that for any set of constant vectors $\{\bm{V}_i\}_{i = 1}^{2}$,
\begin{align*}
4\pi\p_{\bm{X}}\mc{R}_C(\bm{X})\bm{Z} &= \int_{\mbs}\left(\frac{\Delta\bm{X}\cdot\p_{\theta'}\bm{X}'}{|\Delta\bm{X}|^2}-\frac{1}{2}\cot\paren{\frac{\theta-\theta'}{2}} \right)\left(\bm{V}_1 - \p_{\theta'}\bm{Z}' \right)d\theta'\\
&+ \int_{\mbs}\left(\frac{|\Delta\bm{X}|^2\Delta\bm{X}\cdot\p_{\theta'}\bm{Z}' + |\Delta\bm{X}|^2\Delta\bm{Z}\cdot\p_{\theta'}\bm{X}' - 2\Delta\bm{X}\cdot\p_{\theta'}\bm{X}'\Delta\bm{X}\cdot\Delta\bm{Z}}{|\Delta\bm{X}|^4} \right)\left(\bm{V}_2 - \p_{\theta'}\bm{X}' \right)d\theta'\\
\end{align*}
We may take another derivative to obtain:
\begin{align*}
4\pi\p_{\bm{X}}^2\mc{R}_C[\bm{Y}](\bm{W})\bm{Z} &= \int_{\mbs}\left(\frac{|\Delta\bm{Y}|^2\Delta\bm{W}\cdot\p_{\theta'}\bm{Y}' - \Delta\bm{Y}\cdot\p_{\theta'}\bm{Y}'\Delta\bm{Y}\cdot\Delta\bm{W}}{|\Delta\bm{Y}|^4}\right)\left(\bm{V}_1 - \p_{\theta'}\bm{Z}'\right)d\theta'\\
&+\int_{\mbs}\left(\frac{|\Delta\bm{Y}|^2\Delta\bm{Y}\cdot\p_{\theta'}\bm{W}' - \Delta\bm{Y}\cdot\p_{\theta'}\bm{Y}'\Delta\bm{Y}\cdot\Delta\bm{W}}{|\Delta\bm{Y}|^4}\right)\left(\bm{V}_1 - \p_{\theta'}\bm{Z}'\right)d\theta'\\
&+ \int_{\mbs}2\left(\frac{\Delta\bm{Y}\cdot\Delta\bm{W}\Delta\bm{Y}\cdot\p_{\theta'}\bm{Z}' - \Delta\bm{W}\cdot\p_{\theta'}\bm{Y}'\Delta\bm{Y}\cdot\Delta\bm{Z}}{|\Delta\bm{Y}|^4} \right)(\bm{V}_2 - \p_{\theta'}\bm{Y}')d\theta'\\
&+ \int_{\mbs}2\left(\frac{\Delta\bm{Y}\cdot\Delta\bm{W}\Delta\bm{Z}\cdot\p_{\theta'}\bm{Y}' - \Delta\bm{Y}\cdot\p_{\theta'}\bm{W}'\Delta\bm{Y}\cdot\Delta\bm{Z}}{|\Delta\bm{Y}|^4} \right)(\bm{V}_2 - \p_{\theta'}\bm{Y}')d\theta'\\
&+ \int_{\mbs}\left(\frac{\Delta\bm{Y}\cdot\Delta\bm{Y}\Delta\bm{W}\cdot\p_{\theta'}\bm{Z}' - \Delta\bm{Y}\cdot\p_{\theta'}\bm{Y}'\Delta\bm{W}\cdot\Delta\bm{Z}}{|\Delta\bm{Y}|^4} \right)(\bm{V}_2 - \p_{\theta'}\bm{Y}')d\theta'\\
&+ \int_{\mbs}\left(\frac{\Delta\bm{Y}\cdot\Delta\bm{Y}\Delta\bm{Z}\cdot\p_{\theta'}\bm{W}' - \Delta\bm{Y}\cdot\p_{\theta'}\bm{Y}'\Delta\bm{W}\cdot\Delta\bm{Z}}{|\Delta\bm{Y}|^4} \right)(\bm{V}_2 - \p_{\theta'}\bm{Y}')d\theta'\\
&+ \int_{\mbs}4\Delta\bm{Y}\cdot\Delta\bm{W}\left(\frac{\Delta\bm{Y}\cdot\p_{\theta'}\bm{Y}'\Delta\bm{Y}\cdot\Delta\bm{Z} - |\Delta\bm{Y}|^2\Delta\bm{Y}\cdot\p_{\theta'}\bm{Z}'}{|\Delta\bm{Y}|^6} \right)(\bm{V}_2 - \p_{\theta'}\bm{Y}')d\theta'\\
&+ \int_{\mbs}4\Delta\bm{Y}\cdot\Delta\bm{W}\left(\frac{\Delta\bm{Y}\cdot\p_{\theta'}\bm{Y}'\Delta\bm{Y}\cdot\Delta\bm{Z} - |\Delta\bm{Y}|^2\Delta\bm{Z}\cdot\p_{\theta'}\bm{Y}'}{|\Delta\bm{Y}|^6} \right)(\bm{V}_2 - \p_{\theta'}\bm{Y}')d\theta'\\
&+\int_{\mbs}\left( \frac{|\Delta\bm{Y}|^2\Delta\bm{Y}\cdot\p_{\theta'}\bm{Z}' - \Delta\bm{Y}\cdot\p_{\theta'}\bm{Y}'\Delta\bm{Y}\cdot\Delta\bm{Z}}{|\Delta\bm{Y}|^4}\right)(\bm{V}_3 - \p_{\theta'}\bm{W}')d\theta'\\
&+\int_{\mbs}\left( \frac{|\Delta\bm{Y}|^2\Delta\bm{Z}\cdot\p_{\theta'}\bm{Y}' - \Delta\bm{Y}\cdot\p_{\theta'}\bm{Y}'\Delta\bm{Y}\cdot\Delta\bm{Z}}{|\Delta\bm{Y}|^4}\right)(\bm{V}_3 - \p_{\theta'}\bm{W}')d\theta'\\
&=: F_1^{1} + F_1^{2} + F_2^1 + F_2^2 + F_2^3 + F_2^4 + F_2^5 + F_2^6 + F_3^1 + F_3^2,
\end{align*}
where the vectors $\bm{V}_i$ are once again arbitrary. We will show that $F_1^1 \in C^{\floor{2\gamma}, 2\gamma - \floor{2\gamma}}$ as all of these terms use the exact same proof argument. We have

\begin{align*}
F_1^1 (\theta) &= \int_{\mbs}\left(\frac{|\Delta\bm{Y}|^2\Delta\bm{W}\cdot\p_{\theta'}\bm{Y}' - \Delta\bm{Y}\cdot\p_{\theta'}\bm{Y}'\Delta\bm{Y}\cdot\Delta\bm{W}}{|\Delta\bm{Y}|^4}\right)\left(\bm{V}_1 - \p_{\theta'}\bm{Z}'\right)d\theta'.
\end{align*} 
Note that the kernel can be rewritten as

\begin{align}\label{FormOfLemma}
\frac{|\Delta\bm{Y}|^2\Delta\bm{W}\cdot \left(\p_{\theta'}\bm{Y}' - \frac{\Delta\bm{Y}}{\theta - \theta'} \right)}{|\Delta \bm{Y}|^4} + \frac{\Delta\bm{Y}\cdot\Delta\bm{W}\Delta\bm{Y}\cdot \left(\frac{\Delta\bm{Y}}{\theta - \theta'} - \p_{\theta'}\bm{Y}' \right)}{|\Delta\bm{Y}|^4}
\end{align}
which is a sum of functions of the forms used in Lemmas~\ref{GPlog1} and \ref{GPlog2}. Using Lemma~\ref{GPlog1} gives the bound

\begin{align*}
2\frac{\hnorm{\bm{W}}{1}{\gamma}\hnorm{\bm{Y}}{1}{\gamma}}{\starnorm{\bm{Y}}^2}|\theta - \theta'|^{\gamma - 1}.
\end{align*}
Using this and taking $\bm{V}_1 = \p_{\theta}\bm{Z}$ yields

\begin{align*}
|F_1^1 (\theta)| &\leq \int_{\mbs}2\frac{\hnorm{\bm{W}}{1}{\gamma}\hnorm{\bm{Y}}{1}{\gamma}}{\starnorm{\bm{Y}}^2}|\theta - \theta'|^{\gamma - 1}|\p_{\theta'}\bm{Z}' - \p_{\theta}\bm{Z}|d\theta'\\
&\leq 2\frac{\hnorm{\bm{Z}}{1}{\gamma}\hnorm{\bm{W}}{1}{\gamma}\hnorm{\bm{Y}}{1}{\gamma}}{\starnorm{\bm{Y}}^2}\int_{\mbs}|\theta - \theta'|^{2\gamma - 1}d\theta' \leq  C\frac{\hnorm{\bm{Z}}{1}{\gamma}\hnorm{\bm{W}}{1}{\gamma}\hnorm{\bm{Y}}{1}{\gamma}}{\starnorm{\bm{Y}}^2},
\end{align*}
where $C$ depends only on $\gamma$. If $\gamma < 1/2$, we take $\bm{V}_1 = \p_{\theta}\bm{Z}$. Then,

\begin{align*}
\diff_h F_1^1  &= \int_{\mbs} \trl_{h}\left(\frac{|\Delta\bm{Y}|^2\Delta\bm{W}\cdot\p_{\theta'}\bm{Y}' - \Delta\bm{Y}\cdot\p_{\theta'}\bm{Y}'\Delta\bm{Y}\cdot\Delta\bm{W}}{|\Delta\bm{Y}|^2}\right)(\p_{\theta}\bm{Z} - \p_{\theta'}\bm{Z}')d\theta'\\
&- \int_{\mbs}\left(\frac{|\Delta\bm{Y}|^2\Delta\bm{W}\cdot\p_{\theta'}\bm{Y}' - \Delta\bm{Y}\cdot\p_{\theta'}\bm{Y}'\Delta\bm{Y}\cdot\Delta\bm{W}}{|\Delta\bm{Y}|^2}\right)(\p_{\theta}\bm{Z} - \p_{\theta'}\bm{Z}')d\theta'.
\end{align*}
Using our standard sets $\mc{I}_s$ and $\mc{I}_f$ we have

\begin{align*}
\diff_h F_1^1  &= \int_{\mc{I}_s} \trl_{h}\left(\frac{|\Delta\bm{Y}|^2\Delta\bm{W}\cdot\p_{\theta'}\bm{Y}' - \Delta\bm{Y}\cdot\p_{\theta'}\bm{Y}'\Delta\bm{Y}\cdot\Delta\bm{W}}{|\Delta\bm{Y}|^2}\right)(\p_{\theta}\bm{Z} - \p_{\theta'}\bm{Z}')d\theta'\\
&- \int_{\mc{I}_s}\left(\frac{|\Delta\bm{Y}|^2\Delta\bm{W}\cdot\p_{\theta'}\bm{Y}' - \Delta\bm{Y}\cdot\p_{\theta'}\bm{Y}'\Delta\bm{Y}\cdot\Delta\bm{W}}{|\Delta\bm{Y}|^2}\right)(\p_{\theta}\bm{Z} - \p_{\theta'}\bm{Z}')d\theta'\\
&+ \int_{\mc{I}_f}\diff_h \left(\frac{|\Delta\bm{Y}|^2\Delta\bm{W}\cdot\p_{\theta'}\bm{Y}' - \Delta\bm{Y}\cdot\p_{\theta'}\bm{Y}'\Delta\bm{Y}\cdot\Delta\bm{W}}{|\Delta\bm{Y}|^2}\right)(\p_{\theta}\bm{Z} - \p_{\theta'}\bm{Z}')d\theta'.
\end{align*}
It is clear using the bound found above that the first term can be bounded by

\begin{align*}
\left| \int_{\mc{I}_s} \trl_{h}\left(\frac{|\Delta\bm{Y}|^2\Delta\bm{W}\cdot\p_{\theta'}\bm{Y}' - \Delta\bm{Y}\cdot\p_{\theta'}\bm{Y}'\Delta\bm{Y}\cdot\Delta\bm{W}}{|\Delta\bm{Y}|^2}\right)(\p_{\theta}\bm{Z} - \p_{\theta'}\bm{Z}')d\theta'\right|\\
&\hspace{-4in}\leq 2\frac{\hnorm{\bm{W}}{1}{\gamma}\hnorm{\bm{Y}}{1}{\gamma}\hnorm{\bm{Z}}{1}{\gamma}}{\starnorm{\bm{Y}}^2}\int_{\mc{I}_s}|\theta + h - \theta'|^{\gamma - 1}|\theta - \theta'|^{\gamma}d\theta'\\ 
&\hspace{-4in}\leq 2(3/2)^{\gamma}h^{\gamma}\frac{\hnorm{\bm{W}}{1}{\gamma}\hnorm{\bm{Y}}{1}{\gamma}\hnorm{\bm{Z}}{1}{\gamma}}{\starnorm{\bm{Y}}^2}\int_{\mc{I}_s}|\theta + h - \theta'|^{\gamma - 1}d\theta'\\
&\hspace{-4in}\leq Ch^{2\gamma}\frac{\hnorm{\bm{W}}{1}{\gamma}\hnorm{\bm{Y}}{1}{\gamma}\hnorm{\bm{Z}}{1}{\gamma}}{\starnorm{\bm{Y}}^2},
\end{align*}
where $C$ depends only on $\gamma$. The same bound holds for the second term by the same argument. For the third term, we make use of Lemma~\ref{GPlog2} which gives the bound

\begin{align*}
\left|\diff_h \left(\frac{|\Delta\bm{Y}|^2\Delta\bm{W}\cdot\p_{\theta'}\bm{Y}' - \Delta\bm{Y}\cdot\p_{\theta'}\bm{Y}'\Delta\bm{Y}\cdot\Delta\bm{W}}{|\Delta\bm{Y}|^2}\right)\right| &\leq Ch\frac{\hnorm{\bm{Y}}{1}{\gamma}^2\hnorm{\bm{W}}{1}{\gamma}}{\starnorm{\bm{Y}}^3}|\theta - \theta'|^{\gamma - 2} .
\end{align*}
Thus,

\begin{align*}
\left|\int_{\mc{I}_f}\diff_h \left(\frac{|\Delta\bm{Y}|^2\Delta\bm{W}\cdot\p_{\theta'}\bm{Y}' - \Delta\bm{Y}\cdot\p_{\theta'}\bm{Y}'\Delta\bm{Y}\cdot\Delta\bm{W}}{|\Delta\bm{Y}|^2}\right)(\p_{\theta}\bm{Z} - \p_{\theta'}\bm{Z}')d\theta'\right|\\
&\hspace{-4in}\leq Ch\frac{\hnorm{\bm{Y}}{1}{\gamma}^2\hnorm{\bm{W}}{1}{\gamma}}{\starnorm{\bm{Y}}^3}\int_{\mc{I}_f}|\theta - \theta'|^{\gamma - 2}|\p_{\theta}\bm{Z} - \p_{\theta'}\bm{Z}'|d\theta'\\
&\hspace{-4in}\leq Ch\frac{\hnorm{\bm{Y}}{1}{\gamma}^2\hnorm{\bm{W}}{1}{\gamma}\hnorm{\bm{Z}}{1}{\gamma}}{\starnorm{\bm{Y}}^3}\int_{\mc{I}_f}|\theta - \theta'|^{2\gamma - 2}d\theta'\\
&\hspace{-4in}\leq Ch^{2\gamma}\frac{\hnorm{\bm{Y}}{1}{\gamma}^2\hnorm{\bm{W}}{1}{\gamma}\hnorm{\bm{Z}}{1}{\gamma}}{\starnorm{\bm{Y}}^3},
\end{align*}
with $C$ again depending only on $\gamma$. Combining the above and noting that $\frac{\hnorm{\bm{Y}}{1}{\gamma}}{\starnorm{\bm{Y}}} \geq 1$, we have

\begin{align*}
\hnorm{F_1^1 }{0}{2\gamma} &\leq C\frac{\hnorm{\bm{Y}}{1}{\gamma}^2\hnorm{\bm{W}}{1}{\gamma}\hnorm{\bm{Z}}{1}{\gamma}}{\starnorm{\bm{Y}}^3}.
\end{align*}
We now turn to the case where $\gamma > 1/2$. Using the above, we already have a bound on the $C^0$ norm. As in $(\ref{Aopdef})$ we find that

\begin{align*}
(\p_{\theta}F_1^1)(\theta) &= \int_{\mbs}\p_{\theta}\left(\frac{|\Delta\bm{Y}|^2\Delta\bm{W}\cdot\p_{\theta'}\bm{Y}' - \Delta\bm{Y}\cdot\p_{\theta'}\bm{Y}'\Delta\bm{Y}\cdot\Delta\bm{W}}{|\Delta\bm{Y}|^2}\right)\left(\p_{\theta}\bm{Z}_1 - \p_{\theta'}\bm{Z}'\right)d\theta',
\end{align*}
with $\bm{V}_1 = \p_{\theta}\bm{Z}$ as a necessity for equivalence. We have already shown that the kernel of $F_1^1$ can be written as the sum of functions used in Lemmas~\ref{GPlog1} and \ref{GPlog2}. Using Lemma~\ref{GPlog1} and the triangle inequality yields

\begin{align*}
|\p_{\theta}F_1^1 (\theta)| &\leq \int_{\mbs} \left|\p_{\theta}\left(\frac{|\Delta\bm{Y}|^2\Delta\bm{W}\cdot\p_{\theta'}\bm{Y}' - \Delta\bm{Y}\cdot\p_{\theta'}\bm{Y}'\Delta\bm{Y}\cdot\Delta\bm{W}}{|\Delta\bm{Y}|^4}\right)\right||\p_{\theta}\bm{Z} - \p_{\theta}'\bm{Z}|d\theta'\\
&\leq C\frac{\hnorm{\bm{Y}}{1}{\gamma}^2\hnorm{\bm{W}}{1}{\gamma}\hnorm{\bm{Z}}{1}{\gamma}}{\starnorm{\bm{Y}}^3}\int_{\mbs} |\theta - \theta'|^{2\gamma - 2} d\theta'\\
&\leq C\frac{\hnorm{\bm{Y}}{1}{\gamma}^2\hnorm{\bm{W}}{1}{\gamma}\hnorm{\bm{Z}}{1}{\gamma}}{\starnorm{\bm{Y}}^3}
\end{align*}
where $C$ depends only on $\gamma$ and we used the fact that $\gamma > 1/2$ implies that $2\gamma - 1 > 0$. We now bound the $2\gamma - 1$ semi-norm. Without loss of generality we let $\theta \in (0, 2\pi)$ and  $0 < h < 2\pi$. Using our usual sets of integration $\mc{I}_s$ and $\mc{I}_f$ gives

\begin{align*}
\diff_h \p_{\theta}F_1^1 &= \int_{\mbs}\trl_h\left(\p_{\theta}\left(\frac{|\Delta\bm{Y}|^2\Delta\bm{W}\cdot\p_{\theta'}\bm{Y}' - \Delta\bm{Y}\cdot\p_{\theta'}\bm{Y}'\Delta\bm{Y}\cdot\Delta\bm{W}}{|\Delta\bm{Y}|^4}\right)\right)(\p_{\theta}\bm{Z}(\theta + h) - \p_{\theta'}\bm{Z}')d\theta'\\
&- \int_{\mbs}\p_{\theta}\left(\frac{|\Delta\bm{Y}|^2\Delta\bm{W}\cdot\p_{\theta'}\bm{Y}' - \Delta\bm{Y}\cdot\p_{\theta'}\bm{Y}'\Delta\bm{Y}\cdot\Delta\bm{W}}{|\Delta\bm{Y}|^4}\right)(\p_{\theta}\bm{Z}(\theta) - \p_{\theta'}\bm{Z}')d\theta'\\
&= \int_{\mc{I}_s}\trl_h\left(\p_{\theta}\left(\frac{|\Delta\bm{Y}|^2\Delta\bm{W}\cdot\p_{\theta'}\bm{Y}' - \Delta\bm{Y}\cdot\p_{\theta'}\bm{Y}'\Delta\bm{Y}\cdot\Delta\bm{W}}{|\Delta\bm{Y}|^4}\right)\right)(\p_{\theta}\bm{Z}(\theta + h) - \p_{\theta'}\bm{Z}')d\theta'\\
&- \int_{\mc{I}_s}\p_{\theta}\left(\frac{|\Delta\bm{Y}|^2\Delta\bm{W}\cdot\p_{\theta'}\bm{Y}' - \Delta\bm{Y}\cdot\p_{\theta'}\bm{Y}'\Delta\bm{Y}\cdot\Delta\bm{W}}{|\Delta\bm{Y}|^4}\right)(\p_{\theta}\bm{Z}(\theta) - \p_{\theta'}\bm{Z}')d\theta'\\
&+ \int_{\mc{I}_f}\diff_h\left(\p_{\theta}\left(\frac{|\Delta\bm{Y}|^2\Delta\bm{W}\cdot\p_{\theta'}\bm{Y}' - \Delta\bm{Y}\cdot\p_{\theta'}\bm{Y}'\Delta\bm{Y}\cdot\Delta\bm{W}}{|\Delta\bm{Y}|^4}\right) \right)(\p_{\theta}\bm{Z}(\theta + h) - \p_{\theta'}\bm{Z}')d\theta'\\
&+\int_{\mc{I}_f}\p_{\theta}\left(\frac{|\Delta\bm{Y}|^2\Delta\bm{W}\cdot\p_{\theta'}\bm{Y}' - \Delta\bm{Y}\cdot\p_{\theta'}\bm{Y}'\Delta\bm{Y}\cdot\Delta\bm{W}}{|\Delta\bm{Y}|^4}\right)(\p_{\theta}\bm{Z}(\theta + h) - \p_{\theta}\bm{Z}(\theta))d\theta'\\
&=: I + II + III + IV.
\end{align*}
The first two terms, $I$ and $II$ can be bounded using the same method as we used to bound the $C^1$ norm above. Using Lemma~\ref{GPlog1} we have

\begin{align*}
|I| &\leq C\frac{\hnorm{\bm{Y}}{1}{\gamma}^2\hnorm{\bm{W}}{1}{\gamma}\hnorm{\bm{Z}}{1}{\gamma}}{\starnorm{\bm{Y}}^3}\int_{\mc{I}_s}|\theta + h - \theta'|^{2\gamma - 2} d\theta'\\
&\leq Ch^{2\gamma - 1}\frac{\hnorm{\bm{Y}}{1}{\gamma}^2\hnorm{\bm{W}}{1}{\gamma}\hnorm{\bm{Z}}{1}{\gamma}}{\starnorm{\bm{Y}}^3}.
\end{align*}
The same bound holds for term $II$. For term $III$ we use Lemma~\ref{GPlog2}, inequality (\ref{logest4}). By linearity of the derivative and $\diff_h$ operators,

\begin{align*}
|III| &\leq C\frac{\hnorm{\bm{Y}}{1}{\gamma}^3\hnorm{\bm{W}}{1}{\gamma}\hnorm{\bm{Z}}{1}{\gamma}}{\starnorm{\bm{Y}}^{4}}\int_{\mc{I}_f}\left(h^{\gamma}|\theta - \theta'|^{-2} + h|\theta - \theta'|^{\gamma - 3}\right)|\theta + h - \theta'|^{\gamma}d\theta'.
\end{align*}
On the set $\mc{I}_f$, $|\theta + h - \theta'| < |\theta - \theta'|$. Thus,

\begin{align*}
|III| &\leq C\frac{\hnorm{\bm{Y}}{1}{\gamma}^3\hnorm{\bm{W}}{1}{\gamma}\hnorm{\bm{Z}}{1}{\gamma}}{\starnorm{\bm{Y}}^{4}}\int_{\mc{I}_f}h^{\gamma}|\theta - \theta'|^{\gamma-2} + h|\theta - \theta'|^{2\gamma - 3} d\theta'\\
&\leq Ch^{2\gamma - 1}\frac{\hnorm{\bm{Y}}{1}{\gamma}^3\hnorm{\bm{W}}{1}{\gamma}\hnorm{\bm{Z}}{1}{\gamma}}{\starnorm{\bm{Y}}^{4}}.
\end{align*}
Finally, for term $IV$, we apply the estimate from $(\ref{logest2})$ in Lemma~\ref{GPlog1} to find

\begin{align*}
|IV| &\leq C\frac{\hnorm{\bm{Y}}{2}{\gamma}^2\hnorm{\bm{W}}{1}{\gamma}}{\starnorm{\bm{Y}}^3}\int_{\mc{I}_f}|\theta - \theta'|^{\gamma - 2} |\p_{\theta}\bm{Z}(\theta + h) - \p_{\theta}\bm{Z}(\theta)|d\theta'\\
&\leq Ch^{\gamma}\frac{\hnorm{\bm{Y}}{2}{\gamma}^2\hnorm{\bm{W}}{1}{\gamma}\hnorm{\bm{Z}}{1}{\gamma}}{\starnorm{\bm{Y}}^3}\int_{\mc{I}_f}|\theta - \theta'|^{\gamma - 2} d\theta'\\
&\leq Ch^{2\gamma - 1}\frac{\hnorm{\bm{Y}}{2}{\gamma}^2\hnorm{\bm{W}}{1}{\gamma}\hnorm{\bm{Z}}{1}{\gamma}}{\starnorm{\bm{Y}}^3}.
\end{align*}
Assembling terms $I$ - $IV$ with the $C^0$ and $C^1$ norm bounds gives

\begin{align*}
\hflnorm{F_1^1}{2\gamma} &\leq C\frac{\hnorm{\bm{Y}}{1}{\gamma}^3\hnorm{\bm{W}}{1}{\gamma}\hnorm{\bm{Z}}{1}{\gamma}}{\starnorm{\bm{Y}}^4}.
\end{align*}
The arguments above can be used on the remaining terms $F_j^i$ so that each is in $C^{\floor{2\gamma}, 2\gamma - \floor{2\gamma}}$. The dominant bound out of all of these terms comes from the $F_2^i$ terms and is

\begin{align*}
\hnorm{F_{2}^{i}}{\floor{2\gamma}}{2\gamma - \floor{2\gamma}} &\leq C\frac{\hnorm{\bm{Y}}{1}{\gamma}^4\hnorm{\bm{W}}{1}{\gamma}\hnorm{\bm{Z}}{1}{\gamma}}{\starnorm{\bm{Y}}^5}.
\end{align*}
Assembling the bounds for all of the $F_j^i$ terms yields

\begin{align*}
\hnorm{\p_{\bm{X}}^2\mc{R}_C[\bm{Y}](\bm{W})\bm{Z}}{\floor{2\gamma}}{2\gamma - \floor{2\gamma}}  &\leq C\frac{\hnorm{\bm{Y}}{1}{\gamma}^4\hnorm{\bm{W}}{1}{\gamma}\hnorm{\bm{Z}}{1}{\gamma}}{\starnorm{\bm{Y}}^5}.
\end{align*}

Though tedious, $\p_{\bm{X}}\mc{R}_T(\bm{X})\bm{Z}$ can be linearized as well and will yield the same bound of 

\begin{align*}
\hnorm{\p_{\bm{X}}^2\mc{R}_T[\bm{Y}](\bm{W})\bm{Z}}{\floor{2\gamma}}{2\gamma - \floor{2\gamma}}  &\leq C\frac{\hnorm{\bm{Y}}{1}{\gamma}^4\hnorm{\bm{W}}{1}{\gamma}\hnorm{\bm{Z}}{1}{\gamma}}{\starnorm{\bm{Y}}^5}.
\end{align*}
Estimate \eqref{p2XRestimate} is now immediate.
\end{poof} 

Using the above Proposition together with observation \eqref{N0pXN0}, we obtain the following estimate on $\mc{N}$.
\begin{lemma}\label{N_est}
Suppose $\mc{B}$ is a convex set contained in $O^{M, m}$ and $\bm{X}_1, \bm{X}_2, \mc{P} \bm{X}_1, \mc{P} \bm{X}_2\in \mc{B}$, where $\mc{P}$ 
is the projection given in \eqref{projPPi}.
Then, for $\gamma\in(0,1), \gamma\neq 1/2$, 
\begin{align*}
\hnorm{\mc{N}(\bm{X}_1) - \mc{N}(\bm{X}_2)}{\floor{2\gamma}}{2\gamma-\floor{2\gamma}} &\leq C\frac{M^4}{m^5}\left(\hnorm{\Pi\bm{X}_1}{1}{\gamma} + \hnorm{\Pi\bm{X}_2}{1}{\gamma} \right)\hnorm{\bm{X}_1 - \bm{X}_2}{1}{\gamma},
\end{align*}
where $\Pi$ is defined in \eqref{projPPi}. If $\gamma = 1/2$, 
\begin{align*}
\hnorm{\mc{N}(\bm{X}_1) - \mc{N}(\bm{X}_2)}{0}{\alpha} &\leq C\frac{M^4}{m^5}\left(\hnorm{\Pi\bm{X}_1}{1}{\gamma} + \hnorm{\Pi\bm{X}_2}{1}{\gamma} \right)\hnorm{\bm{X}_1 - \bm{X}_2}{1}{\gamma},
\end{align*}
for any $\alpha\in (0, 1)$.
\end{lemma}

\begin{poof}
Let $\bm{X}_1, \bm{X}_2$ be as above. Then,
\begin{align}\label{NX1X2}
\mc{N}(\bm{X}_1) - \mc{N}(\bm{X}_2) &= \int_{0}^{1}\frac{d}{ds}\mc{N}(s\bm{X}_1 + (1 - s)\bm{X}_2)ds = \int_{0}^{1} \p_{\bm{X}}\mc{N}[s\bm{X}_1 + (1 - s)\bm{X}_2]\bm{Y}ds.
\end{align}
where $\bm{Y}=\bm{X}_1-\bm{X}_2$. Note that $\mc{P}\bm{X}_1, \mc{P}\bm{X}_2$ and hence $s\mc{P}\bm{X}_1+(1-s)\mc{P}\bm{X}_2$
are in $\wh{\mc{V}}$ (are circular equilibria).
Given \eqref{N0pXN0}, we have:
\begin{equation*}
\p_{\bm{X}}\mc{N}[\mc{P}(s\bm{X}_1 + (1 - s)\bm{X}_2)]\bm{Y}=0.
\end{equation*}
Note also that, for any $\bm{V}\in C^{1,\gamma}, \starnorm{\bm{V}}>0$ we have:
\begin{equation*}
\p_{\bm X}\mc{N}[\bm{V}]\bm{Y}=\p_{\bm X}\mc{R}[\bm{V}]\bm{Y}-(\mc{L}\bm{Y}-\Lambda\bm{Y}). 
\end{equation*}
In fact, the difference between $\p_{\bm X}\mc{N}$ and $\p_{\bm X}\mc{R}$ can be written in terms of the Hilbert transform as we saw in \eqref{LLambdaB}.
We may thus estimate the integrand in \eqref{NX1X2} as follows:
\begin{equation}
\begin{split}\label{NX1X2int}
&\hnorm{\p_{\bm{X}}\mc{N}[s\bm{X}_1 + (1 - s)\bm{X}_2]\bm{Y}}{\floor{2\gamma}}{2\gamma - \floor{2\gamma}}\\
=&\hnorm{\p_{\bm{X}}\mc{N}[s\bm{X}_1 + (1 - s)\bm{X}_2]\bm{Y}-\p_{\bm{X}}\mc{N}[\mc{P}(s\bm{X}_1 + (1 - s)\bm{X}_2)]\bm{Y}}{\floor{2\gamma}}{2\gamma - \floor{2\gamma}}\\
=&\hnorm{\p_{\bm{X}}\mc{R}[s\bm{X}_1 + (1 - s)\bm{X}_2]\bm{Y}-\p_{\bm{X}}\mc{R}[\mc{P}(s\bm{X}_1 + (1 - s)\bm{X}_2)]\bm{Y}}{\floor{2\gamma}}{2\gamma - \floor{2\gamma}}\\
\leq &C\frac{M^4}{m^5}\hnorm{s\Pi\bm{X}_1+(1-s)\Pi\bm{X}_2}{1}{\gamma}\hnorm{\bm{Y}}{1}{\gamma}\\
\leq &C\frac{M^4}{m^5}(s\hnorm{\Pi\bm{X}_1}{1}{\gamma}+(1-s)\hnorm{\Pi\bm{X}_2}{1}{\gamma})\hnorm{\bm{Y}}{1}{\gamma}
\end{split}
\end{equation}
In the first inequality, we used Proposition \ref{doublePertN}
and the definition of the projection $\Pi$ given in \eqref{projPPi}. 
The desired estimate is now immediate by taking the $C^{\floor{2\gamma},2\gamma-\floor{2\gamma}}$
norm on both sides of \eqref{NX1X2} and using \eqref{NX1X2int}. The second statement follows by applying the second statement of Proposition \ref{doublePertN}. 
\end{poof}

Recall that $\mc{V}$ was the four-dimensional kernel of $\mc{L}$ spanned by \eqref{zero_eigenspace}. Let $\mc{W}$ be the orthogonal complement 
of $\mc{V}$ in $C^{1,\gamma}(\mbs)$:
\begin{equation*}
\mc{W}=\lbrace\bm{w}\in C^{1,\gamma}(\mbs)|\dual{\bm{w}}{\bm{v}}=0 \text{ for all } \bm{v}\in \mc{V}\rbrace.
\end{equation*}
An equivalent definition is that $\mc{W}$ are the elements of $C^{1,\gamma}(\mbs)$ annihilated by the projection $\mc{P}$.
The subspaces $\mc{W}$ and $\mc{V}$ are both closed in $C^{1,\gamma}(\mbs)$ and are thus Banach spaces with respect to the $C^{1,\gamma}$ norm.

To motivate the calculation to follow, we perform the following formal calculation.
Consider \eqref{LNeqn}, apply the projections $\mc{P}$ and $\Pi$ and set $\bm{Y}=\mc{P}\bm{X}$ and $\bm{Z}=\Pi\bm{X}$.
We see that $\bm{Y}$ and $\bm{Z}$ should satisfy:
\begin{align*}
\p_t \bm{Y}&=\mc{L}Y+\Pi\mc{N}(\bm{Y}+\bm{Z}),\\
\p_t \bm{Z}&=\mc{P}\mc{N}(\bm{Y}+\bm{Z}).
\end{align*}
The equation for $\bm{Y}$ is a evolution equation in $\mc{W}$ whereas the equation for $\bm{Z}$ lives in $\mc{V}$.
Given the decay estimate of Proposition \ref{SGDecay} for $e^{t\mc{L}}$ acting on $\mc{W}$, we 
expect $\bm{Y}$ to decay exponentially.

Define the following function spaces:
\begin{align*}
\mc{W}_\sigma&=\lbrace\bm{Y}(t)\in C([0,\infty); \mc{W})|\norm{\bm{Y}}_\sigma<\infty\rbrace, \; \norm{\bm{Y}}_\sigma=\sup_{t\geq 0}e^{t\sigma}\norm{\bm{Y}(t)}_{C^{1,\gamma}},\;\sigma>0,\\
\mc{V}_0&=\lbrace\bm{Z}(t)\in C([0,\infty); \mc{V})\norm{\bm{Z}}_0<\infty\rbrace, \; \norm{\bm{Z}}_0=\sup_{t\geq 0}\norm{\bm{Z}(t)}_{C^{1,\gamma}}.
\end{align*}
Note that $\mc{V}$ is finite dimensional, and thus, all norms are equivalent on $\mc{V}$. 
In inequality \eqref{projected_dynamics_decay} in the statement of Theorem \ref{StabilityTheorem}, 
we used the norm, denoted by $\norm{\cdot}_{\mc{V}}$, induced by the coordinate vectors $\bm{e}_{{\rm r},{\rm t},x,y}$.
This is also the norm we use to computationally check our decay result in Section \ref{sect:numerics}. 
In the estimates to follow, however, it would be more convenient for us to use the $C^{1,\gamma}$ norm.

For a pair of functions $(\bm{Y},\bm{Z})\in \mc{W}_\sigma\times \mc{V}_0$, define the norm to be:
\begin{equation*}
\norm{(\bm{Y},\bm{Z})}_{\mc{W}_\sigma\times \mc{V}_0}=\norm{\bm{Y}}_\sigma+\norm{\bm{Z}}_0.
\end{equation*}

Let $\bm{Z}_\star\in \wh{\mc{V}}$ be a uniformly parametrized circle of radius $1$ centered at the origin. 
Define the set of functions $\mc{B}_{M_{\bm{Y}}, M_{\bm{Z}}}$ as the set of all $(\bm{Y}, \bm{Z}) \in \mc{W}_\sigma\times \mc{V}_0$ 
satisfying
\begin{align*}
\norm{\bm{Y}}_\sigma\leq M_{\bm{Y}} \text{ and } \norm{\bm{Z}-\bm{Z}_\star}_0\leq M_{\bm{Z}}
\end{align*}
where $M_{\bm{Y}}$, $M_{\bm{Z}}$ are some positive constants. In a slight abuse of notation, we let
$\bm{Z}_\star$ above denote a function of $t$ that is constant in time with value $\bm{Z}_\star$. 

\begin{proposition}\label{S_contract}
Let $0 < \sigma < 1/4$ and $\bm{Y}_0\in \mc{W}$. 
Then there exist positive constants $\rho_0, M_{\bm{Y}}$ and $M_{\bm{Z}}$ with the following properties. 
If $\norm{\bm{Y}_0}_{C^{1, \gamma}}\leq \rho_0$, the map
\begin{align*}
\mc{S}\begin{pmatrix}
\bm{Y} \\ \bm{Z}
\end{pmatrix} &= \begin{pmatrix}
e^{\mc{L}t}\bm{Y}_0 + \int_{0}^{t}e^{\mc{L}(t - s)}\Pi \mc{N}(\bm{Y}(s)+\bm{Z}(s))ds\\
\bm{Z}_\star + \int_{0}^{t}\mc{P}\mc{N}(\bm{Y}(s)+\bm{Z}(s))ds
\end{pmatrix}
\end{align*}
maps $\mc{B}_{M_{\bm{Y}}, M_{\bm{Z}}}$ to itself and is a contraction. The resulting fixed point satisfies:
\begin{equation}\label{YZfixedpointbound}
\norm{\bm{Y}}_\sigma \leq C\norm{\bm{Y}_0}_{C^{1,\gamma}}, \quad \norm{\bm{Z}-\bm{Z}_\star}_0 \leq C\norm{\bm{Y}_0}_{C^{1,\gamma}}
\end{equation}
where $C$ is a constant that only depends on $\sigma$ and $\gamma$.
\end{proposition}

\begin{poof}
First, we make $M_{\bm{Y}}$ and $M_{\bm{Z}}$ small enough so that perturbations away from $\bm{Z}_\star$ remains non-degenerate
and non-self-intersecting. 
Consider the curve $\bm{Z}_\star+\bm{W}$ where $\bm{W}\in C^{1,\gamma}(\mbs)$. We have:
\begin{equation*}
\abs{\bm{Z}_\star(\theta)+\bm{W}(\theta)-(\bm{Z}_\star(\theta')+\bm{W}(\theta'))}
\geq \abs{\bm{Z}_\star(\theta)-\bm{Z}_\star(\theta')}-\abs{\bm{W}(\theta)-\bm{W}(\theta')}
\geq (\starnorm{\bm{Z}_\star}-\snorm{\bm{W}}_{C^1})\abs{\theta-\theta'}.
\end{equation*}
Dividing through by $\abs{\theta-\theta'}$ and taking the infimum on the left hand side, we have:
\begin{equation*}
\starnorm{\bm{Z}_\star+\bm{W}}\geq \starnorm{\bm{Z}_\star}-\snorm{\bm{W}}_{C^1}\geq \frac{2}{\pi}-\norm{\bm{W}}_{C^{1,\gamma}},
\end{equation*}
where we used $\starnorm{\bm{Z}_\star}=2/\pi$. If $\norm{\bm{W}}_{C^{1,\gamma}}\leq 1/\pi$, say, we are assured that $\starnorm{\bm{Z}_\star+\bm{W}}\geq 1/\pi$.
We thus choose $M_{\bm{Y}}$ and $M_{\bm{Z}}$ small enough so that
\begin{equation*}
\norm{\bm{Y}+\bm{Z}-\bm{Z}_\star}_{C^{1,\gamma}}\leq \norm{\bm{Y}}_{C^{1,\gamma}}+\norm{\bm{Z}-\bm{Z}_\star}_{C^{1,\gamma}}\leq M_{\bm{Y}}+M_{\bm{Z}}\leq \frac{1}{\pi}.
\end{equation*}
So long as $M_{\bm{Y}}$ and $M_{\bm{Z}}$ are chosen in this way, for $(\bm{Y}(t),\bm{Z}(t))\in \mc{B}_{M_{\bm{Y}},M_{\bm{Z}}}$ we have:
\begin{equation}\label{Mmbound}
\starnorm{\bm{Y}(t)+\bm{Z}(t)}\geq \frac{1}{\pi}=m, \quad \norm{\bm{Y}(t)+\bm{Z}(t)}_{C^{1,\gamma}}\leq \norm{\bm{Z}_\star}_{C^{1,\gamma}}+\frac{1}{\pi}=M.
\end{equation}
Thus, for any $(\bm{Y},\bm{Z})\in \mc{B}_{M_{\bm{Y}},M_{\bm{Z}}}$, $\bm{X}(t)=\bm{Y}(t)+\bm{Z}(t)\in O^{M,m}$ with the constants $M$ and $m$ as above.

Let $\bm{X}_1 = \bm{Y}_1 + \bm{Z}_1$ and $\bm{X}_2 = \bm{Y}_2 + \bm{Z}_2$ with $\bm{Y}_i, \bm{Z}_i \in\mc{B}_{M_{\bm{Y}}, M_{\bm{Z}}}$ for $i = 1 , 2$. Then, 
\begin{align*}
\mc{S}\begin{pmatrix}\bm{Y}_1 \\ \bm{Z}_1\end{pmatrix} - \mc{S}\begin{pmatrix} \bm{Y}_2 \\ \bm{Z}_2 \end{pmatrix} &= \begin{pmatrix}
\int_{0}^{t}e^{\mc{L}(t - s)}\Pi\left( \mc{N}(\bm{X}_1(s)) - \mc{N}(\bm{X}_2 (s))\right)ds \\ \int_{0}^{t}\mc{P}\left( \mc{N}(\bm{X}_1 (s)) - \mc{N}(\bm{X}_2 (s))\right)ds 
\end{pmatrix} =: \begin{pmatrix} \bm{W}_1 - \bm{W}_2 \\ \bm{V}_1 - \bm{V}_2 \end{pmatrix}.
\end{align*}
We will show that $\mc{S}$ is a contraction in both components. For the first component,
\begin{align*}
\norm{\bm{W}_1 - \bm{W}_2}_{\sigma} &\leq \sup_{t \geq 0}e^{\sigma t}\int_{0}^{t} \hnorm{e^{\mc{L}(t - s)}\Pi\left( \mc{N}(\bm{X}_1(s)) - \mc{N}(\bm{X}_2(s))\right)}{1}{\gamma}ds\\
&\leq \sup_{t \geq 0}C e^{\sigma t} \int_{0}^{t} e^{-(t - s)/4}\left((t - s)^{-\gamma} + 1 \right)\hnorm{\mc{N}(\bm{X}_1(s)) - \mc{N}(\bm{X}_2(s))}{\floor{2\gamma}}{2\gamma - \floor{2\gamma}}ds,
\end{align*}
where we used Proposition \ref{SGDecay} in the second inequality.
Using Lemma \ref{N_est} and \eqref{Mmbound}, we find
\begin{align}\label{N_diff_bound}
\hnorm{\mc{N}(\bm{X}_1) - \mc{N}(\bm{X}_2)}{\floor{2\gamma}}{2\gamma - \floor{2\gamma}} &\leq C\left( \hnorm{\bm{Y}_1}{1}{\gamma} + \hnorm{\bm{Y}_2}{1}{\gamma}\right)\norm{\bm{X}_1-\bm{X}_2}_{C^{1,\gamma}},
\end{align}
But for $i = 1, 2$,
\begin{align*}
\hnorm{\bm{Y}_i(s)}{1}{\gamma} &\leq e^{-\sigma s}\norm{\bm{Y}_i}_{C^{1, \gamma}, \sigma} \leq e^{-\sigma s}M_{\bm{Y}}.
\end{align*}
Thus,
\begin{align*}
\norm{\bm{W}_1 - \bm{W}_2}_\sigma &\leq \sup_{t \geq 0} CM_{\bm{Y}}e^{\sigma t}\int_{0}^t e^{-(t - s)/4}e^{-\sigma s}\left( (t - s)^{-\gamma} + 1 \right) 
\hnorm{(\bm{X}_1 - \bm{X}_2)(s)}{1}{\gamma} ds\\
&\leq \sup_{t \geq 0} CM_{\bm{Y}}\int_{0}^t e^{-(t - s)(1/4-\sigma)}\left( (t - s)^{-\gamma} + 1 \right)ds
\paren{\norm{\bm{Y}_1-\bm{Y}_2}_\sigma+\norm{\bm{Z}_1-\bm{Z}_2}_0}.
\end{align*}
Using $\sigma<1/4$, we may bound the integral above as follows:
\begin{align*}
&\int_{0}^{t} e^{-(t - s)(1/4 - \sigma)}\left((t - s)^{-\gamma} + 1 \right)ds=\int_{0}^{t}e^{-u(1/4 - \sigma)}\left( u^{-\gamma} + 1 \right)du\\
\leq &\int_0^1 (u^{-\gamma} + 1)du+\int_1^\infty 2e^{-u(1/4 - \sigma)}du=\frac{1}{1-\gamma}+1+\frac{8}{1-4\sigma}.
\end{align*}
We thus have:
\begin{align*}
\norm{\bm{W}_1 - \bm{W}_2}_\sigma \leq CM_{\bm{Y}}\left(\norm{\bm{Y}_1 - \bm{Y}_2}_\sigma + \norm{\bm{Z}_1 - \bm{Z}_2}_0 \right).
\end{align*}
We may shrink $M_{\bm{Y}}$ as much as we wish so that 
\begin{align}\label{W_contract}
\norm{\bm{W}_1 - \bm{W}_2}_\sigma\leq \frac{1}{2}\left(\norm{\bm{Y}_1 - \bm{Y}_2}_\sigma + \norm{\bm{Z}_1 - \bm{Z}_2}_0 \right).
\end{align}
We now show that $\mc{S}$ is a contraction in the second component as well. We compute
\begin{align*}
\norm{\bm{V}_1 - \bm{V}_2}_0 &\leq \sup_{t \geq 0}\int_{0}^{t}\norm{\mc{P}\left(\mc{N}(\bm{X}_1 (s)) - \mc{N}(\bm{X}_2 (s)) \right)}_{C^{1,\gamma}}ds\\
&\leq \sup_{t \geq 0}C\int_{0}^{t}\hnorm{\mc{N}(\bm{X}_1 (s)) - \mc{N}(\bm{X}_2 (s))}{\floor{2\gamma}}{2\gamma - \floor{2\gamma}}ds,
\end{align*}
where we have used the fact that $\mc{P}$ is a bounded operator mapping to $\mc{V}$
and that all norms are equivalent on a finite dimensional space. 
Then, using Lemma \ref{N_est} gives
\begin{align*}
\norm{\bm{V}_1 - \bm{V}_2}_0 &\leq \sup_{t \geq 0} C\int_{0}^{t}\left(\hnorm{\bm{Y}_1 (s)}{1}{\gamma} + \hnorm{\bm{Y}_2 (s)}{1}{\gamma} \right)\left(\hnorm{(\bm{Y}_1 - \bm{Y}_2)(s)}{1}{\gamma} + \hnorm{(\bm{Z}_1 - \bm{Z}_2)(s)}{1}{\gamma} \right)ds\\
&\leq \sup_{t\geq 0} CM_{\bm{Y}}\int_{0}^{t}e^{-\sigma s}\left(\hnorm{(\bm{Y}_1 - \bm{Y}_2)(s)}{1}{\gamma} + \hnorm{(\bm{Z}_1 - \bm{Z}_2)(s)}{1}{\gamma}\right)ds\\
&\leq \sup_{t\geq 0}CM_{\bm{Y}} \int_{0}^{t} \paren{e^{-2s\sigma}\norm{\bm{Y}_1 - \bm{Y}_2}_\sigma + e^{-s\sigma}\norm{\bm{Z}_1 - \bm{Z}_2}_0}ds\\
&\leq CM_{\bm{Y}} \left(\norm{\bm{Y}_1 - \bm{Y}_2}_\sigma + \norm{\bm{Z}_1 - \bm{Z}_2}_0 \right).
\end{align*}
Shrinking $M_{\bm{Y}}$ again as needed, we conclude
\begin{align}\label{V_contract}
\norm{\bm{V}_1 - \bm{V}_2}_0 \leq \frac{1}{2}\left(\norm{\bm{Y}_1 - \bm{Y}_2}_\sigma + \norm{\bm{Z}_1 - \bm{Z}_2}_0 \right).
\end{align}

We now choose $\norm{\bm{Y}_0}_{C^{1,\gamma}}$ small enough so that $\mc{S}$ maps $\mc{B}_{M_{\bm{Y}},M_{\bm{Z}}}$ to itself.
\begin{align*}
\begin{pmatrix}
\bm{W} \\ \bm{V}
\end{pmatrix} = \mc{S}\begin{pmatrix}
\bm{Y}\\ \bm{Z}
\end{pmatrix} &= 
\bigg( \mc{S}\begin{pmatrix}
\bm{Y}\\ \bm{Z}
\end{pmatrix} - \mc{S}\begin{pmatrix}
0 \\ \bm{Z}
\end{pmatrix}\bigg)
 + \mc{S}\begin{pmatrix}
0 \\ \bm{Z}
\end{pmatrix}.
\end{align*}
Note that, since $\mc{N}(\bm{Z})=0$ (see \eqref{N0pXN0}), we have:
\begin{equation*}
\mc{S}\begin{pmatrix}
0 \\ \bm{Z}
\end{pmatrix}=
\begin{pmatrix}
e^{t\mc{L}}\bm{Y}_0\\ \bm{Z}_\star
\end{pmatrix}.
\end{equation*}
Equation (\ref{W_contract}) then implies that 
\begin{equation}\label{Wmapsback}
\begin{split}
\norm{\bm{W}}_\sigma &\leq \frac{1}{2}\norm{\bm{Y}}_\sigma + \sup_{t\geq 0}e^{t\sigma}\norm{e^{t\mc{L}}\bm{Y}_0}_{C^{1,\gamma}}
\leq \frac{1}{2}\norm{\bm{Y}}_\sigma+\sup_{t\geq 0} C e^{-t(1/4 - \sigma)}\hnorm{\bm{Y}_0}{1}{\gamma}\\
&\leq \frac{1}{2}\norm{\bm{Y}}_\sigma + C\hnorm{\bm{Y}_0}{1}{\gamma} \leq \frac{1}{2}M_{\bm{Y}}+C\hnorm{\bm{Y}_0}{1}{\gamma},
\end{split}
\end{equation}
where we used Lemma \ref{SGDecay} in the second inequality and the fact that $\sigma<1/4$ in the third.
Equation (\ref{V_contract}) implies
\begin{equation}\label{Vmapsback}
\norm{\bm{V} - \bm{Z}_\star}_0 \leq \frac{1}{2}\norm{\bm{Y}}_\sigma \leq  \frac{1}{2}M_{\bm{Y}}.
\end{equation}
If we thus take $\hnorm{\bm{Y}_0}{1}{\gamma}\leq \rho_0=M_{\bm Y}/2C$ and $M_{\bm Y}/2\leq M_{\bm{Z}}$,
$\mc{S}$ maps $\mc{B}_{M_{\bm{Y}},M_{\bm{Z}}}$ to itself.

Finally, let $(\bm{Y}, \bm{Z})$ be the fixed point of this map. 
Substituting $\bm{W}=\bm{Y}$ and $\bm{V}=\bm{Z}$ into \eqref{Wmapsback} and \eqref{Vmapsback} respectively gives \eqref{YZfixedpointbound}.
\end{poof}

\begin{poof}[Proof of item \ref{i:1gamma_decay} of Theorem \ref{StabilityTheorem}]
Let $(\bm{Y}(t),\bm{Z}(t))$ be the fixed point of the map $\mc{S}$ considered in Proposition \ref{S_contract}:
\begin{align}
\label{Y_is_fixedpoint}
\bm{Y}(t)&=e^{t\mc{L}}\bm{Y}_0+\int_0^t e^{(t-s)\mc{L}}\Pi \mc{N}(\bm{Y}(s)+\bm{Z}(s))ds,\\
\label{Z_is_fixedpoint}
\bm{Z}(t)&=\bm{Z}_\star+\int_0^t \mc{P}\mc{N}(\bm{Y}(s)+\bm{Z}(s))ds.
\end{align}
We continue to use the notation used in the proof of Proposition \ref{S_contract}.
From \eqref{YZfixedpointbound}, we see that:
\begin{equation}\label{Ydecayw/sigma}
\norm{\bm{Y}(t)}_{C^{1,\gamma}}\leq C\norm{\bm{Y}_0}_{C^{1,\gamma}}e^{-\sigma t}, \; 0<\sigma<1/4.
\end{equation}
We first show that we may replace the above exponential decay with $e^{-t/4}$ with a possible adjustment of the constant $C$.
Take the $C^{1,\gamma}$ norm on both sides of \eqref{Y_is_fixedpoint}.
\begin{align*}
\norm{\bm{Y}(t)}_{C^{1,\gamma}}&\leq \norm{e^{t\mc{L}}\bm{Y}_0}_{C^{1,\gamma}}
+\int_0^t \norm{e^{(t-s)\mc{L}}\Pi \mc{N}(\bm{Y}(s)+\bm{Z}(s))}_{C^{1,\gamma}}ds\\
&\leq C\norm{\bm{Y}_0}_{C^{1,\gamma}}e^{-t/4}
+\int_0^t Ce^{-(t-s)/4}\paren{\frac{1}{(t-s)^{\gamma}}+1}\hnorm{\mc{N}(\bm{Y}(s)+\bm{Z}(s))}{\floor{2\gamma}}{2\gamma-\floor{2\gamma}}ds.
\end{align*}
where we used Lemma \eqref{SGDecay} in the second inequality.
Observe that, using Lemma \ref{N_est}, we have the following estimate:
\begin{equation}
\begin{aligned}\label{NYZ}
\hnorm{\mc{N}(\bm{Y}(t)+\bm{Z}(t))}{\floor{2\gamma}}{2\gamma - \floor{2\gamma}}&= 
\hnorm{\mc{N}(\bm{Y}(t)+\bm{Z}(t)) - \mc{N}(\bm{Z}(t))}{\floor{2\gamma}}{2\gamma - \floor{2\gamma}} \\
&\leq C\hnorm{\bm{Y}(t)}{1}{\gamma}^2
\leq C\norm{\bm{Y}_0}_{C^{1,\gamma}}^2 e^{-2\sigma t},  
\end{aligned}
\end{equation}
where we used \eqref{N0pXN0} in the equality and \eqref{Ydecayw/sigma} in the last inequality. 
We thus obtain:
\begin{equation*}
\norm{\bm{Y}(t)}_{C^{1,\gamma}}
\leq C\norm{\bm{Y}_0}_{C^{1,\gamma}}e^{-t/4}+C\int_0^te^{-(t-s)/4}e^{-2\sigma s}\paren{\frac{1}{(t-s)^{\gamma}}+1}ds\norm{\bm{Y}_0}_{C^{1,\gamma}}^2.
\end{equation*}
If we take $\sigma>1/8$, for $t\geq 1$,
\begin{equation}\label{intboundsigma>1/8}
\begin{split}
&\int_0^te^{-(t-s)/4}e^{-2\sigma s}\paren{\frac{1}{(t-s)^{\gamma}}+1}ds=
e^{-2\sigma t}\int_0^t e^{(2\sigma-1/4)u}\paren{\frac{1}{u^\gamma}+1}du\\
\leq&e^{-2\sigma t}\paren{e^{(2\sigma-1/4)}\int_0^1\paren{\frac{1}{u^\gamma}+1}du+2\int_1^t e^{(2\sigma-1/4)u}du}\\
=&e^{-2\sigma t}\paren{e^{(2\sigma-1/4)}\paren{\frac{1}{1-\gamma}+1}+\frac{8}{8\sigma-1}\paren{e^{(2\sigma-1/4)t}-1}}\leq Ce^{-t/4}.
\end{split}
\end{equation}
For $t<1$, the second integral in the second line above is not needed.
We thus have,
\begin{equation}\label{Y1/4decay}
\norm{\bm{Y}(t)}_{C^{1,\gamma}}
\leq C(\norm{\bm{Y}_0}_{C^{1,\gamma}}+\norm{\bm{Y}_0}_{C^{1,\gamma}}^2)e^{-t/4}\leq C(1+\rho_0)\norm{\bm{Y}_0}_{C^{1,\gamma}}e^{-t/4},
\end{equation}
where we used the assumption on $\norm{\bm{Y}_0}_{C^{1,\gamma}}$ stated in Proposition \ref{S_contract} in the last inequality.

We now turn to the component $\bm{Z}$. Take the norm on both sides of \eqref{Z_is_fixedpoint}. 
\begin{equation}\label{Z_finite}
\begin{split}
\norm{\bm{Z}(t)}_{C^{1,\gamma}} &\leq \norm{\bm{Z}_\star} + \int_{0}^{t}\hnorm{\mc{N}(\bm{Y}(s), \bm{Z}(s))}{\floor{2\gamma}}{2\gamma - \floor{2\gamma}}ds\\
&\leq \norm{\bm{Z}_\star} + C\norm{\bm{Y}_0}_{C^{1,\gamma}}^2\int_{0}^{t}e^{-s/2}ds\leq \norm{\bm{Z}_\star}_{C^{1,\gamma}}+2C\norm{\bm{Y}_0}_{C^{1,\gamma}}^2,
\end{split}
\end{equation}
where we used \eqref{Y1/4decay} in the second inequality.
Thus, the following is well-defined.
\begin{align*}
\bm{Z}_\infty &:= \bm{Z}_\star + \int_{0}^{\infty}\mc{P}\mc{N}(\bm{Y}(s), \bm{Z}(s))ds.
\end{align*}
It is clear that $\bm{Z}_\infty\in \mc{V}$ and since it does not degenerate to a point ($\starnorm{\bm{Z}(t)}\geq m=1/\pi$ by construction, see \eqref{Mmbound})
$\bm{Z}_\infty\in \wh{\mc{V}}$ is a uniformly parametrized circle. 
We have:
\begin{equation}\label{Z1/2decay}
\begin{split}
\norm{\bm{Z}(t) - \bm{Z}_\infty}_{C^{1,\gamma}} &\leq \int_{t}^{\infty}\hnorm{\mc{N}(\bm{Y}(s), \bm{Z}(s))}{\floor{2\gamma}}{2\gamma - \floor{2\gamma}}ds\\
&\leq C\norm{\bm{Y}_0}_{C^{1,\gamma}}^2\int_{t}^{\infty}e^{-s/2}ds=2C\norm{\bm{Y}_0}_{C^{1,\gamma}}^2e^{-t/2}.
\end{split}
\end{equation}
Thus, $\bm{Z}$ converges to $\bm{Z}_\infty$ exponentially with the above rate.

We must still show that $\bm{X}=\bm{Y}+\bm{Z}$ is indeed a solution to the Peskin problem. 
Suppose $\bm{X}_0\in h^{1,\gamma}(\mbs)$ is our initial data with $\Pi\bm{X}_0=\bm{Y}_0$ and $\mc{P}\bm{X}_0=\bm{Z}_0$.
We first assume that $\bm{Z}_0$ is a uniformly parametrized unit circle centered at the origin, $\bm{Z}_0=\bm{Z}_{\star}$. 
This restriction will be later lifted with a scaling argument.
From \eqref{Y_is_fixedpoint} and \eqref{Z_is_fixedpoint}, by a simple adaptation 
of the proof of Lemma \ref{mild_is_strong}, we see that 
$\bm{Y}\in C([0, T); \mc{W})\cap C^1((0, T); C^{\gamma}(\mbs))$ for any $T>0$ 
($\bm{Z}$ lives in a finite dimensional space so its differentiability is automatic)
and that $\bm{Y}$ and $\bm{Z}$ satisfy the strong form of the equations for $t>0$
\begin{align*}
\p_{t}\bm{Y} &= \mc{L}\bm{Y} + \Pi\mc{N}(\bm{Y}+\bm{Z}),\\
\p_{t}\bm{Z} &=  \mc{P}\mc{N}(\bm{Y}+\bm{Z}),
\end{align*}
with $(\bm{Y}(t),\bm{Z}(t))\to (\bm{Y}_0,\bm{Z}_0)$ in the $C^{1,\gamma}$ norm.
Adding these two equations together and letting $\bm{X} = \bm{Y} + \bm{Z}$, we see that $\bm{X}$ satisfies:
\begin{align*}
\p_t \bm{X} &= \mc{L}\bm{X} + \mc{N}(\bm{X}), \; \bm{Y}=\Pi \bm{X}, \; \bm{Z}=\mc{P}\bm{X},\\
\bm{X}(t) &\to \bm{X}_0 \text{ in } C^{1,\gamma}(\mbs).
\end{align*}
Recall from \eqref{LNeqn} that 
\begin{align*}
\mc{L}\bm{X} + \mc{N}(\bm{X}) &= \mc{R}(\bm{X}) = \Lambda\bm{X} + \mc{R}(\bm{X}).
\end{align*} 
We thus see that $\bm{X}=\bm{Y}+\bm{Z}$ is in fact a strong solution of the Peskin problem, 
and is consequently also a mild solution of the Peskin problem thanks to Lemma \ref{strong_is_mild}.
By the uniqueness result for mild solutions (Theorem \ref{LWPTheorem}), $\bm{X}=\bm{Y}+\bm{Z}$ is {\em the} mild solution to the 
Peskin initial value problem with $\bm{X}=\bm{X}_0$. 

We finally lift the restriction that $\mc{P}\bm{X}_0=\bm{Z}_0$ is a uniformly parametrized unit circle centered at the origin.
Take any initial data $\bm{X}_0\in h^{1,\gamma}(\mbs)$ and its mild solution $\bm{X}(t)$. Let $R$ be the 
radius of the circle $\bm{Z}_0$. Then, if we set:
\begin{equation}\label{Xnormalized}
\wh{\bm{X}}_0=\frac{1}{R}\paren{\bm{X}_0-\bm{p}_{\rm ctr}}, \; \wh{\bm{X}}(t)=\frac{1}{R}\paren{\bm{X}-\bm{p}_{\rm ctr}},
\quad \bm{p}_{\rm ctr}=\mc{P}_{\rm trl}\bm{X}_0,
\end{equation}
where $\mc{P}_{\rm trl}$ was defined in \eqref{proj_trl} and $\bm{p}_{\rm ctr}$ is simply the center point of the circle $\mc{P}\bm{X}_0$.
By dilation and translation invariance, $\wh{\bm{X}}(t)$ is a solution to the Peskin problem with initial data $\wh{\bm{X}}_0$. 
By design, $\mc{P}\wh{\bm{X}_0}$
is a uniformly parametrized circle centered at the origin. 
We may thus apply the results \eqref{Y1/4decay} and \eqref{Z1/2decay}
to obtain the estimate:
\begin{equation}\label{whYwhZestimates}
\norm{\wh{\bm{Y}}(t)}_{C^{1,\gamma}}\leq C\norm{\wh{\bm{Y}}_0}_{C^{1,\gamma}}e^{-t/4}, 
\quad \norm{\wh{\bm{Z}}(t)-\wh{\bm{Z}}_\infty}_{C^{1,\gamma}}\leq C\norm{\wh{\bm{Y}}_0}_{C^{1,\gamma}}^2e^{-t/2}
\text{ if } \norm{\wh{\bm{Y}}_0}_{C^{1,\gamma}}\leq \rho_0,
\end{equation}
where $\wh{\bm{Y}}(t)=\Pi\wh{\bm{X}}(t), \wh{\bm{Z}}(t)=\mc{P}\wh{\bm{X}}(t)$, $\wh{\bm{Y}}_0=\Pi\wh{\bm{X}}_0, \wh{\bm{Z}}_0=\mc{P}\wh{\bm{X}_0}$
and $\wh{\bm{Z}}_\infty$ is the point to which $\wh{\bm{Z}}(t)$ converges (whose existence was guaranteed above).
From \eqref{Xnormalized}, we see that
\begin{equation*}
\wh{\bm{Y}}(t)=\frac{1}{R}\bm{Y}(t), \; \wh{\bm{Y}}_0=\frac{1}{R}\bm{Y}_0, \; \wh{\bm{Z}}(t)=\frac{1}{R}\paren{\bm{Z}(t)-\bm{p}_{\rm ctr}}.
\end{equation*}
By plugging in the above into \eqref{whYwhZestimates}, we obtain the inequalities 
\eqref{difference_to_circle_decay} and \eqref{projected_dynamics_decay} 
by setting $\bm{Z}_\infty=R\wh{\bm{Z}}_\infty+\bm{p}_{\rm ctr}$.

Finally, note that 
\begin{equation*}
\norm{\bm{X}-\bm{Z}_\infty}_{C^{1,\gamma}}=
\norm{\Pi\bm{X}+\mc{P}\bm{X}-\bm{Z}_\infty}_{C^{1,\gamma}}\leq \norm{\Pi\bm{X}}_{C^{1,\gamma}}+\norm{\mc{P}\bm{X}-\bm{Z}_\infty}_{C^{1,\gamma}}.
\end{equation*}
Inequality  \eqref{e:XdecaytoZhigher} is a direct consequence of this and \eqref{difference_to_circle_decay} and \eqref{projected_dynamics_decay}.
\end{poof}

To obtain exponential decay in higher H\"older norms, we state a result that is a direct consequence
of Lemma \ref{l:CarryRegularity} and the proof of Proposition \ref{p:spatialsmooth}. 
We omit the proof since it will be an almost exact reiteration of the proof of 
Proposition \ref{p:spatialsmooth}.
\begin{lemma}\label{Cnboundlemma}
Suppose $\bm{X}(t)$ is a mild solution to the Peskin problem up to $t=\epsilon$ satisfying 
\begin{equation*}
\norm{\bm{X}(t)}_{C^{1,\gamma}}\leq M \text{ and } \starnorm{{\bm{X}}(t)}\geq m \text{ for } 0\leq t\leq \epsilon.
\end{equation*}
Then, for any $n\in \mathbb{N}$ and $0<\alpha<1$, we have:
\begin{equation*}
\norm{\bm{X}(\epsilon)}_{C^{n,\alpha}}\leq M_{n,\alpha}
\end{equation*}
where $M_{n,\alpha}$ is a constant that depends only on $n,\alpha,\epsilon, \gamma, M$ and $m$.
\end{lemma}

\begin{poof}[Proof of item \ref{i:higher_norm_decay} of Theorem \ref{StabilityTheorem}]
Like the above proof of item \ref{i:1gamma_decay} in Theorem \ref{StabilityTheorem}, we first assume that 
$\mc{P}\bm{X}_0$ is a unit circle centered at the origin and  later use a scaling argument to obtain the result for general initial data.
First, recall from \eqref{Mmbound} that 
\begin{align*}
\sup_{t\geq 0} \hnorm{\bm{X}(t)}{1}{\gamma} &= \norm{\bm{Z}_0}_{C^{1,\gamma}}+\frac{1}{\pi} < \infty ,\\
\inf_{t\geq 0} \starnorm{\bm{X}(t)} &\geq \frac{1}{\pi} > 0.
\end{align*}
Applying Lemma \eqref{Cnboundlemma} to the solution $\bm{X}(t+\tau), t\geq 0$ ($\tau\geq 0$ considered a parameter),
\begin{equation*}
\norm{\bm{X}(\tau+\epsilon)}_{C^{n,\alpha}}\leq M_{n,\alpha} \text{ for any } \tau\geq 0,
\end{equation*}
where $M_n$ depends only on $n,\alpha, \epsilon$ and $\gamma$. Thus, 
\begin{equation*}
\sup_{\tau\geq \epsilon} \norm{\bm{X}(t)}_{C^{n,\alpha}}\leq M_{n,\alpha}, \; n\in \mathbb{N}, 0<\alpha<1.
\end{equation*}
Lemma \ref{l:CarryRegularity} implies that:
\begin{equation*}
\sup_{\tau\geq \epsilon} \norm{\mc{R}(\bm{X}(\tau))}_{C^{n,\alpha}}\leq M^{\mc R}_{n,\alpha}, \; n\in\mathbb{N}, 0<\alpha<1.
\end{equation*}
where the constant  $M^{\mc R}_{n,\alpha}$ again depends only on $n,\alpha,\epsilon$ and $\gamma$.
Recall that:
\begin{equation*}
\mc{N}(\bm{X})=-\mc{Q}\bm{X}+\mc{R}(\bm{X}), \; \mc{Q}\bm{X}=\mc{L}\bm{X}-\Lambda \bm{X}.
\end{equation*}
Using the fact that $\mc{Q}$ is a bounded operator from $C^{n,\alpha}(\mbs)$ to itself,
(see \eqref{LLambdaB}; the Hilbert transform is bounded from $C^{n,\alpha}$ to itself as long as $0<\alpha<1$),
we have:
\begin{equation}\label{Nnalphabound}
\norm{\mc{N}(\bm{X}(t))}_{C^{n,\alpha}}\leq M^{\mc{N}}_{n,\alpha} \text{ for } t\geq \epsilon.
\end{equation}
where the above constant depends only on $n,\alpha,\epsilon$ and $\gamma$.
On the other hand, combining \eqref{NYZ} and \eqref{Y1/4decay}, we have:
\begin{equation}\label{N2gammabound}
\hnorm{\mc{N}(\bm{X}(t))}{\floor{2\gamma}}{2\gamma - \floor{2\gamma}}\leq C\norm{\bm{Y}_0}_{C^{1,\gamma}}^2 e^{-t/2},
\end{equation}
Interpolating the two estimates \eqref{Nnalphabound} and \eqref{N2gammabound}, for any $0<\sigma<1/4$, we obtain 
(see Chapter 1 of \cite{lunardi} for results on interpolation in H\"older spaces):
\begin{equation}\label{Nbetabound}
\hnorm{\mc{N}(\bm{X}(t))}{\floor{\beta}}{\beta-\floor{\beta}}\leq C\norm{\bm{Y}_0}_{C^{1,\gamma}}^{8\sigma}e^{-2\sigma t}, \text{ for } t\geq \epsilon, \; 
\beta=8\gamma\sigma +(n+\alpha)(1-4\sigma)
\end{equation}
where the above constant $C$ depends only on $n,\alpha,\epsilon, \gamma$ and $\sigma$. Since $n$ can be made arbitrarily 
large, the above estimate is true for any $\beta\geq 1+\gamma$ and  $0<\sigma<1/4$.

Suppose we know that:
\begin{equation}\label{Y1/4betadecay}
\hnorm{\bm{Y}(t)}{\floor{\beta}}{\beta-\floor{\beta}}\leq   C\norm{\bm{Y}_0}_{C^{1,\gamma}}e^{-t/4}, \; t\geq k\epsilon, \; \text{ for some } k\in \mathbb{N}
\end{equation}
where the above constant $C$ that depends only on $\beta,\gamma$ and $k\epsilon$. The above is true if $\beta=1+\gamma$ and $k=1$
by \eqref{Y1/4decay}. Let $\beta<\kappa<1+\beta$ and take the $C^{\floor{\kappa},\kappa-\floor{\kappa}}$ norm on both sides of \eqref{Y_is_fixedpoint}.
We have:
\begin{align*}
&\hnorm{\bm{Y}(t+k\epsilon)}{\floor{\kappa}}{\kappa-\floor{\kappa}}\leq \hnorm{e^{t\mc{L}}\bm{Y}(k\epsilon)}{\floor{\kappa}}{\kappa-\floor{\kappa}}
+\int_{0}^t \hnorm{e^{(t-s)\mc{L}}\Pi\mc{N}(\bm{X}(s+k\epsilon))}{\floor{\kappa}}{\kappa-\floor{\kappa}}ds\\
\leq &\frac{C}{t^{\kappa-\beta}}e^{-t/4}\hnorm{\bm{Y}(k\epsilon)}{\floor{\beta}}{\beta-\floor{\beta}}
+\int_{0}^t Ce^{-t/4}\paren{\frac{1}{(t-s)^{\kappa-\beta}}+1}\hnorm{\mc{N}(\bm{X}(s+k\epsilon))}{\floor{\beta}}{\beta-\floor{\beta}}ds\\
\leq &\frac{C}{t^{\kappa-\beta}}e^{-t/4}\hnorm{\bm{Y}_0}{1}{\gamma}
+\int_{0}^t Ce^{-(t-s)/4}\paren{\frac{1}{(t-s)^{\kappa-\beta}}+1}\norm{\bm{Y}_0}_{C^{1,\gamma}}^{8\sigma}e^{-2\sigma t}ds,
\end{align*}
where we used Lemma \ref{SGDecay} in the second inequality and \eqref{Nbetabound} and \eqref{Y1/4betadecay} in the last inequality.
Letting $\sigma>1/8$, we see from \eqref{intboundsigma>1/8} that:
\begin{equation*}
\begin{split}
\hnorm{\bm{Y}(t+k\epsilon)}{\floor{\kappa}}{\kappa-\floor{\kappa}}&\leq 
C\paren{\frac{1}{t^{\kappa-\beta}}\hnorm{\bm{Y}_0}{1}{\gamma} +\norm{\bm{Y}_0}_{C^{1,\gamma}}^{8\sigma}}e^{-t/4}\\
&\leq C\paren{\frac{1}{\epsilon^{\kappa-\beta}}+\rho_0^{8\sigma-1}}\norm{\bm{Y}_0}_{C^{1,\gamma}}e^{-t/4},
\end{split}
\end{equation*}
where we used the assumption that $\norm{\bm{Y}_0}\leq \rho_0$ (see statement of Proposition \ref{S_contract}) and $8\sigma-1>0$
in the last inequality. Therefore, 
\begin{equation*}
\hnorm{\bm{Y}(t)}{\floor{\kappa}}{\kappa-\floor{\kappa}}\leq   C\norm{\bm{Y}_0}_{C^{1,\gamma}}e^{-t/4}, \; t\geq (k+1)\epsilon,
\end{equation*}
where the above constant $C$ that depends only on $\kappa,\gamma$ and $(k+1)\epsilon$. 
Starting with $\beta=1+\gamma$, we may iterate this process indefinitely to find that the above estimate is true for any $\kappa$ with a suitably large $k$.
Since $\epsilon$ can be taken arbitrarily small and $k$ is finite, we may replace $k\epsilon$ with $\epsilon$ in the above by making the constant $C$
larger if necessary, to obtain:
\begin{equation}
\hnorm{\bm{Y}(t)}{\floor{\kappa}}{\kappa-\floor{\kappa}}\leq   C\norm{\bm{Y}_0}_{C^{1,\gamma}}e^{-t/4}, \; t\geq\epsilon,
\end{equation}
where $C$ depends only on $\gamma,\kappa$ and $\epsilon$. To obtain the same bound for the $C^n$ norm, as stated in 
\eqref{higher_norm_decay}, simply 
note that the $C^{n,\alpha}$ norm for $0<\alpha<1$ dominates the $C^n$ norm.

Finally, to obtain \eqref{higher_norm_decay} for general initial conditions,
we may use the same scaling argument used at the end of the proof of item \ref{i:1gamma_decay} of Theorem \ref{StabilityTheorem} above.
Note also that 
\begin{equation*}
\norm{\bm{X}-\bm{Z}_\infty}_{C^{n,\gamma}}\leq 
\norm{\Pi\bm{X}}_{C^{n,\gamma}}+\norm{\mc{P}\bm{X}-\bm{Z}_\infty}_{C^{n,\gamma}}
\leq \norm{\Pi\bm{X}}_{C^{n,\gamma}}+C\norm{\mc{P}\bm{X}-\bm{Z}_\infty}_{\mc{V}}
\end{equation*}
where the last inequality follows from the equivalence of all norms in finite dimensional spaces. Thus, \eqref{e:XdecaytoZhigher} 
is thus a direct consequence of \eqref{projected_dynamics_decay} and \eqref{higher_norm_decay}.
\end{poof}
\subsection{Computational Verification}\label{sect:numerics}

Here, we computationally verify the exponential decay rate to the circle.
We developed a numerical scheme to simulate the Peskin problem based on the 
small scale decomposition \eqref{SSD}. The numerical scheme is second order 
accurate in time $t$ and spectrally accurate in $\theta$. We point out that the 
second order accuracy in time (as opposed to a first-order scheme) 
turned out to be crucial in computationally verifying the asymptotic decay rate. 

We first give a description of the numerical scheme. We use equation \eqref{SSD}
and \eqref{mild_soln} as the basis for our algorithm. Discretize $\mbs$ with $N$
points so that $Nh=2\pi$, where $h$ is the grid spacing.
Let $\theta=\theta_k=k\Delta \theta, k=0,1,\cdots N-1$ be the grid locations, and $\bm{X}_{h,k}=(X_{h,k},Y_{h,k})$
be the numerically approximated value of $\bm{X}(\theta_k)$. We let $N$ be even.
For a function $W$ defined on the discrete $\theta$ grid, define $\mc{F}_h$ to be the discrete Fourier transform:
\begin{equation*}
(\mc{F}_h W)_k=\sum_{l=0}^{N-1} \exp(-2\pi ikl/N)W_l, \; k=-N/2+1, \cdots, -1, 0, 1, \cdots N/2.
\end{equation*}
Define the approximation to the derivative $\mc{D}_h$ and the semigroup $\mc{S}_h(t)$ as follows: 
\begin{align*}
\mc{D}_h W&=\mc{F}_h^{-1}\wh{\mc{D}}_h(k)\mc{F}_h W, \;
\wh{\mc{D}}_h(k)=
\begin{cases}
ik &\text{ if } k\neq N/2,\\
0 &\text{ if } k=N/2,
\end{cases}\\
\mc{S}_h(t) W&=\mc{F}_h^{-1}\wh{\mc{S}}_h(t,k)\mc{F}_h W, \;
\wh{\mc{S}}_h(t,k)=
\begin{cases}
e^{-t\abs{k}/4} &\text{ if } k\neq N/2,\\
0 &\text{ if } k=N/2.
\end{cases}
\end{align*}
Recall from \eqref{SSD} that $\mc{R}$ can be written as:
\begin{equation*}
\mc{R}=-\frac{1}{4\pi}\int_{\mbs}\p_\theta'\paren{-\log\paren{\frac{\abs{\Delta \bm{X}}}{2\abs{\sin((\theta-\theta')/2)}}}I
+\frac{\Delta \bm{X}\otimes \Delta \bm{X}}{\abs{\Delta \bm{X}}^2}}\p_\theta'\bm{X}'d\theta'.
\end{equation*}
We now approximate $\mc{R}$. Define:
\begin{align*}
H_{kl}&=-\log\paren{\frac{\abs{\bm{X}_{h,k}-\bm{X}_{h,l}}}{2\abs{\sin((\theta_k-\theta_l)/2)}}}I
+\frac{(\bm{X}_{h,k}-\bm{X}_{h,l})\otimes(\bm{X}_{h,k}-\bm{X}_{h,l})}{\abs{\bm{X}_{h,k}-\bm{X}_{h,l}}^2} \text{ for } l\neq k, \; \\
H_{kk}&=-\log\abs{(\mc{D}_h\bm{X}_h)_k}I+\frac{(\mc{D}_h\bm{X}_h)_k\otimes (\mc{D}_h\bm{X}_h)_k}{\abs{(\mc{D}_h\bm{X}_h)_k}^2}.
\end{align*}
We let the approximation of $\mc{R}$ at $\theta=\theta_k$ be:
\begin{equation*}
\mc{R}_{h,k}(\bm{X}_h)=-\frac{1}{4\pi}\sum_{l=0}^{N-1} (\mc{D}_{h,j}H_{kj})_l(\mc{D}_h \bm{X}_h)_l h
\end{equation*}
where $\mc{D}_{h,j}$ means that the operator $\mc{D}_h$ is acting on the grid function with argument $j$.

Let $t_n=n\triangle t$ where $\triangle t$ is the time step. 
We use the following Runge-Kutta type approximation scheme.
Before we describe our time-stepping scheme, note from \eqref{mild_soln} that:
\begin{equation*}
\begin{split}
\bm{X}(t_{n+1})=&e^{\triangle t\Lambda}\bm{X}(t_n)+\int_{t_n}^{t_{n+1}}e^{(t_{n+1}-s)\Lambda}\mc{R}(\bm{X}(s))ds\\
\approx&e^{\triangle t\Lambda}\bm{X}(t_n)+e^{\triangle t \Lambda/2}\mc{R}(\bm{X}(t_{n+1/2}))\triangle t, \quad t_{n+1/2}=t_n+\triangle t/2.
\end{split}
\end{equation*}
In order to approximate $\bm{X}(t_{n+1/2})$, we just use:
\begin{equation*}
\bm{X}(t_{n+1/2})\approx e^{\triangle t\Lambda/2}\bm{X}(t_n)+e^{\triangle t \Lambda/2}\mc{R}(\bm{X}(t_n))\triangle t/2.
\end{equation*}
We use the spatially discrete version of the above for our time-stepping.
Let $\bm{X}_h^n$ be the discrete approximation to $\bm{X}$ at time $t=t_n=n\triangle t$.
We let:
\begin{align*}
\bm{X}_h^{n+1/2}&=\mc{S}_h(\triangle t/2)\paren{\bm{X}_h^n+\mc{R}_h(\bm{X}_h^n)\triangle t/2},\\
\bm{X}_h^{n+1}&=\mc{S}_h(\triangle t) \bm{X}_h^n+\mc{S}_h(\triangle t/2)\mc{R}_h(\bm{X}_h^{n+1/2})\triangle t.
\end{align*}
This concludes our description of the numerical scheme.

To define the discrete projection operators, define the following discrete inner product for grid functions $\bm{V}$ and $\bm{W}$:
\begin{equation*}
\dual{\bm{V}}{\bm{W}}_h=\sum_{k=0}^{N-1} (\bm{V}_k\cdot \bm{W}_k)h.
\end{equation*}
Let $\bm{e}_{x,h}, \bm{e}_{y,h}, \bm{e}_{{\rm r},h}$ and $\bm{e}_{{\rm t},h}$ simply be the evaluation of the vectors 
$\bm{e}_{x,y}, \bm{e}_{\rm r,t}$ in \eqref{circular_equilibria} evaluated at the grid points $\theta_k$. Define the discrete 
projection operators:
\begin{equation*}
\mc{P}_h \bm{V}=\frac{1}{2\pi}\paren{\dual{\bm{V}}{\bm{e}_{x,h}}_h\bm{e}_{x,h}+\dual{\bm{V}}{\bm{e}_{y,h}}_h\bm{e}_{y,h}+\dual{\bm{V}}{\bm{e}_{{\rm r},h}}_h\bm{e}_{{\rm r},h}
+\dual{\bm{V}}{\bm{e}_{{\rm t},h}}_h\bm{e}_{{\rm t},h}}, \; \Pi_h\bm{V}=\bm{V}-\mc{P}_h\bm{V}.
\end{equation*}
Let us also define the discrete $C^1$ norm as follows:
\begin{equation*}
\norm{\bm{V}}_{C^1_h}=\sup_k\abs{\bm{V}_k}+\sup_k \abs{(\mc{D}_h\bm{V})_k}.
\end{equation*}

\begin{figure}
\begin{center}
\includegraphics[width=\textwidth]{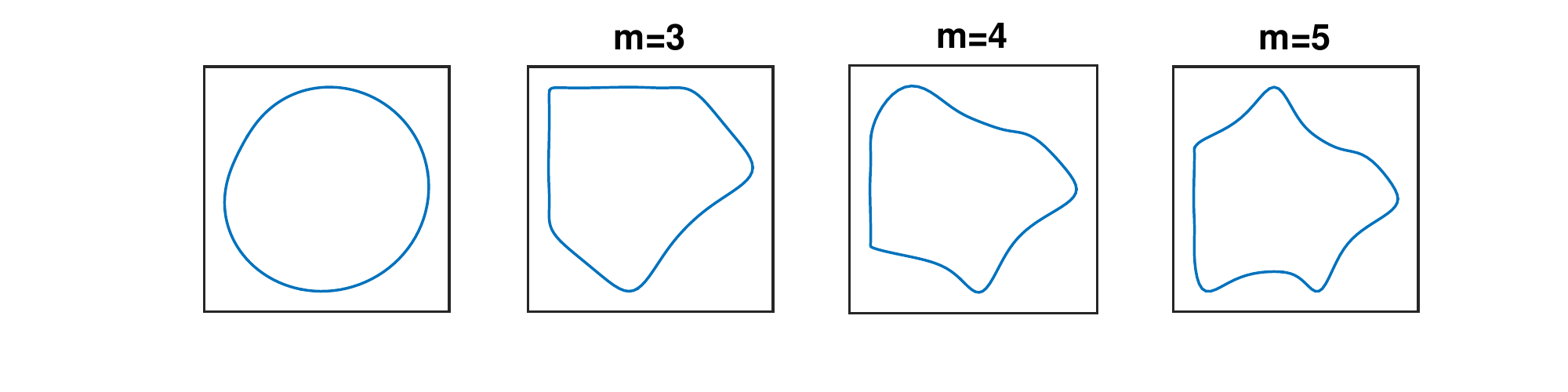}
\end{center}
\caption{\label{fig:init_curves}Initial Curves. The unlabeled curve corresponds to equation \eqref{unlabeled_curve}, where as the 
other three correspond to $m=3,4,5$ respectively in equation \eqref{labeled_curve}.}
\end{figure}

To numerically compute the decay rate to the circle, we take four different initial conditions. They are:
\begin{equation}\label{unlabeled_curve}
\bm{X}_0=
\begin{pmatrix} 
\cos(\theta)+\cos(2\theta)/5-\sin(2\theta)/10\\
\sin(\theta)+\sin(2\theta)/5+\cos(2\theta)/10
\end{pmatrix}
\end{equation} 
and 
\begin{equation}\label{labeled_curve}
\bm{X}_0=
\begin{pmatrix}
(1+\exp(\cos(3\theta))/4)\cos(\theta)\\
(1+\exp(\sin(m\theta))/4)\sin(\theta)
\end{pmatrix}, \quad m=3,4,5.
\end{equation}
The configurations of the above initial curves can be found in Figure \ref{fig:init_curves}.
The initial curve \eqref{unlabeled_curve} corresponds to a perturbation of the unit circle proportional to 
the primary decay mode \eqref{primary_decay_mode}.

\begin{figure}
\begin{center}
\includegraphics[width=0.8\textwidth]{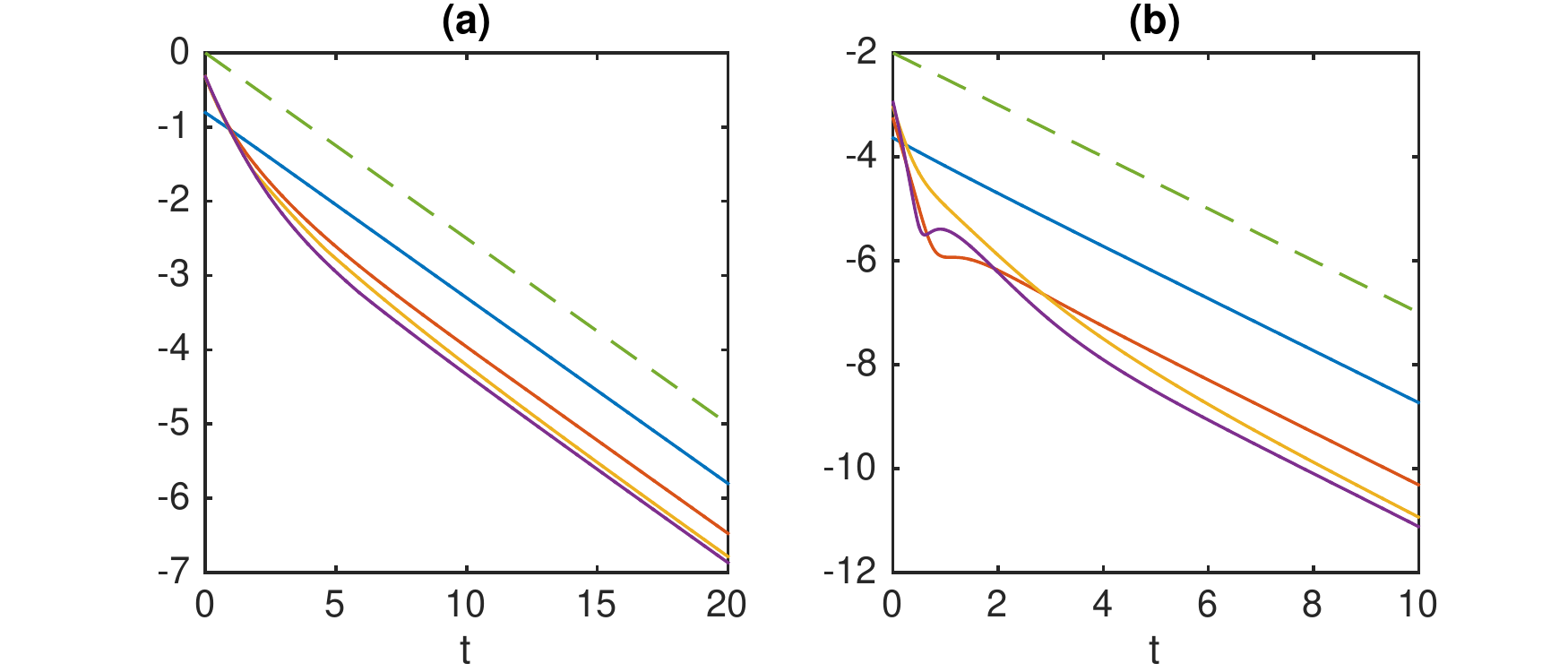}
\end{center}
\caption{\label{fig:decay_rate} In Figure (a), $\log(\norm{\Pi_h\bm{X}_h^n}_{C^1_h})$ is plotted against $t=n\triangle t$ for 
the four different initial data given in \eqref{unlabeled_curve} and \eqref{labeled_curve}. 
The dashed line has a slope $-1/4$, indicating the theoretical decay rate. 
In Figure (b), $\log(\abs{(\mc{D}_t\bm{a}_h)^{n+1/2}})$ is plotted against $t=(n+1/2)\triangle t$ with the four 
different initial data. The dashed line has a slope of $-1/2$, the theoretical decay rate.
The solid curves that are almost straight in both figures correspond to the initial data \eqref{unlabeled_curve}.
Note that we have plotted $\log(\norm{\Pi_h\bm{X}_h^n}_{C^1_h})$ up to $t=20$ whereas 
$\log(\abs{(\mc{D}_t\bm{a}_h)^{n+1/2}})$ is plotted up to $t=10$, in accordance with the fact that 
the latter decays twice as fast as the former.
}
\end{figure}

We simulated the dynamics with the above initial data with $N=128$ and $\triangle t=0.01$.
Let us recall from Theorem \ref{StabilityTheorem} that $\Pi\bm{X}$ decays to $0$ at an exponential rate 
of $e^{-t/4}$ and $\mc{P}\bm{X}$ decays to some uniformly parametrized circle at an exponential rate 
of $e^{-t/2}$. 
To computationally verify the decay result for $\Pi\bm{X}$, we compute
\begin{equation}\label{PiXC1h}
\norm{\Pi_h\bm{X}_h^n}_{C^1_h}.
\end{equation}
For $\mc{P}\bm{X}$, the circle to which $\bm{X}$ converges is unknown and thus
we instead compute the decay of the time derivative.  For $\mc{P}\bm{X}$, we have:
\begin{equation*}
\mc{P}\bm{X}(t)=a_x(t)\bm{e}_x+a_y(t)\bm{e}_y+a_{\rm r}(t)\bm{e}_{\rm r}+a_{\rm t}(t)\bm{e}_{\rm t}.
\end{equation*}
The proof of Theorem \ref{StabilityTheorem}  implies that 
\begin{equation*}
\abs{d\bm{a}/dt}, \; \bm{a}=(a_x,a_y,a_{\rm r}, a_{\rm t})^{\rm T} \text{ decays like } e^{-t/2},
\end{equation*}
where the absolute value above is the Euclidean norm in $\mathbb{R}^4$.
Let:
\begin{equation*}
\mc{P}_h\bm{X}_h^n=a_{x,h}^n\bm{e}_{x,h}+a_{y,h}^n\bm{e}_{y,h}+a_{{\rm r},h}^n\bm{e}_{{\rm r},h}+a_{{\rm t},h}^n\bm{e}_{{\rm t},h}.
\end{equation*}
Set:
\begin{equation*}
(\mc{D}_t\bm{a}_h)^{n+1/2}=\frac{\bm{a}_h^{n+1}-\bm{a}_h^n}{\triangle t}, \; \bm{a}_h=(a_{x,h},a_{y,h},a_{{\rm r},h}, a_{{\rm t},h})^{\rm T}.
\end{equation*}
We thus compute
\begin{equation}\label{dtahR4}
\abs{(\mc{D}_t\bm{a}_h)^{n+1/2}}.
\end{equation}
The norms \eqref{PiXC1h} and \eqref{dtahR4} are plotted in Figure \ref{fig:decay_rate}.
It is clearly seen that the asymptotic decay rate conforms to the theory.
The decay with initial data \eqref{unlabeled_curve} is almost precisely exponential at the theoretical rate, 
which is to be expected given that the perturbation is proportional to the primary decay modes in \eqref{primary_decay_mode}.
For the other initial data, the decay rate asymptotically approaches the theoretically predicted value.

\section{Global Behavior}\label{sect:global}

Recall that the area and energy identities \eqref{area_conservation} and \eqref{energy_identity} are satisfied 
for mild solutions, as proved in Proposition \ref{equiv_formulations}.
Take area conservation. Viewing the interior fluid area $\abs{\Omega_{\rm i}}$ as a function of $t$, we have:
\begin{equation*}
\abs{\Omega_{\rm i}}(t)=\abs{\Omega}_{\rm i}(\epsilon) \text{ for } 0<\epsilon<t
\end{equation*}
so long as the solution exists up to time $t$. Given the expression for $\abs{\Omega_{\rm i}}$ given in 
\eqref{area_conservation}, it is clear that the $\abs{\Omega_{\rm i}}$ is a continuous functional of 
$\bm{X}\in C^{1,\gamma}(\mbs)$. Since our mild solution $\bm{X}(t)\in C([0,T]; C^{1,\gamma}(\mbs))$, 
we may take the limit as $\epsilon\to 0$ in the above equation to find that:
\begin{equation*}
\abs{\Omega_{\rm i}}(t)=\abs{\Omega_{\rm i}}(0).
\end{equation*}

Let us now turn to energy conservation. From \eqref{energy_identity}, we have:
\begin{equation*}
\mc{E}(t)-\mc{E}(\epsilon)=-\int_\epsilon^t \mc{D}(s)ds \text{ for } 0<\epsilon<t.
\end{equation*}
The energy functional $\mc{E}$ is continuous with respect to $\bm{X}\in C^{1,\gamma}(\mbs)$, 
and therefore, we may take the limit as $\epsilon\to 0$ in the above to find:
\begin{equation}\label{EtE0Ds}
\mc{E}(t)-\mc{E}(0)=-\int_0^t \mc{D}(s)ds.
\end{equation}
so long as the solution exists up to time $t$. Note that the right hand side does not need to 
be interpreted as an improper integral since the integrand is non-negative.

We now state a simple observation that is a consequence of the above.
\begin{lemma}\label{apriorinormbounds}
Suppose we have a mild solution $\bm{X}(t)\in C([0,T]; C^{1,\gamma}(\mbs)), 0<\gamma<1$. 
Then, $\starnorm{\bm{X}}$ has an upper bound and $\norm{\bm{X}}_{\dot{C}^{1,\gamma}}$ has a lower bound.
More concretely, 
\begin{align}
\label{starnormupperbound}
\starnorm{\bm{X}(t)}&\leq \sqrt{\frac{\mc{E}(0)}{\pi}},\\
\label{C1gammanormlowerbound}
\norm{\p_\theta \bm{X}(t)}_{C^{0,\gamma}}&\geq \sqrt{\frac{\abs{\Omega_{\rm i}(0)}}{\pi}}. 
\end{align}
\end{lemma}
\begin{poof}
We first consider the \eqref{starnormupperbound}. From \eqref{EtE0Ds}, $\mc{E}(t)\leq \mc{E}(0)$, and therefore,
\begin{equation*}
\mc{E}(0)\geq \mc{E}(t)=\frac{1}{2}\int_{\mbs}\abs{\p_\theta\bm{X}}^2d\theta
\geq \pi\inf_{\theta\in \mbs} \abs{\p_\theta \bm{X}}^2\geq \pi\starnorm{\bm{X}}^2.
\end{equation*}
For \eqref{C1gammanormlowerbound}, we use the isoperimetric inequality:
\begin{equation*}
\sqrt{4\pi |\Omega_{\rm i}(0) |}\leq \int_{\mbs}\abs{\p_\theta \bm{X}}d\theta
\leq 2\pi \sup_{\theta\in \mbs}\abs{\p_\theta \bm{X}}\leq 2\pi \norm{\p_\theta\bm{X}}_{C^{0,\gamma}}.
\end{equation*}
\end{poof}

We are now ready to prove Theorem \ref{GlobalBehaviorTheorem}.

\begin{poof}[Proof of Theorem \ref{GlobalBehaviorTheorem}]
We first make the following observation. 
Suppose we have:
\begin{equation}\label{defratiobound}
\sup_{0\leq t<\tau_{\rm max}(\bm{X}_0)} \varrho_\alpha(\bm{X}(t))=K<\infty.
\end{equation}
By Lemma \ref{apriorinormbounds}, we immediately have:
\begin{equation}\label{apriori+ratiobound}
\sqrt{\frac{\abs{\Omega_{\rm i}(0)}}{\pi}}\leq \norm{\p_\theta \bm{X}}_{C^{0,\alpha}}\leq K\starnorm{\bm{X}}\leq K\sqrt{\frac{\mc{E}(0)}{\pi}}.
\end{equation}
A bound on the ratio thus immediately results in upper and lower bounds for $\norm{\p_\theta\bm{X}}_{C^{0,\alpha}}$ and $\starnorm{\bm{X}}$.

We first consider item \ref{GBT1}.
Suppose that, for some $\alpha>0$, 
\begin{equation*}
\limsup_{t\to \tau_{\max}(\bm{X}_0)} \varrho_\alpha(\bm{X}(t))<\infty.
\end{equation*}
We then show we can then extend the mild solution beyond $\tau_{\rm max}(\bm{X}_0)$.
Since $\varrho_\alpha>\varrho_\alpha'$ for any $\alpha>\alpha'$, we may as well assume that $0<\alpha\leq \gamma$ and $\alpha<1/2$.
Then, \eqref{defratiobound} holds, and thus we have \eqref{apriori+ratiobound}. In particular, we have:
\begin{equation}\label{starboundtotaumax}
m=\frac{1}{K}\sqrt{\frac{\abs{\Omega_{\rm i}(0)}}{\pi}}\leq \starnorm{\bm{X}} \text{ for } 0\leq t<\tau_{\rm max}(\bm{X}_0).
\end{equation}
By \eqref{FCCdot} and \eqref{FTCdot}, this implies that there is a constant $M_{\mc{R}}$ such that:
\begin{equation*}
\chnorm{\mc{R}(\bm{X}(t))}{0,2\alpha}\leq C\frac{\chnorm{\p_\theta \bm{X}(t)}{0,\alpha}^4}{\starnorm{\bm{X}(t)}^3}\leq M_{\mc{R}} \text{ for } 0\leq t<\tau_{\rm max}(\bm{X}_0).
\end{equation*}
Thus,
\begin{equation*}
\begin{split}
\chnorm{\bm{X}(t)}{1,\alpha}
&\leq C\chnorm{\bm{X}_0}{1,\alpha}+\int_0^t \norm{e^{(t-s)\Lambda}\mc{R}(\bm{X}(s))}ds\\
&\leq C\chnorm{\bm{X}_0}{1,\alpha}+\int_0^t \frac{C}{(t-s)^{1-\alpha}}M_{\mc{R}}ds\leq M \text{ for } 0\leq t<\tau_{\rm max}(\bm{X}_0). 
\end{split}
\end{equation*}
From this, we know from Lemma \ref{Cnboundlemma} that, for each $0<\beta<1$ there are constants $M_\beta>0$ such that 
\begin{equation}\label{C2betabound}
\chnorm{\bm{X}(t)}{2,\beta}\leq M_\beta\;\text{ for all }\; \frac{1}{2}\tau_{\rm max}(\bm{X}_0) \leq t<\tau_{\rm max}(\bm{X}_0).
\end{equation}
From the strong form of our equation and Lemma \ref{l:CarryRegularity}
we see that, for some constant $\wt{M}_\beta$,
\begin{equation*}
\norm{\p_t\bm{X}}_{C^{1,\beta}}\leq \wt{M}_\beta\; \text{ for all }\; \frac{1}{2}\tau_{\rm max}(\bm{X}_0) \leq t<\tau_{\rm max}(\bm{X}_0).
\end{equation*}
Applying this to $\beta=\alpha$ in particular, we see that $\bm{X}(t)$ is uniformly bounded in $C^1([\tau_{\rm max}(\bm{X}_0)/2,\tau_{\rm max}(\bm{X}_0)),C^{1,\alpha}(\mbs))$.
This implies that $\bm{X}(t)\in C([0,\tau_{\rm max}(\bm{X}_0));C^{1,\alpha}(\mbs))$ is uniformly continuous in time. 
Thus, the following limit exists in $C^{1,\alpha}(\mbs)$:
\begin{equation*}
\lim_{t\to \tau_{\rm max}(\bm{X}_0)} \bm{X}(t)=\bm{X}_\star.
\end{equation*}
By \eqref{C2betabound}, $\bm{X}_\star\in C^{2,\beta}(\mbs)\subset h^{1,\alpha}(\mbs)$, and by \eqref{starboundtotaumax}, we have $\starnorm{\bm{X}_\star}\geq m>0$.
This means that we may continue the solution on from $\tau_{\rm max}(\bm{X}_0)$ using our local existence theorem Theorem \ref{LWPTheorem}.
A mild solution in $C^{1,\alpha}(\mbs)$ is a mild solution in $C^{1,\gamma}(\mbs)$ by our regularity results established in Theorem \ref{SmoothTheorem}.

Let us next consider item \ref{GBT2}. Our method is to consider the $\omega$-limit set of the global solution 
with bounded deformation ratio, although we will not 
explicitly use the terminology associated with $\omega$-limit sets (see \cite{sell2013dynamics} for example)
since our setting is quite simple.
Define:
\begin{equation*}
\wh{\bm{X}}=\Pi_{\rm trl}\bm{X},
\end{equation*}
where $\Pi_{\rm trl}$ is the projection operator defined in \eqref{proj_trl}.
Clearly, 
\begin{equation*}
\norm{\p_\theta\wh{\bm{X}}}_{C^{0,\beta}}=\norm{\p_\theta \bm{X}}_{C^{0,\beta}}, \; 0<\beta<1.
\end{equation*}
Furthermore, it is easily seen that the following Poincar\'e type inequality holds 
\begin{equation}\label{whXequivnorm}
\norm{\wh{\bm{X}}}_{C^{1,\beta}}\leq C\norm{\p_\theta \wh{\bm{X}}}_{C^{0,\beta}},
\end{equation}
for some constant $C$ independent of $\bm{X}$. Indeed, let $\bm{\wh{X}}=(\wh{X},\wh{Y})$.
By the definition of $\Pi_{\rm trl}$, we have: 
\begin{equation*}
\int_{\mbs} \wh{X}d\theta=0.
\end{equation*} 
There must thus be a point $\theta_*\in \mbs$ at which $\wh{X}(\theta_*)=0$.
Thus, 
\begin{equation*}
\abs{\wh{X}(\theta)}\leq \int_{\theta_*}^\theta \abs{\p_\theta \wh{X}}d\theta\leq 2\pi \norm{\p_\theta\wh{X}}_{C^0}.
\end{equation*}
A similar bound can be found for $\wh{Y}$.

Suppose bound \eqref{defratiobound} is satisfied with $\tau_{\rm max}=\infty$. 
By \eqref{apriori+ratiobound} and \eqref{whXequivnorm}, this implies that:
\begin{equation}\label{C1alphainfbound}
\sup_{t\geq 0} \norm{\wh{\bm{X}}}_{C^{1,\alpha}}\equiv M_\alpha<\infty.
\end{equation}
Take any sequence $t_1<t_2<t_3<\cdots t_k\to \infty$ and let $\bm{X}_k=\bm{X}(t_k)$. The set consisting of 
$\bm{\wh{X}}_k=\Pi_{\rm trl}\bm{X}_k$ is precompact in $C^{1,\beta}(\mbs), 0<\beta<\alpha$, thanks to \eqref{C1alphainfbound}.
We may thus extract a subseqence of the time points above, which we shall continue 
to call $t=t_k$, so that $\wh{\bm{X}}_k\to \wh{\bm{X}}_\star$ in $C^{1,\beta}(\mbs)$ for some $\wh{\bm{X}}_\star\in C^{1,\beta}(\mbs)$.
In fact, $\wh{\bm{X}}_\star\in h^{1,\beta}(\mbs)$ since it is in the completion of $C^{1,\alpha}(\mbs)$ in $C^{1,\beta}(\mbs)$.

Consider the mild solutions $\bm{W}_k(t)$ and $\bm{W}_\star(t)$ with initial data $\wh{\bm{X}}_k$ and $\wh{\bm{X}}_\star$ respectively.
Note that there is a local mild solution with initial data $\wh{\bm{X}}_\star$ thanks to \eqref{apriori+ratiobound}; $\starnorm{\wh{\bm{X}}_\star}$ is bounded from below.
By continuity with respect to initial data established in Theorem \ref{LWPTheorem}, 
there exists a $T>0$ such that, for $k$ sufficiently large, $\bm{W}_k(t)$ is well-defined for $0\leq t\leq T$, and
\begin{equation}\label{WkWstar}
\bm{W}_k(t)\to \bm{W}_\star(t) \text{ in } C([0,T]; C^{1,\beta}(\mbs)).
\end{equation}

We now argue that the energy $\mc{E}$ is constant on the solution $\bm{W}_\star(t)$.
Let us use the notation $\mc{E}(\bm{Z})$ to mean the energy evaluated at configuration $\bm{Z}$.
Consider the original mild solution $\bm{X}(t)$. The energy is monotone decreasing and is non-negative, and therefore, $\mc{E}(\bm{X}(t))$
converges to some value $\mc{E}_\star$. Since $\mc{E}(\bm{X}_k)=\mc{E}(\wh{\bm{X}}_k)$, we have $\mc{E}(\wh{\bm{X}}_\star)=\mc{E}_\star$
(the energy is clearly continuous with respect to the $C^{1,\beta}$ norm).
The same argument can be made for $\mc{E}(\bm{W}_\star(t)), \; 0\leq t\leq \tau$. Indeed,
\begin{equation*}
\mc{E}(\bm{W}_k(t))=\mc{E}(\bm{X}(t+t_k))\to \mc{E}_\star \text{ as } k\to \infty.
\end{equation*}
Thus, by \eqref{WkWstar},
\begin{equation*}
\mc{E}(\bm{W}_\star(t))=\mc{E}_\star, \; 0\leq t\leq \tau.
\end{equation*}
From \eqref{energy_identity}, this implies that the dissipation $\mc{D}$ is $0$ along $\bm{W}_\star(t)$, which in turn implies that the velocity field $\bm{u}$
is identically $0$. 
By the arguments leading to Proposition \ref{equilibria_are_circles}, this shows that 
that $\bm{W}_\star(t)=\wh{\bm{X}}_\star$ is a uniformly parametrized stationary circle.
We have thus shown that:
\begin{equation*}
\norm{\wh{\bm{X}}_k-\wh{\bm{X}}_\star}_{C^{1,\beta}}=\norm{\bm{X}_k-\bm{Z}_k}_{C^{1,\beta}} \to 0 \text{ as } k\to \infty, \; 
\bm{Z}_k=\wh{\bm{X}}_\star+\mc{P}_{\rm trl}\bm{X}_k,
\end{equation*}
where $\mc{P}_{\rm trl}$ was defined in \eqref{proj_trl}. Note that the configurations $\bm{Z}_k$ are just the circle $\bm{X}_\star$
translated by $\mc{P}_{\rm trl}\bm{X}_k$. 
Now, note that:
\begin{equation*}
\norm{\Pi \bm{X}_k}_{C^{1,\beta}}=\norm{\bm{X}_k-\mc{P}\bm{X}_k}_{C^{1,\beta}}\leq  
\norm{\bm{X}_k-\bm{Z}_k}_{C^{1,\beta}}+\norm{\mc{P}(\bm{X}_k-\bm{Z}_k)}_{C^{1,\beta}}\leq C\norm{\bm{X}_k-\bm{Z}_k}_{C^{1,\beta}},
\end{equation*}
where we used $\mc{P}\bm{Z}_k=\bm{Z}_k$ in the second inequality and the boundedness of $\mc{P}$ in the $C^{1,\beta}$ norm in the 
last inequality. Thus, 
\begin{equation*}
\lim_{k\to \infty}\norm{\Pi\bm{X}(t_k)}=0.
\end{equation*}
Take $k$ sufficiently large, so that $\norm{\Pi\bm{X}(t_k)}_{C^{1,\beta}}$ is small enough to apply Theorem \ref{StabilityTheorem} with initial 
data $\bm{X}(t_k)$. This concludes the proof.
\end{poof}

\bibliographystyle{plain}
\bibliography{paper}{}

\end{document}